\documentclass[12pt, reqno,oneside]{amsart}
\allowdisplaybreaks
\usepackage{amssymb, amsmath, amsthm, wasysym, mathrsfs}
\usepackage[backref]{hyperref}
\usepackage[alphabetic,backrefs,lite]{amsrefs}
\usepackage{fullpage}
\usepackage{setspace}
\usepackage{color}
\usepackage{mathtools}
\usepackage{enumitem}

\DeclareMathOperator{\Br}{Br}

\DeclareMathOperator{\Div}{Div}
\DeclareMathOperator{\DIv}{div}
\DeclareMathOperator{\Eff}{Eff}

\DeclareMathOperator{\Gal}{Gal}

\DeclareMathOperator{\PDiv}{PDiv}
\DeclareMathOperator{\Pic}{Pic}
\DeclareMathOperator{\pr}{pr}

\DeclareMathOperator{\Res}{Res}
\DeclareMathOperator{\Spec}{Spec}

\DeclareMathOperator{\tr}{tr}

\theoremstyle{definition}
\newtheorem{dfn}{Definition}[section]

\newtheorem{lem}[dfn]{Lemma}
\newtheorem{prop}[dfn]{Proposition}
\newtheorem{co}[dfn]{Corollary}
\newtheorem{rem}[dfn]{Remark}

\newtheorem{nota}[dfn]{Notation}
\newtheorem{con}[dfn]{Condition}
\theoremstyle{plain}
\newtheorem{theorem}[dfn]{Theorem}
\begin{document}
\title{Semi-integral points of bounded height on vector group compactifications}
\date{}
\author{Haruki Ito}
\address{Graduate School of Mathematics, Nagoya University, Furocho, Chikusaku, Nagoya 464-8602, Japan}
\email{haruki.ito.c8@math.nagoya-u.ac.jp}

\begin{abstract}
In this paper, we obtain the asymptotic formula for the counting function of Darmon points of bounded height on equivariant compactifications of vector groups using ideas similar to those in \cite{PSTVA21}.
We also show that this asymptotic formula is compatible with Manin's conjecture for $\mathcal{M}$-points formulated in \cite{Moe26a}.
\end{abstract}

\maketitle

\tableofcontents

\section{Introduction}

\subsection{Main results}
In arithmetic geometry, the problem of counting rational points of bounded height is central.
Our main result is a Darmon analogue of the asymptotic formula for Campana points on equivariant compactifications of vector groups proved in \cite{PSTVA21}.
Roughly speaking, Darmon points are rational points satisfying divisibility conditions on intersection multiplicities with boundary divisors; see Definition~\ref{Darmon} for the precise definition.

Let $F$ be a number field, let $G=\mathbb{G}_{a,F}^n$ be the $n$-dimensional vector group over $F$, and let $X$ be a smooth projective equivariant compactification of $G$.
Let $L$ be a big line bundle on $X$ equipped with a smooth adelic metric.
As explained in \S~3.3, $L$ induces a height function
\[
H_L\colon G(F)\to\mathbb{R}_{>0}.
\]

Let $X\setminus G=\bigcup_{\alpha\in\mathcal{A}}D_\alpha$ be the decomposition into irreducible components.
By \cite[PROPOSITION~5.1]{PSTVA21} and the bigness of $L$, there exists a tuple of positive integers $(\lambda_\alpha)_{\alpha\in\mathcal{A}}$ such that $L=\sum_{\alpha\in\mathcal{A}}\lambda_\alpha D_\alpha$.
For each $\mathbf{m}=(m_\alpha)_{\alpha\in\mathcal{A}}\in\mathbb{Z}_{>0}^{\mathcal{A}}$, we set
\[
D_{\mathbf{m}}:=\sum_{\alpha\in\mathcal{A}}\left(1-\frac{1}{m_\alpha}\right)D_\alpha.
\]
Let $S$ be a finite set of places of $F$ containing all infinite places, let $\mathcal{O}_{F,S}$ be the ring of $S$-integers of $F$, and let $(\mathcal{X},\mathcal{D}_{\mathbf{m}})$ be a good integral model of $(X,D_{\mathbf{m}})$ away from $S$.
For each $\alpha\in\mathcal{A}$, we denote the Zariski closure of $D_\alpha$ in $\mathcal{X}$ by $\mathcal{D}_\alpha$.
We say that a rational point $P\in G(F)$ is a Darmon $\mathcal{O}_{F,S}$-point on $(\mathcal{X},\mathcal{D}_{\mathbf{m}})$ if for every $\alpha\in\mathcal{A}$ and every finite place $v\in\Omega_F^{<\infty}\setminus S$, the intersection multiplicity of $P$ with $\mathcal{D}_\alpha$ at $v$ is a multiple of $m_\alpha$.
We denote the set of Darmon $\mathcal{O}_{F,S}$-points on $(\mathcal{X},\mathcal{D}_{\mathbf{m}})$ by $(\mathcal{X},\mathcal{D}_{\mathbf{m}})^{\mathrm{D}}(\mathcal{O}_{F,S})$.
Let $(\rho_\alpha)_{\alpha\in\mathcal{A}}$ be a tuple of integers such that $-K_X=\sum_{\alpha\in\mathcal{A}}\rho_\alpha D_\alpha$, where $K_X$ is the canonical divisor class of $X$.
We prove the following asymptotic formula for the counting function of Darmon $\mathcal{O}_{F,S}$-points on $(\mathcal{X},\mathcal{D}_{\mathbf{m}})$, which is an analogue of \cite[THEOREM~1.2]{PSTVA21}.

\begin{theorem}\label{main}
Keep the above notation.
Define constants $a$ and $b$ by
\begin{align*}
a&=\max_{\alpha\in\mathcal{A}}\frac{\rho_\alpha-1+m_\alpha^{-1}}{\lambda_\alpha},
\\
b&=\#\left\{\alpha\in\mathcal{A}\,\middle|\,a=\frac{\rho_\alpha-1+m_\alpha^{-1}}{\lambda_\alpha}\right\}.
\end{align*}
Assume that the $\mathbb{Q}$-divisor $aL+K_X+D_{\mathbf{m}}$ on $X$ is rigid.
Then there exists a positive constant $c$ such that the asymptotic formula
\[
\#\{P\in(\mathcal{X},\mathcal{D}_{\mathbf{m}})^{\mathrm{D}}(\mathcal{O}_{F,S})\mid H_L(P)\leq B\}\sim\frac{c}{a(b-1)!}B^a(\log B)^{b-1}\quad(B\to\infty)
\]
holds.
Moreover, the constant $c$ is given in Proposition~\ref{c} by
\[
c=\lim_{s\to a}(s-a)^b\widehat{\mathsf{H}}_{\mathbf{m}}(\mathbf{0},sL).
\]
Here $\widehat{\mathsf{H}}_{\mathbf{m}}(\mathbf{0},sL)$ is the height integral at $\mathbf{a}=\mathbf{0}$ defined in Definition~\ref{HI}~(2).
\end{theorem}

We note that the constants $a$ and $b$ coincide with those defined in \cite[\S~9]{PSTVA21}.
On the other hand, the leading constant $c$ differs in general from the corresponding constant in \cite{PSTVA21}.
Furthermore, in the last section, we show that this asymptotic formula is compatible with Manin's conjecture for $\mathcal{M}$-points formulated in \cite[Conjecture~1.4]{Moe26a}.

\subsection{Background}
The problem of counting rational points of bounded height has been studied extensively in the context of Manin's conjecture.
This conjecture concerns the distribution of rational points on smooth projective Fano varieties over number fields.
It was proposed by Y.~Manin and his collaborators in the late 1980s and has been refined in \cite{FMT89, BM90, Pey95, BT98a, Pey03, Pey17, LST22}.
Manin-type counting problems have also been studied for integral points; see, for example, \cite{BO12, CLT12, CLT12b, TBT13, TT15, Cho19, Wil24, San25}.
Recently, the distribution of semi-integral points, such as Campana points, Darmon points, and $\mathcal{M}$-points, has attracted increasing attention.

Let $X$ be a projective variety over a number field $F$, and let $(D_\alpha)_{\alpha\in\mathcal{A}}$ be a finite family of distinct prime divisors on $X$.
Fix a finite set $S$ of places of $F$ containing all infinite places, and let $\mathcal{X}$ be an integral model of $X$ over $\mathcal{O}_{F,S}$.
For each $\alpha\in\mathcal{A}$, we denote by $\mathcal{D}_\alpha$ the
Zariski closure of $D_\alpha$ in $\mathcal{X}$.
Let
\[
\mathcal{D}_{\mathbf{m}}=\sum_{\alpha\in\mathcal{A}}\left(1-\frac{1}{m_\alpha}\right)\mathcal{D}_\alpha
\]
be a weighted $\mathbb{Q}$-divisor on $\mathcal{X}$, where each $m_\alpha$ is a positive integer.

Campana points on $(\mathcal{X},\mathcal{D}_{\mathbf{m}})$ are rational points that are integral with respect to $\mathcal{D}_{\mathbf{m}}$ in the sense of Campana.
More precisely, we say that $P\in (X\setminus\bigcup_{\alpha\in\mathcal{A}}D_\alpha)(F)$ is a Campana point on $(\mathcal{X},\mathcal{D}_{\mathbf{m}})$ if, for all $\alpha\in\mathcal{A}$ and all finite places $v\notin S$, the intersection multiplicity of $P$ with $\mathcal{D}_\alpha$ at $v$ is either zero or at least $m_\alpha$.
Campana points were defined in \cite[Definition~4.1]{Cam05} and \cite[Definition~2.4.17]{Abr09}, and the distribution of Campana points has been studied in \cite{BVV12, VV12, BY21, PSTVA21, Shu22, Str22, Xia22, BBK24, PS24, CLTBT26}.
Manin's conjecture for Campana points is formulated in
\cite[CONJECTURE~1.1]{PSTVA21} and \cite[Conjecture~8.3]{CLTBT26}.

We say that a rational point $P\in (X\setminus\bigcup_{\alpha\in\mathcal{A}}D_\alpha)(F)$ is a Darmon $\mathcal{O}_{F,S}$-point on $(\mathcal{X},\mathcal{D}_{\mathbf{m}})$ if, for all $\alpha\in\mathcal{A}$ and all finite places $v\notin S$, the intersection multiplicity of $P$ with $\mathcal{D}_\alpha$ at $v$ is a multiple of $m_\alpha$.
Darmon points were introduced in \cite[Definition~2.10]{MNS24}, based on ideas arising from H.~Darmon's work \cite{Dar97} on generalized Fermat equations.
We note that related notions were considered in \cite[Definition~2.4]{Cam11} and \cite[\S~7.6]{Cam15}.
In \cite[Theorem~1.2]{AP25}, the distribution of solutions to generalized Fermat equations was studied.
As explained in \cite[\S~9.3]{Moe26a}, this can be interpreted as a result on the distribution of Darmon points on projective space.
Moreover, \cite[Theorem~1.1]{SS24} studies the distribution of geometric Darmon points on toric varieties, which may be regarded as geometric analogues of Darmon points.
In \cite[\S~8]{Moe24}, Darmon points are related to integral points on root stacks constructed as in \cite[Definition~2.2.4]{Cad07}.

The notion of $\mathcal{M}$-points, which generalizes these notions, was introduced in \cite[Definition~3.12]{Moe24} and \cite[Definition~3.12]{Moe26a}. Manin's conjecture for $\mathcal{M}$-points was proposed in \cite[Conjecture~1.4]{Moe26a}.
In \cite[Theorem~1.1]{Moe26b}, B. Moerman proved Manin's conjecture for $\mathcal{M}$-points on smooth proper split toric varieties over $\mathbb{Q}$ under the assumption that the pair $(X,M)$ is smooth and proper.
As explained in \S~5, the invariants $a$ and $b$ appearing in Theorem~\ref{main} are compatible with the corresponding invariants in \cite[Definition~4.37]{Moe26a}.
Moreover, the leading constant agrees with the leading constant defined in \cite[Definition~8.5]{Moe26a}.

Asymptotic formulae for rational points, integral points, and Campana points on equivariant compactifications of vector groups were established in \cite{CLT02, CLT12, PSTVA21}, respectively.
In this paper, we study the corresponding problem for Darmon points in this setting.
For a fixed weighted divisor, the invariants $a$ and $b$ in Theorem~\ref{main} coincide with the corresponding invariants in \cite{PSTVA21}, whereas the leading constant differs in general.

\subsection{Strategy of the proof}
We keep the notation introduced in \S~1.1.
The proof of the main theorem follows the strategy of \cite{PSTVA21}.
By the Tauberian theorem proved in \cite{Del54}, the asymptotic formula in Theorem~\ref{main} follows from the analytic properties of the series
\[
\mathsf{Z}_{\mathbf{m}}(sL):=\sum_{\mathbf{x}\in(\mathcal{X},\mathcal{D}_{\mathbf{m}})^{\mathrm{D}}(\mathcal{O}_{F,S})}\frac{1}{H_L(\mathbf{x})^s}.
\]
By an argument similar to that in \cite[Proposition~4.5]{CLT02}, the above series converges absolutely for $\Re s$ sufficiently large.
For $\Re s$ sufficiently large, the Poisson summation formula gives
\[
\mathsf{Z}_{\mathbf{m}}(sL)=\sum_{\mathbf{a}\in G(F)}\widehat{\mathsf{H}}_{\mathbf{m}}(\mathbf{a},sL)=\widehat{\mathsf{H}}_{\mathbf{m}}(\mathbf{0},sL)+\sum_{\mathbf{a}\in G(F)\setminus\{\mathbf{0}\}}\widehat{\mathsf{H}}_{\mathbf{m}}(\mathbf{a},sL),
\]
where, for each $\mathbf{a}\in G(F)$,
\[
\widehat{\mathsf{H}}_{\mathbf{m}}(\mathbf{a},sL)=\int_{G(\mathbb{A}_F)}\delta_{\mathbf{m}}(\mathbf{x})\mathsf{H}(\mathbf{x},-sL)\psi_{\mathbf{a}}(\mathbf{x})\,\mathrm{d}\mathbf{x}
\]
is the Fourier transform at $\mathbf{a}$ of the function $\mathbf{x}\mapsto\delta_{\mathbf{m}}(\mathbf{x})\mathsf{H}(\mathbf{x},-sL)$.
In Proposition~\ref{c}, we show that the function
\[
s\mapsto\widehat{\mathsf{H}}_{\mathbf{m}}(\mathbf{0},sL)
\]
has a rightmost pole at $s=a$ of order $b$, and the constant
\[
c=\lim_{s\to a}(s-a)^b\widehat{\mathsf{H}}_{\mathbf{m}}(\mathbf{0},sL)
\]
is positive.
On the other hand, in Proposition~\ref{Hzero} and Proposition~\ref{Haneq0}, we show that, under the assumption that the $\mathbb{Q}$-divisor $aL+K_X+D_{\mathbf{m}}$ on $X$ is rigid, the function
\[
s\mapsto\sum_{\mathbf{a}\in G(F)\setminus\{\mathbf{0}\}}\widehat{\mathsf{H}}_{\mathbf{m}}(\mathbf{a},sL)
\]
does not contribute to the rightmost pole of $\mathsf{Z}_{\mathbf{m}}(sL)$.
Therefore, the function $s\mapsto\mathsf{Z}_{\mathbf{m}}(sL)$ has a rightmost pole at $s=a$ of order $b$, and we have
\[
\lim_{s\to a}(s-a)^b\mathsf{Z}_{\mathbf{m}}(sL)=c>0.
\]
Applying the Tauberian theorem then gives Theorem~\ref{main}.
\subsection{Structure of the paper}

In \S~2, we set up the notation.
In \S~3, we recall equivariant compactifications of vector groups, Darmon points, and height functions, in order to state the main theorem precisely.
In \S~4, we prove the main theorem.
We first study the analytic properties of the height integral in \S~4.1.
In \S~4.2, we recall reduction maps, which are used to compute local integrals at good places.
In \S~4.3, following \cite[\S~7]{PSTVA21}, we compute the height integral $\widehat{\mathsf{H}}_{\mathbf{m}}(\mathbf{a},\mathbf{s})$ in the case $\mathbf{a}=\mathbf{0}$, and in \S~4.4 we treat the case $\mathbf{a}\neq\mathbf{0}$.
In \S~5, we show that Theorem~\ref{main} is compatible with Manin's conjecture for $\mathcal{M}$-points formulated in \cite[Conjecture~1.4]{Moe26a}.
\subsection{Acknowledgements}

The author would like to thank his advisor Sho Tanimoto for continuous support and encouragement.
The author is grateful to Fumiya Okamura, Shuhei Katsuta, Rikuto Ito, Marta Pieropan, Sam Streeter, and Boaz Moerman for their helpful comments and suggestions, which improved the quality of this manuscript.
This paper is based on the author's master's thesis \cite{Ito25}. 
This work was financially supported by JST SPRING, Grant Number JPMJSP2125, through the THERS Make New Standards Program for the Next Generation Researchers.

\section{Notation}
\begin{nota}
For any number field $K$, we use the following notation.

Let $\mathcal{O}_K$ be the ring of integers of $K$, let $\mathbb{A}_K$ be the adele ring of $K$, let $\Omega_K$ be the set of places of $K$, let $\Omega_K^\infty$ be the set of infinite places of $K$, and let $\Omega_K^{<\infty}$ be the set of finite places of $K$.
For each finite subset $S$ of $\Omega_K$ containing $\Omega_K^\infty$, let $\mathcal{O}_{K,S}$ be the ring of $S$-integers of $K$.
We denote by $K_v$ the completion of $K$ at $v$ for each $v\in\Omega_K$.
For each $v\in\Omega_K^{<\infty}$, let $\mathcal{O}_v$ be the valuation ring of $K_v$, let $\mathfrak{m}_v$ be the maximal ideal of $\mathcal{O}_v$, let $k_v=\mathcal{O}_v/\mathfrak{m}_v$, let $q_v=\#k_v$, and let $\pi_v$ be a prime element of $\mathcal{O}_v$.
For each $v\in\Omega_K^{<\infty}$, let $p_v$ be the unique prime number such that the restriction of the $v$-adic topology to $\mathbb{Q}$ is equivalent to the $p_v$-adic topology, and let $\mathfrak{D}_v$ be the different of the field extension $K_v/\mathbb{Q}_{p_v}$ with norm $N\mathfrak{D}_v$.
As in \cite[\S~2.1]{Tat67}, for each $v\in\Omega_K$, we fix a measure $\mathrm{d}x_v$ on $K_v$, also denoted by $\mu_v$, as follows:
\begin{enumerate}
\item if $v$ is a real place, $\mathrm{d}x_v$ is the ordinary Lebesgue measure,
\item if $v$ is a complex place, $\mathrm{d}x_v$ is twice the ordinary Lebesgue measure,
\item if $v$ is a finite place, $\mathrm{d}x_v$ is normalized so that $\mathcal{O}_v$ has volume $(N\mathfrak{D}_v)^{-1/2}$.
\end{enumerate}
We endow $\mathbb{A}_K$ with the restricted product topology and the Haar measure induced by the measures $\mathrm{d}x_v$ on $K_v$; we endow $K_v^n$ with the product topology and the Haar measure induced by the measures $\mathrm{d}x_v$ on $K_v$, and we endow $\mathbb{A}_K^n$ with the product topology and the Haar measure induced by the measures on $\mathbb{A}_K$.
For each $v\in\Omega_K$ and each $x\in K_v$, we define $|x|_v$ to be the unique non-negative real number satisfying
\[
\mu_v(xB)=|x|_v\mu_v(B)
\]
for all measurable subsets $B$ of $K_v$.
In particular, for each finite place $v$, we have $|\pi_v|_v=q_v^{-1}$.

Let $F$ be a number field and let $G$ be the $n$-dimensional vector group over $F$.
For each $\mathbf{x}\in G(\mathbb{A}_F)$ and each $v\in\Omega_F$, $\mathbf{x}_v$ denotes the image of $\mathbf{x}$ under the projection $G(\mathbb{A}_F)\to G(F_v)$.
We endow $G(F_v)$ and $G(\mathbb{A}_F)$ with the topology and Haar measure induced by the isomorphisms $G(F_v)\cong F_v^n$ and $G(\mathbb{A}_F)\cong\mathbb{A}_F^n$, respectively.
\end{nota}

\begin{nota}
\leavevmode
\begin{enumerate}
\item For each finite extension $K'/K$, $\tr_{K'/K}\colon K'\to K$ denotes the trace map of $K'/K$.
\item For a finite set $\mathcal{A}$, we write tuples indexed by $\mathcal{A}$ in boldface. Thus, for
\(\mathbf{s}\in\mathbb{C}^{\mathcal{A}}\), we write
\[
\mathbf{s}=(s_\alpha)_{\alpha\in\mathcal{A}},
\]
where the components are denoted by the same letter with indices.
We set $\Re\mathbf{s}=\min_{\alpha\in\mathcal{A}}\Re s_\alpha$ and $|\mathbf{s}|=\max_{\alpha\in\mathcal{A}}|s_\alpha|$ for each $\mathbf{s}=(s_\alpha)_{\alpha\in\mathcal{A}}\in\mathbb{C}^{\mathcal{A}}$.
\item For a smooth projective rationally connected variety $X$ over $F$, we denote the Picard group of $X$, the effective cone of $X$, and the pseudo-effective cone of $X$ by $\Pic(X)$, $\Eff^1(X)$, and $\overline{\Eff}^1(X)$, respectively.
We identify line bundles and divisors on $X$ with their corresponding classes in $\Pic(X)$.
\item For a set $A$, we denote its indicator function by $\mathbf{1}_A$.
\item We set $\overline{\mathbb{N}}:=\mathbb{Z}_{\geq0}\cup\{\infty\}$.
\end{enumerate}
\end{nota}
\section{Preliminaries}

\subsection{Equivariant compactifications of vector groups}

In this subsection, we recall equivariant compactifications of vector groups studied in \cite{HT99} and \cite{CLT02}.
Manin's conjecture for rational points, integral points, and Campana points on equivariant compactifications of vector groups was proved in \cite{CLT02}, \cite{CLT12}, and \cite{PSTVA21}, respectively.
In this framework, we study the distribution of Darmon points on equivariant compactifications of vector groups.

We use the following notation throughout this paper. 
Let $F$ be a number field, let $G=\mathbb{G}_{a,F}^n$ be the $n$-dimensional vector group over $F$, and let $X$ be a smooth projective equivariant compactification of $G$ over $F$; that is, $X$ is equipped with an action of $G$ and contains $G$ as an open dense orbit.
Let $X\setminus G=\bigcup_{\alpha\in\mathcal{A}}D_\alpha$ be the decomposition into irreducible components.

We recall some results on the Picard group, the effective cone, and the anticanonical divisor of $X$.

\begin{rem}\label{Eff}
\leavevmode
\begin{enumerate}
\item By \cite[PROPOSITION~5.1]{PSTVA21}, for all $\alpha\in\mathcal{A}$, $D_\alpha$ is a prime divisor on $X$, and we have
\[
\Pic(X)=\bigoplus_{\alpha\in\mathcal{A}}\mathbb{Z}D_\alpha,\qquad\Eff^1(X)=\overline{\Eff}^1(X)=\bigoplus_{\alpha\in\mathcal{A}}\mathbb{R}_{\geq0}D_\alpha.
\]
\item We define a finite tuple of integers $\boldsymbol{\rho}=(\rho_\alpha)_{\alpha\in\mathcal{A}}\in\mathbb{Z}^{\mathcal{A}}$ by
\[
-K_X=\sum_{\alpha\in\mathcal{A}}\rho_\alpha D_\alpha,
\]
where $-K_X$ is the anticanonical divisor class of $X$.
By \cite[Lemma~2.4]{CLT02}, we have $\rho_\alpha\geq2$ for all $\alpha\in\mathcal{A}$.
\end{enumerate}
\end{rem}

\subsection{Darmon points}
In this subsection, we recall the definitions of Campana orbifolds and Darmon points.
Campana orbifolds were introduced in \cite[\S~1.2.1]{Cam04} and \cite[Definition~2.3]{Cam11}, and Darmon points were defined in \cite[Definition~2.10]{MNS24}.
We keep the notation introduced in \S~2 and \S~3.1.
We assume that the boundary divisor
\[
D:=\sum_{\alpha\in\mathcal{A}}D_\alpha
\]
is a strict normal crossings divisor.

\begin{dfn}
For each finite tuple of positive integers $\mathbf{m}=(m_\alpha)_{\alpha\in\mathcal{A}}$, we define a $\mathbb{Q}$-divisor $D_{\mathbf{m}}$ on $X$ by
\[
D_{\mathbf{m}}=\sum_{\alpha\in\mathcal{A}}\left(1-\frac{1}{m_\alpha}\right)D_\alpha.
\]
We call the pair $(X,D_{\mathbf{m}})$ a \textbf{Campana orbifold} over $F$.
\end{dfn}

Following \cite[after REMARK~3.2]{PSTVA21}, we first introduce integral models of Campana orbifolds in order to define Darmon points.

\begin{dfn}
Let $S$ be a finite set of places of $F$ containing all infinite places, let $\mathbf{m}=(m_\alpha)_{\alpha\in\mathcal{A}}$ be a tuple of positive integers, and let $(X,D_{\mathbf{m}})$ be a Campana orbifold over $F$ associated with $\mathbf{m}$.
A \textbf{good integral model} of $(X,D_{\mathbf{m}})$ away from $S$ is a pair $(\mathcal{X},\mathcal{D}_{\mathbf{m}})$ satisfying the following conditions:
\begin{enumerate}
\item $\mathcal{X}$ is a regular scheme, proper and flat over $\mathcal{O}_{F,S}$, and its generic fiber is isomorphic to $X$; that is,
\[
\mathcal{X}\otimes_{\mathcal{O}_{F,S}}F\cong X.
\]
\item $\mathcal{D}_{\mathbf{m}}$ is a $\mathbb{Q}$-Cartier divisor on $\mathcal{X}$ of the form
\[
\mathcal{D}_{\mathbf{m}}=\sum_{\alpha\in\mathcal{A}}\left(1-\frac{1}{m_\alpha}\right)\mathcal{D}_\alpha,
\]
where for each $\alpha\in\mathcal{A}$, $\mathcal{D}_\alpha$ is the Zariski closure of $D_\alpha$ in $\mathcal{X}$.
\end{enumerate}
\end{dfn}

For the remainder of this subsection, we fix a finite set $S$ of places of $F$ containing all infinite places, a tuple $\mathbf{m}=(m_\alpha)_{\alpha\in\mathcal{A}}$ of positive integers, and a good integral model $(\mathcal{X},\mathcal{D}_{\mathbf{m}})$ of $(X,D_{\mathbf{m}})$ away from $S$.
For each $\alpha\in\mathcal{A}$, $\mathcal{D}_\alpha$ denotes the Zariski closure of $D_\alpha$ in $\mathcal{X}$.
We then recall intersection multiplicities with the boundary divisors of an integral model.

\begin{dfn}
Let $\alpha\in\mathcal{A}$, let $v\in\Omega_F^{<\infty}\setminus S$, and let $P\in X(F_v)$.
We define the \textbf{intersection multiplicity} $n_v(\mathcal{D}_\alpha,P)$ of $P$ with $\mathcal{D}_\alpha$ at $v$ as follows.
The rational point $P$ uniquely induces an $\mathcal{O}_v$-rational point $\mathcal{P}_v\in\mathcal{X}(\mathcal{O}_v)$ by the valuative criterion for properness.
If $P\in D_\alpha(F_v)$, we set $n_v(\mathcal{D}_\alpha,P)=\infty$.
If $P\notin D_\alpha(F_v)$, let $\Spec(\mathcal{O}_v/\pi_v^N\mathcal{O}_v)$ be the closed subscheme of $\Spec\mathcal{O}_v$ determined by the pullback of $\mathcal{D}_\alpha$ via $\mathcal{P}_v$.
We define $n_v(\mathcal{D}_\alpha,P)=N$. 
\end{dfn}

Following \cite[Definition~2.10]{MNS24}, we now define Darmon points in this setting.

\begin{dfn}\label{Darmon}
Let $P\in(X\setminus\bigcup_{\alpha\in\mathcal{A}}D_\alpha)(F)$.
We say that $P$ is a \textbf{Darmon-}$\mathcal{O}_{F,S}$ \textbf{point} on $(\mathcal{X},\mathcal{D}_{\mathbf{m}})$ if
\[
n_v(\mathcal{D}_\alpha,P)\in m_\alpha\mathbb{Z}
\]
holds for all $\alpha\in\mathcal{A}$ and all $v\in\Omega_F^{<\infty}\setminus S$.
We denote the set of Darmon $\mathcal{O}_{F,S}$-points on $(\mathcal{X},\mathcal{D}_{\mathbf{m}})$ by $(\mathcal{X},\mathcal{D}_{\mathbf{m}})^{\mathrm{D}}(\mathcal{O}_{F,S})$, and its indicator function on $X(F)$ by $\delta_{\mathbf{m}}$.
\end{dfn}

\begin{rem}
\leavevmode
\begin{enumerate}
\item By the definition of Darmon points, we have
\[
(\mathcal{X},\mathcal{D}_{\mathbf{m}})^{\mathrm{D}}(\mathcal{O}_{F,S})\subseteq\left(X\setminus\bigcup_{\alpha\in\mathcal{A}}D_\alpha\right)(F)=G(F).
\]
\item Let
\[
\mathfrak{M}:=\prod_{\alpha\in\mathcal{A}}m_\alpha\mathbb{Z}_{\geq0}\subseteq\overline{\mathbb{N}}^{\mathcal{A}}
\]
and let
\[
\mathcal{M}=((\mathcal{D}_\alpha)_{\alpha\in\mathcal{A}},\mathfrak{M}).
\]
Then $(\mathcal{X},\mathcal{D}_{\mathbf{m}})^{\mathrm{D}}(\mathcal{O}_{F,S})$ is the set of $\mathcal{M}$-points on $(\mathcal{X},\mathcal{M})$ over $\mathcal{O}_{F,S}$, where $\mathcal{M}$-points were introduced in \cite[Definition~3.12]{Moe26a}.
\end{enumerate}
\end{rem}

\subsection{Heights}

We continue to use the notation fixed in the previous subsection.
In particular, $F$ is a number field and $X$ is a smooth projective equivariant compactification of the $n$-dimensional vector group $G=\mathbb{G}_{a,F}^n$ over $F$.
Moreover,
\[
X\setminus G=\bigcup_{\alpha\in\mathcal{A}}D_\alpha
\]
is the decomposition into irreducible components, $(X,D_{\mathbf{m}})$ is a Campana orbifold over $F$ associated with $\mathbf{m}\in\mathbb{Z}_{>0}^{\mathcal{A}}$, $S$ is a finite set of places of $F$ containing all infinite places, and $(\mathcal{X},\mathcal{D}_{\mathbf{m}})$ is a good integral model of $(X,D_{\mathbf{m}})$ away from $S$.
For each $\alpha\in\mathcal{A}$, let $\mathcal{O}(D_\alpha)$ be a line bundle on $X$ associated with $D_\alpha$ and let $\mathsf{f}_\alpha$ be a section of $\mathcal{O}(D_\alpha)$.
We fix a smooth adelic metric on $\mathcal{O}(D_\alpha)$ for each $\alpha\in\mathcal{A}$.

Following \cite[\S~4]{CLT02}, we recall the height functions on equivariant compactifications of vector groups.
We note that
\[
\Pic(X)=\bigoplus_{\alpha\in\mathcal{A}}\mathbb{Z}D_\alpha,\quad\Pic(X)_{\mathbb{C}}=\bigoplus_{\alpha\in\mathcal{A}}\mathbb{C}D_\alpha.
\]

\begin{dfn}
\leavevmode
\begin{enumerate}
\item For each $v\in\Omega_F$, we define the \textbf{local height pairing} $\mathsf{H}_v\colon G(F_v)\times\Pic(X)_{\mathbb{C}}\to\mathbb{C}^{\times}$ at $v$ by
\[
\mathsf{H}_v\left(\mathbf{x}_v,\sum_{\alpha\in\mathcal{A}}s_\alpha D_\alpha\right)=\prod_{\alpha\in\mathcal{A}}\|\mathsf{f}_\alpha(\mathbf{x}_v)\|_v^{-s_\alpha}.
\]
\item We define the \textbf{global height pairing} $\mathsf{H}\colon G(\mathbb{A}_F)\times\Pic(X)_{\mathbb{C}}\to\mathbb{C}^{\times}$ by
\[
\mathsf{H}(\mathbf{x},\mathbf{s})=\prod_{v\in\Omega_F}\mathsf{H}_v(\mathbf{x}_v,\mathbf{s}),
\]
where $\mathbf{x}_v$ denotes the image of $\mathbf{x}$ under the projection $G(\mathbb{A}_F)\to G(F_v)$.
\item Let $L\in\Pic(X)$.
We define the height function $H_L\colon G(F)\to\mathbb{R}_{>0}$ associated with $L$ by $H_L(P)=\mathsf{H}(P,L)$.
\end{enumerate}
\end{dfn}

Let $\mathbf{x}\in G(\mathbb{A}_F)$.
As in the proof of \cite[Lemma~2.1~(3)]{CLT10}, for all but finitely many finite places $v$ of $F$ and for all $\alpha\in\mathcal{A}$, we have $\|\mathsf{f}_\alpha(\mathbf{x}_v)\|_v=1$.
Hence, the infinite product defining the global height pairing is well defined.

If a line bundle $L$ on $X$ is big, then, by \cite[Proposition~4.3]{CLT02},
\[
\#\{P\in(\mathcal{X},\mathcal{D}_{\mathbf{m}})^{\mathrm{D}}(\mathcal{O}_{F,S})\mid H_L(P)\leq B\}
\]
is finite for every real number $B>0$.
In \S~4, we study the asymptotic behavior of the counting function.

\section{Proof of the main theorem}
In this section, we prove Theorem~\ref{main} using the method of \cite{PSTVA21}.

We continue to use the notation fixed in the previous section.
In particular, $F$ is a number field and $X$ is a smooth projective equivariant compactification of the $n$-dimensional vector group $G=\mathbb{G}_{a,F}^n$ over $F$.
Moreover,
\[
X\setminus G=\bigcup_{\alpha\in\mathcal{A}}D_\alpha
\]
is the decomposition into irreducible components, and we assume that the boundary divisor
\[
D=\sum_{\alpha\in\mathcal{A}}D_\alpha
\]
is a strict normal crossings divisor.
We fix a tuple $\mathbf{m}=(m_\alpha)_{\alpha\in\mathcal{A}}\in\mathbb{Z}_{>0}^{\mathcal{A}}$ of positive integers and define a $\mathbb{Q}$-divisor $D_{\mathbf{m}}$ on $X$ by
\[
D_{\mathbf{m}}=\sum_{\alpha\in\mathcal{A}}\left(1-\frac{1}{m_\alpha}\right)D_\alpha.
\]
We fix a finite set $S$ of places of $F$ containing all infinite places and a good integral model $(\mathcal{X},\mathcal{D}_{\mathbf{m}})$ of $(X,D_{\mathbf{m}})$ away from $S$.
For each $\alpha\in\mathcal{A}$, we fix a smooth adelic metric on the line bundle $\mathcal{O}(D_\alpha)$ associated with $D_\alpha$ and a section $\mathsf{f}_\alpha$ of $\mathcal{O}(D_\alpha)$.

\subsection{The height zeta function}
Following the method of \cite{PSTVA21}, we prove the main theorem by applying a Tauberian theorem.
Let $L$ be a big line bundle on $X$ and let $\{r_k\}_{k=1}^{\infty}$ be the sequence of distinct real numbers that occur as heights of Darmon $\mathcal{O}_{F,S}$-points on $(\mathcal{X},\mathcal{D}_{\mathbf{m}})$, arranged in increasing order.
Let
\[
a_k=\#\{P\in(\mathcal{X},\mathcal{D}_{\mathbf{m}})^{\mathrm{D}}(\mathcal{O}_{F,S})\mid H_L(P)=r_k\}
\]
for each $k$.
As in \cite{PSTVA21}, we prove Theorem~\ref{main} by applying the Tauberian theorem \cite{Del54} in the form stated in \cite[Theorem~6.1.4]{Tsc09}.

\begin{theorem}
Let $a$ and $\delta$ be positive real numbers, let $b$ be a positive integer, and let $\{a_k\}_{k=1}^\infty$ and $\{r_k\}_{k=1}^\infty$ be sequences of positive real numbers with $\lim_{k\to\infty}r_k=\infty$.
Assume that the following conditions hold:
\begin{enumerate}
\item for all $k\geq1$, the inequality $r_k<r_{k+1}$ holds, 
\item the Dirichlet series $\mathsf{Z}(s):=\sum_{k=1}^\infty a_k/r_k^s$ converges for $\Re s>a$,
\item $\mathsf{Z}(s)$ extends to a meromorphic function in the domain $\Re s>a-\delta$,
\item $\mathsf{Z}(s)$ has a rightmost pole of order $b$ at $s=a$, and
\item $c:=\lim_{s\to a}(s-a)^b\mathsf{Z}(s)$ is positive.
\end{enumerate}
Then the asymptotic formula
\[
\mathsf{N}(B):=\sum_{r_k\leq B}a_k\sim\frac{c}{a(b-1)!}B^a(\log B)^{b-1}\quad(B\to\infty)
\]
holds.
\end{theorem}

In order to apply the Tauberian theorem, we consider the Dirichlet series
\[
\mathsf{Z}_{\mathbf{m}}(sL):=\sum_{k=1}^\infty\frac{a_k}{r_k^s}.
\]
By an argument similar to that in the proof of \cite[Proposition~4.5]{CLT02}, for $\Re s$ sufficiently large, the series
\[
\sum_{\mathbf{x}\in(\mathcal{X},\mathcal{D}_{\mathbf{m}})^{\mathrm{D}}(\mathcal{O}_{F,S})}\frac{1}{\mathsf{H}(\mathbf{x},sL)}=\sum_{\mathbf{x}\in(\mathcal{X},\mathcal{D}_{\mathbf{m}})^{\mathrm{D}}(\mathcal{O}_{F,S})}\mathsf{H}(\mathbf{x},-sL)=\sum_{\mathbf{x}\in G(F)}\delta_{\mathbf{m}}(\mathbf{x})\mathsf{H}(\mathbf{x},-sL)
\]
converges absolutely.
We recall that $\delta_{\mathbf{m}}\colon G(F)\to\{0,1\}$ is the indicator function of $(\mathcal{X},\mathcal{D}_{\mathbf{m}})^{\mathrm{D}}(\mathcal{O}_{F,S})$.

Following \cite[\S~4]{CLT02}, we introduce the height zeta function.

\begin{dfn}
We define the \textbf{height zeta function} $\mathsf{Z}_\mathbf{m}(\mathbf{s})$ by
\[
\mathsf{Z}_\mathbf{m}(\mathbf{s}):=\sum_{\mathbf{x}\in G(F)}\delta_{\mathbf{m}}(\mathbf{x})\mathsf{H}(\mathbf{x},-\mathbf{s}).
\]
\end{dfn}

To analyze the height zeta function, we recall the Poisson summation formula.

\begin{dfn}
We fix the following notation and conventions.
\begin{enumerate}
\item For each infinite place $v$ of $F$, we define the local additive character $\psi_v\colon F_v\to\mathbb{C}^\times$ at $v$ by
\[
\psi_v(x):=\exp(-2\pi\sqrt{-1}\tr_{F_v/\mathbb{R}}(x)).
\]
\item For each prime number $p$, we define the local additive character $\psi_p\colon\mathbb{Q}_p\to\mathbb{C}^\times$ at $p$ by
\[
\psi_p(x):=\exp(2\pi\sqrt{-1}\overline{x}),
\]
where $\overline{x}$ is the image of $x\in\mathbb{Q}_p$ under the canonical map $\mathbb{Q}_p\to\mathbb{Q}_p/\mathbb{Z}_p\to\mathbb{Q}/\mathbb{Z}$.
\item For each finite place $v$ of $F$, we define the local additive character $\psi_v\colon F_v\to\mathbb{C}^\times$ at $v$ by the composition
\[
\psi_v:=\psi_{p_v}\circ\tr_{F_v/\mathbb{Q}_{p_v}}.
\]
\item We define a function $\psi_F\colon\mathbb{A}_F\to\mathbb{C}^\times$ by
\[
\psi_F(x):=\prod_{v\in\Omega_F}\psi_v(x_v),
\]
where, for each $v\in\Omega_F$, $x_v$ denotes the image of $x$ under the projection $\mathbb{A}_F\to F_v$.
This product is well defined, since $x_v\in\mathcal{O}_v$ for all but finitely many finite places $v$ and $\psi_v$ is trivial on $\mathcal{O}_v$.
\item Let $R$ be an $F$-algebra.
For each $\mathbf{a}=(a_i)_{i=1}^n\in G(R)=R^n$, we define a function $f_{\mathbf{a}}\colon G(R)\to R$ by
\[
f_{\mathbf{a}}((x_i)_{i=1}^n):=\sum_{i=1}^na_ix_i.
\]
\item Let $v\in\Omega_F$.
For each $F_v$-rational point $\mathbf{a}\in G(F_v)$, we define a function $\psi_{\mathbf{a},v}\colon G(F_v)\to\mathbb{C}^\times$ by the composition
\[
\psi_{\mathbf{a},v}:=\psi_v\circ f_{\mathbf{a}}.
\]
\item For each adelic point $\mathbf{a}\in G(\mathbb{A}_F)$, we define a function $\psi_{\mathbf{a}}\colon G(\mathbb{A}_F)\to\mathbb{C}^\times$ by
\[
\psi_{\mathbf{a}}:=\psi_F\circ f_{\mathbf{a}}.
\]
\item Let $\Phi\colon G(\mathbb{A}_F)\to\mathbb{C}$ be an absolutely integrable function.
We define the \textbf{Fourier transform} $\widehat{\Phi}\colon G(\mathbb{A}_F)\to\mathbb{C}$ of $\Phi$ by
\[
\widehat{\Phi}(\mathbf{a}):=\int_{G(\mathbb{A}_F)}\Phi(\mathbf{x})\psi_{\mathbf{a}}(\mathbf{x})\,\mathrm{d}\mathbf{x},
\]
where $\mathrm{d}\mathbf{x}$ is the Haar measure on $G(\mathbb{A}_F)=\mathbb{A}_F^n$ fixed in \S~2.
\end{enumerate}
\end{dfn}

With the above notation and conventions, we use the following adelic form of the Poisson summation formula, stated in \cite[LEMMA~4.2.4]{Tat67}.

\begin{theorem}\label{Poisson}
Let $\Phi\colon G(\mathbb{A}_F)\to\mathbb{C}$ be a function which satisfies the following conditions:
\begin{enumerate}
\item the function $\Phi$ is continuous and absolutely integrable on $G(\mathbb{A}_F)$,
\item the series $\sum_{\mathbf{x}\in G(F)}\Phi(\mathbf{x}+\mathbf{y})$ converges absolutely and uniformly as $\mathbf{y}$ ranges over a fundamental domain for
$G(\mathbb{A}_F)/G(F)$, and
\item the series $\sum_{\mathbf{a}\in G(F)}\widehat{\Phi}(\mathbf{a})$ converges absolutely.
\end{enumerate}
Then we have
\[
\sum_{\mathbf{a}\in G(F)}\widehat{\Phi}(\mathbf{a})=\sum_{\mathbf{x}\in G(F)}\Phi(\mathbf{x}).
\]
\end{theorem}

To apply the Poisson summation formula to the function $\delta_{\mathbf{m}}(\mathbf{x})\mathsf{H}(\mathbf{x},-\mathbf{s})$, we extend the indicator function $\delta_{\mathbf{m}}\colon G(F)\to\{0,1\}$ to a function on $G(\mathbb{A}_F)$.

\begin{dfn}\label{deltam}
For each $v\in\Omega_F$, we define a function $\delta_{\mathbf{m},v}\colon X(F_v)\to\{0,1\}$ as follows.
If $v\notin S$, we set
\[
\delta_{\mathbf{m},v}(\mathbf{x})=\left\{
\begin{array}{ll}
1&\text{if }n_v(\mathcal{D}_\alpha,\mathbf{x})\in m_\alpha\mathbb{Z}\text{ for all }\alpha\in\mathcal{A},\\
0&\text{otherwise}.
\end{array}\right.
\]
If $v\in S$, we set
\[
\delta_{\mathbf{m},v}(\mathbf{x})=\left\{
\begin{array}{ll}
1&\text{if }\mathbf{x}\in G(F_v),\\
0&\text{otherwise}.
\end{array}\right.
\]
We define a function $\delta_{\mathbf{m}}\colon G(\mathbb{A}_F)\to\{0,1\}$ by
\[
\delta_{\mathbf{m}}(\mathbf{x})=\prod_{v\in\Omega_F}\delta_{\mathbf{m},v}(\mathbf{x}_v),
\]
where $\mathbf{x}_v$ denotes the image of $\mathbf{x}$ under the projection $G(\mathbb{A}_F)\to G(F_v)$.
\end{dfn}

We note that for $v\in\Omega_F\setminus S$ and for $\mathbf{x}\in X(F_v)\setminus G(F_v)$, there exists $\alpha\in\mathcal{A}$ such that $\mathbf{x}\in D_\alpha(F_v)$.
It follows that $n_v(\mathcal{D}_\alpha,\mathbf{x})=\infty\notin m_\alpha\mathbb{Z}$, and hence we have $\delta_{\mathbf{m},v}(\mathbf{x})=0$.

\begin{rem}
By construction, the restriction of $\delta_{\mathbf{m}}$ to $G(F)$ is the indicator function of $(\mathcal{X},\mathcal{D}_{\mathbf{m}})^{\mathrm{D}}(\mathcal{O}_{F,S})$.
\end{rem}

\begin{rem}
By an argument similar to that in the proof of \cite[Lemma~4.5]{CLT02}, the function
\[
\Phi_{\mathbf{s}}(\mathbf{x}):=\delta_{\mathbf{m}}(\mathbf{x})\mathsf{H}(\mathbf{x},-\mathbf{s})
\]
satisfies the first two conditions of Theorem~\ref{Poisson} for $\Re\mathbf{s}$ sufficiently large.
In Proposition~\ref{conv}, we show that $\Phi_{\mathbf{s}}$ satisfies the third condition of Theorem~\ref{Poisson} for such~$\mathbf{s}$.
\end{rem}

As in \cite[\S~7]{PSTVA21}, we introduce the height integrals.
We recall that $\mathrm{d}\mathbf{x}_v$ is the Haar measure on $G(F_v)$ fixed in \S~2 for each $v\in\Omega_F$, and $\mathrm{d}\mathbf{x}$ is the Haar measure on $G(\mathbb{A}_F)=\mathbb{A}_F^n$ induced by the restricted product of the measures $\mathrm{d}\mathbf{x}_v$.

\begin{dfn}\label{HI}
Let $\mathbf{s}\in\Pic(X)_{\mathbb{C}}$.
\begin{enumerate}
\item Let $v\in\Omega_F$ and let $\mathbf{a}\in G(F_v)$.
We define the \textbf{local height integral} $\widehat{\mathsf{H}}_{\mathbf{m},v}(\mathbf{a},\mathbf{s})$ by
\[
\widehat{\mathsf{H}}_{\mathbf{m},v}(\mathbf{a},\mathbf{s}):=\int_{G(F_v)}\delta_{\mathbf{m},v}(\mathbf{x}_v)\mathsf{H}_v(\mathbf{x}_v,-\mathbf{s})\psi_{\mathbf{a},v}(\mathbf{x}_v)\,\mathrm{d}\mathbf{x}_v.
\]
\item Let $\mathbf{a}\in G(\mathbb{A}_F)$.
We define the \textbf{height integral} $\widehat{\mathsf{H}}_{\mathbf{m}}(\mathbf{a},\mathbf{s})$ by
\[
\widehat{\mathsf{H}}_{\mathbf{m}}(\mathbf{a},\mathbf{s}):=\int_{G(\mathbb{A}_F)}\delta_{\mathbf{m}}(\mathbf{x})\mathsf{H}(\mathbf{x},-\mathbf{s})\psi_{\mathbf{a}}(\mathbf{x})\,\mathrm{d}\mathbf{x}.
\]
\end{enumerate}
\end{dfn}

We now first show that
\[
\widehat{\mathsf{H}}_{\mathbf{m}}(\mathbf{a},\mathbf{s})=\prod_{v\in\Omega_F}\widehat{\mathsf{H}}_{\mathbf{m},v}(\mathbf{a}_v,\mathbf{s}).
\]
As in \cite[LEMMA~6.1]{PSTVA21}, we consider the following conditions, which will allow us to treat separately the places $v$ that satisfy them and those that do not.

\begin{con}\label{con}
For each $v\in\Omega_F$, we consider the following conditions.
\begin{enumerate}
\item The place $v$ is not contained in $S$.
\item For all $\alpha\in\mathcal{A}$, the metric on $D_\alpha\otimes_FF_v$ is induced by the integral model $(\mathcal{X}\otimes_{\mathcal{O}_{F,S}}\mathcal{O}_v,\mathcal{D}_\alpha\otimes_{\mathcal{O}_{F,S}}\mathcal{O}_v)$ (see \cite[\S~2.1.5]{CLT10} for the definition).
\item The scheme $\mathcal{X}\otimes_{\mathcal{O}_{F,S}}\mathcal{O}_v$ is smooth over $\mathcal{O}_v$, and $\mathcal{D}\otimes_{\mathcal{O}_{F,S}}\mathcal{O}_v$ is a relative strict normal crossings divisor (see \cite[\S~2]{IT14} for the definition).
Moreover, $\mathcal{X}\otimes_{\mathcal{O}_{F,S}}\mathcal{O}_v$ is equipped with an action of the $\mathcal{O}_v$-group scheme $\mathbb{G}_{a,\mathcal{O}_v}^n$ compatible with the action of $G$ on $X$.
\item For all $\alpha\in\mathcal{A}$, the unique linearisation of $\mathcal{O}(D_\alpha)$ extends to $\mathcal{O}(\mathcal{D}_\alpha)$.
\end{enumerate}
\end{con}

As in \cite[\S~3.1.8]{CLT12}, Condition~\ref{con} is satisfied for all but finitely many $v\in\Omega_F^{<\infty}\setminus S$.

\begin{rem}\label{Kv}
Let $v$ be a finite place of $F$.
By \cite[Proposition~4.2]{CLT02} and \cite[LEMMA~6.2]{PSTVA21}, there exists a compact open subgroup $K_v$ of $G(\mathcal{O}_v)$ such that, for every $\mathbf{s}\in\Pic(X)_{\mathbb{C}}$, the functions $\delta_{\mathbf{m},v}$ and $\mathsf{H}_v(\cdot,\mathbf{s})$ are $K_v$-invariant.
Moreover, if $v$ satisfies Condition~\ref{con}, we can take $K_v=G(\mathcal{O}_v)$.
\end{rem}

\begin{prop}
Let $\mathbf{a}\in G(\mathbb{A}_F)$, and let $\mathbf{s}\in\Pic(X)_{\mathbb{C}}$ with $\Re\mathbf{s}$ sufficiently large.
For all but finitely many $v\in\Omega_F^{<\infty}$ and for all $\mathbf{x}=(x_1,\cdots,x_n)\in G(\mathcal{O}_v)=\mathcal{O}_v^n$, we have
\[
\mathsf{H}_v(\mathbf{x},-\mathbf{s})=\psi_{\mathbf{a}_v,v}(\mathbf{x})=\delta_{\mathbf{m},v}(\mathbf{x})=\mu_v(\mathcal{O}_v)=1,
\]
where $\mathbf{a}_v$ is the image of $\mathbf{a}$ under the projection $G(\mathbb{A}_F)\to G(F_v)$.
In particular, we have
\[
\widehat{\mathsf{H}}_{\mathbf{m}}(\mathbf{a},\mathbf{s})=\prod_{v\in\Omega_F}\widehat{\mathsf{H}}_{\mathbf{m},v}(\mathbf{a}_v,\mathbf{s}).
\]
\end{prop}

\begin{proof}
The condition $\mu_v(\mathcal{O}_v)=1$ is equivalent to the condition that the prime ideal corresponding to $v$ does not divide the different of $F/\mathbb{Q}$, and hence it is satisfied for all but finitely many $v$.
We will show that
\[
\delta_{\mathbf{m},v}(\mathbf{x})=\mathsf{H}_v(\mathbf{x},-\mathbf{s})=\psi_{\mathbf{a}_v,v}(\mathbf{x})=1,
\]
for every finite place $v$ satisfying Condition~\ref{con} such that $\mu_v(\mathcal{O}_v)=1$ and $\mathbf{a}_v\in G(\mathcal{O}_v)$.
Let $\mathbf{x}=(x_1,\dots,x_n)\in G(\mathcal{O}_v)$.
By \cite[Lemma~2.1~(3)]{CLT10}, we have $\mathsf{H}_v(\mathbf{x},-\mathbf{s})=1$.
For each $1\leq i\leq n$ and each $v\in\Omega_F^{<\infty}$, we denote by $a_{i,v}$ the image of $\mathbf{a}_v$ under the $i$-th projection map $G(F_v)=F_v^n\to F_v$.
Since $\sum_{i=1}^na_{i,v}x_i\in\mathcal{O}_v$ and the prime ideal corresponding to $v$ does not divide the different of $F/\mathbb{Q}$, we have
\[
\psi_{\mathbf{a}_v,v}(\mathbf{x})=\psi_v\left(\sum_{i=1}^na_{i,v}x_i\right)=1.
\]
Let $P=(0,\cdots,0)\in G(\mathcal{O}_v)$ be the origin.
Then $P$ induces an $\mathcal{O}_v$-rational point $\mathcal{P}_v\colon\Spec\mathcal{O}_v\to\mathbb{G}_{a,\mathcal{O}_v}^n$.
It follows that, for all $\alpha\in\mathcal{A}$,
\[
\mathbb{G}_{a,\mathcal{O}_v}^n\times_{\mathcal{X}}\mathcal{D}_\alpha=\varnothing.
\]
Hence we have $\delta_{\mathbf{m},v}(P)=1$.
Since $\delta_{\mathbf{m},v}$ is $G(\mathcal{O}_v)$-invariant, we obtain $\delta_{\mathbf{m},v}|_{G(\mathcal{O}_v)}=1$.
\end{proof}

\begin{nota}
The following notation is based on \cite[\S~6]{PSTVA21}.
\begin{enumerate}
\item For each finite place $v$ of $F$ satisfying Condition~\ref{con}, set $K_v:=G(\mathcal{O}_v)$.
\item For each finite place $v$ of $F$ not satisfying Condition~\ref{con}, fix a compact subgroup $K_v$ of $G(\mathcal{O}_v)$ such that the two functions $\delta_{\mathbf{m},v}$ and $\mathsf{H}_v(\cdot,\mathbf{s})$ are $K_v$-invariant for all $\mathbf{s}\in\Pic(X)_{\mathbb{C}}$.
\item We define a subset $\mathbf{K}$ of $G(\mathbb{A}_F)$ by
\[
\mathbf{K}:=\prod_{v\in\Omega_F^{<\infty}}K_v\times\prod_{v\in\Omega_F^\infty}\{0\}.
\]
\item We define a subset $\Lambda_X$ of $G(F)$ by
\[
\Lambda_X:=\{\mathbf{a}\in G(\mathbb{A}_F)\mid \psi_{\mathbf{a}}|_{\mathbf{K}}=1\}\cap G(F).
\]
\end{enumerate}
\end{nota}

By an argument similar to that in \cite[PROPOSITION~5.3]{PSTVA21}, we have the following statement.

\begin{prop}\label{Hzero}
Let $\mathbf{a}\in G(F)\setminus\Lambda_X$ and let $\mathbf{s}\in\Pic(X)_{\mathbb{C}}$ such that $\mathsf{H}(\cdot,\mathbf{s})^{-1}$ is integrable on $G(\mathbb{A}_F)$.
Then we have $\widehat{\mathsf{H}}_{\mathbf{m}}(\mathbf{a},\mathbf{s})=0$.
\end{prop}

\begin{proof}
Since $\mathbf{a}\in G(F)\setminus\Lambda_X$, there exists an element $\mathbf{b}=(\mathbf{b}_v)_{v\in\Omega_F}\in\mathbf{K}$ such that
\[
1\neq\psi_{\mathbf{a}}(\mathbf{b})=\prod_{v\in\Omega_F}\psi_v(f_{\mathbf{a}}(\mathbf{b}_v)).
\]
Hence, there is a place $v$ of $F$ such that $\psi_v(f_{\mathbf{a}}(\mathbf{b}_v))\neq1$.
Since $\mathrm{d}\mathbf{x}_v$ is a Haar measure on $G(F_v)$, we obtain
\[
\widehat{\mathsf{H}}_{\mathbf{m},v}(\mathbf{a},\mathbf{s})=\int_{G(F_v)}\mathsf{H}_v(\mathbf{x}_v+\mathbf{b}_v,-\mathbf{s})\delta_{\mathbf{m},v}(\mathbf{x}_v+\mathbf{b}_v)\psi_v(f_{\mathbf{a}}(\mathbf{x}_v+\mathbf{b}_v))\,\mathrm{d}\mathbf{x}_v
\]
by translating the variables.
Since the two functions $\mathsf{H}_v$ and $\delta_{\mathbf{m},v}$ are $K_v$-invariant, we have
\begin{align*}
\widehat{\mathsf{H}}_{\mathbf{m},v}(\mathbf{a},\mathbf{s})&=\int_{G(F_v)}\mathsf{H}_v(\mathbf{x}_v,-\mathbf{s})\delta_{\mathbf{m},v}(\mathbf{x}_v)\psi_v(f_{\mathbf{a}}(\mathbf{x}_v+\mathbf{b}_v))\,\mathrm{d}\mathbf{x}_v
\\
&=\int_{G(F_v)}\mathsf{H}_v(\mathbf{x}_v,-\mathbf{s})\delta_{\mathbf{m},v}(\mathbf{x}_v)\psi_v(f_{\mathbf{a}}(\mathbf{x}_v)+f_{\mathbf{a}}(\mathbf{b}_v))\,\mathrm{d}\mathbf{x}_v
\\
&=\int_{G(F_v)}\mathsf{H}_v(\mathbf{x}_v,-\mathbf{s})\delta_{\mathbf{m},v}(\mathbf{x}_v)\psi_v(f_{\mathbf{a}}(\mathbf{x}_v))\psi_v(f_{\mathbf{a}}(\mathbf{b}_v))\,\mathrm{d}\mathbf{x}_v
\\
&=\psi_v(f_{\mathbf{a}}(\mathbf{b}_v))\int_{G(F_v)}\mathsf{H}_v(\mathbf{x}_v,-\mathbf{s})\delta_{\mathbf{m},v}(\mathbf{x}_v)\psi_v(f_{\mathbf{a}}(\mathbf{x}_v))\,\mathrm{d}\mathbf{x}_v
\\
&=\psi_v(f_{\mathbf{a}}(\mathbf{b}_v))\widehat{\mathsf{H}}_{\mathbf{m},v}(\mathbf{a},\mathbf{s}).
\end{align*}
Since $\psi_v(f_{\mathbf{a}}(\mathbf{b}_v))\neq1$, we have $\widehat{\mathsf{H}}_{\mathbf{m},v}(\mathbf{a},\mathbf{s})=0$, and hence $\widehat{\mathsf{H}}_{\mathbf{m}}(\mathbf{a},\mathbf{s})=0$.
\end{proof}

\subsection{The reduction maps}
We continue to use the notation fixed in the previous subsection.
In this subsection, following \cite[\S~6.2]{PSTVA21}, for each finite place $v\in\Omega_F^{<\infty}\setminus S$, we recall the reduction map $\eta_v\colon X(F_v)\to\mathcal{X}(k_v)$ in order to compute the local height integral
\[
\widehat{\mathsf{H}}_{\mathbf{m},v}(\mathbf{a},\mathbf{s})=\int_{G(F_v)}\delta_{\mathbf{m},v}(\mathbf{x}_v)\mathsf{H}_v(\mathbf{x}_v,-\mathbf{s})\psi_{\mathbf{a},v}(\mathbf{x}_v)\,\mathrm{d}\mathbf{x}_v
\]
explicitly.

\begin{dfn}
Let $v\in\Omega_F^{<\infty}\setminus S$ and let $P\in X(F_v)$.
By the valuative criterion for properness, the $F_v$-rational point $P$ uniquely induces an $\mathcal{O}_v$-rational point $\mathcal{P}_v\in\mathcal{X}(\mathcal{O}_v)$.
We denote by $\eta_v(P)\in\mathcal{X}(k_v)$ the composition of the morphism $\mathcal{P}_v\colon\Spec\mathcal{O}_v\to\mathcal{X}$ and the closed immersion $\Spec k_v\to\Spec\mathcal{O}_v$.
The map $\eta_v\colon X(F_v)\to\mathcal{X}(k_v)$ defined above is called the \textbf{reduction map} at $v$.
\end{dfn}

As in \cite[THEOREM~7.1]{PSTVA21}, we compute the local height integrals over $G(F_v)$ by decomposing the domain of integration into the fibers of the reduction map $\eta_v\colon X(F_v)\to\mathcal{X}(k_v)$.
Let $v\in\Omega_F^{<\infty}\setminus S$ be a finite place satisfying Condition~\ref{con}.
For $y\in\mathcal{X}(k_v)$, the integral over the fiber $\eta_v^{-1}(y)$ is determined by how $y$ meets the geometric irreducible components of the divisors $\mathcal{D}_\alpha\otimes k_v$.
To make this precise, we recall the terminology for the stratification of $\mathcal{X}(k_v)$, following \cite[\S~5--7]{PSTVA21}.

\begin{nota}\label{D}
Throughout this paper, we use the following notation for each $v\in\Omega_F$.
\begin{enumerate}
\item Let $D\otimes_FF_v=\bigcup_{\beta\in\mathcal{A}_v}D_{v,\beta}$ be the decomposition into irreducible components.
\item We fix an algebraic closure $\overline{F}$ of $F$.
For each $\alpha\in\mathcal{A}$, let $F_\alpha$ be the field of definition of one of the geometrically irreducible components of $D_\alpha$; equivalently, $F_\alpha$ is the algebraic closure of $F$ in the function field $F(D_\alpha)$.
\item We fix an algebraic closure $\overline{F_v}$ of $F_v$.
For each $\beta\in\mathcal{A}_v$, $\alpha(\beta)$ denotes the unique element $\alpha\in\mathcal{A}$ such that $\beta$ lies over $\alpha$, and $F_{v,\beta}$ denotes the field of definition of one of the geometrically irreducible components of $D_{v,\beta}$; equivalently, $F_{v,\beta}$ is the algebraic closure of $F_v$ in the function field $F_v(D_{v,\beta})$.
\item For each non-empty subset $B\subseteq\mathcal{A}_v$, we set
\[
 D_{v,B}=\bigcap_{\beta\in B}D_{v,\beta},\quad D_{v,B}^\circ=D_{v,B}\setminus\left(\bigcup_{B\subsetneq B'\subseteq\mathcal{A}_v}D_{v,B'}\right),
\]
with the convention that $D_{v,\varnothing}=X\otimes_FF_v$ and $D_{v,\varnothing}^\circ=G\otimes_FF_v$.
\end{enumerate}
\end{nota}

\begin{nota}
Throughout this paper, we use the following notation for each $v\in\Omega_F^{<\infty}\setminus S$.
\begin{enumerate}
\item For each $\beta\in\mathcal{A}_v$, $\mathcal{D}_{v,\beta}$ denotes the closure of $D_{v,\beta}$ in $\mathcal{X}\otimes_{\mathcal{O}_{F,S}}\mathcal{O}_v$.
\item For each subset $B$ of $\mathcal{A}_v$, let $\mathcal{D}_{v,B}$ be the closure of $D_{v,B}$ in $\mathcal{X}\otimes_{\mathcal{O}_{F,S}}\mathcal{O}_v$ and let
\[
\mathcal{D}_{v,B}^\circ=\mathcal{D}_{v,B}\setminus\left(\bigcup_{B\subsetneq B'\subseteq\mathcal{A}_v}\mathcal{D}_{v,B'}\right).
\]
\end{enumerate}
\end{nota}

Let $v\in\Omega_F^{<\infty}\setminus S$.
Since $\mathcal{X}(k_v)=\coprod_{B\subseteq\mathcal{A}_v}\mathcal{D}_{v,B}^\circ(k_v)$, we have the disjoint union
\[
X(F_v)=\eta_v^{-1}(\mathcal{X}(k_v))=\coprod_{B\subseteq\mathcal{A}_v}\eta_v^{-1}(\mathcal{D}_{v,B}^\circ(k_v))=\coprod_{B\subseteq\mathcal{A}_v}\coprod_{y\in\mathcal{D}_{v,B}^\circ(k_v)}\eta_v^{-1}(y).
\]
Using this disjoint union and the fact that the function $\delta_{\mathbf{m},v}$ vanishes on $X(F_v)\setminus G(F_v)$, we compute the local height integral $\widehat{\mathsf{H}}_{\mathbf{m},v}$ by decomposing it into contributions from the fibers of the reduction map.
By \cite[Theorem~2.13]{Sal98}, \cite[\S~6.2]{PSTVA21}, and the proof of \cite[THEOREM~7.1]{PSTVA21}, the local height integrals over the fibers of the reduction map can be calculated as follows.
\begin{rem}\label{Sal}
Let $v\in\Omega_F^{<\infty}\setminus S$ be a finite place of $F$ satisfying Condition~\ref{con}, and let $y\in\mathcal{X}(k_v)$.
Then there exist local analytic coordinates $(z_1,\dots,z_n)\colon\eta_v^{-1}(y)\to F_v^n$ such that the following conditions are satisfied.
\begin{enumerate}
\item Let $\tau_{X,v}$ be the Tamagawa measure on $X(F_v)$ (see \cite[\S~2.1.8]{CLT10} for the definition).
We equip $\eta_v^{-1}(y)$ with the measure induced by $\tau_{X,v}$, and we equip $\mathfrak{m}_v^n$ with the measure induced by the Haar measure on $F_v^n$ fixed in \S~2.
Then $(z_1,\cdots,z_n)$ induces a measure-preserving local analytic isomorphism $\eta_v^{-1}(y)\cong\mathfrak{m}_v^n$.
\item Let $B=\{\beta_1,\dots,\beta_r\}$ be a subset of $\mathcal{A}_v$ such that $y\in\mathcal{D}_{v,B}^\circ(k_v)$.
Then, for each $1\leq i\leq r$, the set $\eta_v^{-1}(y)\cap D_{v,\beta_i}(F_v)$ is cut out by $z_i=0$.
\end{enumerate}
For such $v\in\Omega_F^{<\infty}\setminus S$ and $y\in\mathcal{X}(k_v)$, we fix local analytic coordinates satisfying the above conditions.
We denote by $\mathrm{d}\mathbf{z}$ the measure on $\mathfrak{m}_v^n$ induced by the Haar measure on $F_v^n$ fixed in \S~2.
For each $1\leq i\leq n$, we denote by $\mathrm{d}z_i$ the measure on the $i$-th factor $\mathfrak{m}_v$, so that
\[
\mathrm{d}\mathbf{z}=\prod_{i=1}^n\mathrm{d}z_i.
\]
\end{rem}

We introduce notation and recall some results to describe how local heights, intersection multiplicities, and characters are expressed in the above local analytic coordinates.

\begin{nota}
Let $\mathbf{a}=(a_1,\dots,a_n)\in G(F)\setminus\{\mathbf{0}\}$.
We recall that the function $f_{\mathbf{a}}\colon G(F)\to F$ is defined by
\[
f_{\mathbf{a}}(x_1,\dots,x_n)=\sum_{i=1}^na_ix_i.
\]
We regard $f_{\mathbf{a}}$ as a rational function on $X$.
The following notation is based on \cite[\S~1]{CLT02} and \cite[\S~8]{PSTVA21}.
\begin{enumerate}
\item We denote by $E(f_{\mathbf{a}})$ the closure in $X$ of the hyperplane in $G$ defined by $f_{\mathbf{a}}=0$.
\item We denote by $\mathcal{E}(f_{\mathbf{a}})$ the closure in $\mathcal{X}$ of $E(f_{\mathbf{a}})$.
\item We define the tuple $(d_\alpha(f_{\mathbf{a}}))_{\alpha\in\mathcal{A}}$ of integers by
\[
E(f_{\mathbf{a}})=\sum_{\alpha\in\mathcal{A}}d_\alpha(f_{\mathbf{a}})D_\alpha+\DIv(f_{\mathbf{a}}).
\]
\item For each $v\in\Omega_F^{<\infty}$, we set
\[
j_v(\mathbf{a})=\min\{v(a_1),\dots,v(a_n)\}.
\]
\end{enumerate}
\end{nota}

\begin{rem}\label{Hnp}
Let $v\in\Omega_F^{<\infty}\setminus S$ be a finite place satisfying Condition~\ref{con}, let $\mathbf{x}\in G(F_v)$, let $\mathbf{s}=\sum_{\alpha\in\mathcal{A}}s_\alpha D_\alpha\in\Pic(X)_{\mathbb{C}}$, let $y=\eta_v(\mathbf{x})$, let $\mathbf{a}\in G(F)$, and let $B=\{\beta_1,\dots,\beta_r\}\subseteq\mathcal{A}_v$ such that $y\in\mathcal{D}_{v,B}^\circ(k_v)$.
For each $1\leq i\leq r$, we assume that $D_{v,\beta_i}(F_v)\cap\eta_v^{-1}(y)$ is defined by $z_i=0$.
\begin{enumerate}
\item We have
\[
\mathsf{H}_v(\mathbf{x},\mathbf{s})=\left\{
\begin{array}{ll}
1&\text{if }B=\varnothing, \\
\prod_{i=1}^r|z_i(\mathbf{x})|_v^{-s_{\alpha(\beta_i)}}&\text{if }B\neq\varnothing.
\end{array}\right.
\]
\item For all $\alpha\in\mathcal{A}$, we have
\[
\#(\mathcal{A}_v(\alpha)\cap B)=\left\{
\begin{array}{ll}
1&\text{if }n_v(\mathcal{D}_\alpha,\mathbf{x})>0, \\
0&\text{if }n_v(\mathcal{D}_\alpha,\mathbf{x})=0.
\end{array}\right.
\]
Moreover, we have $n_v(\mathcal{D}_\alpha,\mathbf{x})=v(z_i(\mathbf{x}))$ if $\mathcal{A}_v(\alpha)\cap B=\{\beta_i\}$.
\item We assume that $r=1$ and $\mathbf{a}\neq\mathbf{0}$.
If $y\notin\mathcal{E}(f_{\mathbf{a}})(k_v)$, then there exists an element $u\in\mathcal{O}_v^\times$ such that, for all $\mathbf{x}\in\eta_v^{-1}(y)$, the following equality holds:
\[
\psi_{\mathbf{a},v}(\mathbf{x})=\psi_v\left(\frac{u\pi_v^{j_v(\mathbf{a})}}{z_1(\mathbf{x})^{d_{\alpha(\beta_1)}(f_{\mathbf{a}})}}\right).
\]
\end{enumerate}
\end{rem}

\begin{proof}
(1) follows from the proof of \cite[THEOREM~7.1]{PSTVA21}.
The first assertion of (2) follows from Hensel's lemma and the fact that $D_\alpha$ is smooth over $F$.
The second assertion of (2) follows from \cite[\S~6.2]{PSTVA21} and the proof of \cite[Proposition~7.1]{MNS24}.
(3) follows from the proof of \cite[Proposition~10.2]{CLT02}.
\end{proof}

\subsection{The height integral when $\mathbf{a}=\mathbf{0}$}

We continue to use the notation fixed in the previous subsection.
In this subsection, we study the local height integral
\[
\widehat{\mathsf{H}}_{\mathbf{m},v}(\mathbf{0},\mathbf{s})=\int_{G(F_v)}\delta_{\mathbf{m},v}(\mathbf{x}_v)\mathsf{H}_v(\mathbf{x}_v,-\mathbf{s})\,\mathrm{d}\mathbf{x}_v
\]
at $v\in\Omega_F$ and its Euler product
\[
\widehat{\mathsf{H}}_{\mathbf{m}}(\mathbf{0},\mathbf{s})=\prod_{v\in\Omega_F}\widehat{\mathsf{H}}_{\mathbf{m},v}(\mathbf{0},\mathbf{s})
\]
following \cite[\S~7]{PSTVA21}.
We identify $\Pic(X)_{\mathbb{C}}\cong\bigoplus_{\alpha\in\mathcal{A}}\mathbb{C}D_\alpha$ with $\mathbb{C}^{\mathcal{A}}$.
As stated in the previous subsection, if $v\in\Omega_F^{<\infty}\setminus S$ satisfies Condition~\ref{con}, then the local height integral at $v$ can be calculated using the reduction map $\eta_v$.

First, we establish the following finiteness lemma in order to choose a positive integer $N$ such that every finite place $v$ with $q_v>N$ is good.

\begin{lem}\label{N}
The following properties hold for all but finitely many finite places $v$ of $F$:
\begin{enumerate}
\item $v$ satisfies Condition~\ref{con},
\item $p_v\neq2,3$,
\item the field extension $F_v/\mathbb{Q}_{p_v}$ is unramified, and
\item $j_v(\mathbf{a})\geq0$ for all $\mathbf{a}\in\Lambda_X\setminus\{\mathbf{0}\}$.
\end{enumerate}
\end{lem}

\begin{proof}
The assertions (1), (2), and (3) are immediate.
(4) follows from the fact that $\Lambda_X$ is a finitely generated $\mathcal{O}_F$-module.
\end{proof}

\begin{nota}
\leavevmode
\begin{enumerate}
\item We fix a positive integer $N$ such that all finite places $v$ with $q_v>N$ satisfy the conditions of Lemma~\ref{N}.
\item For each real number $c$, we define a subset $\mathsf{T}_{>c}$ of $\Pic(X)_{\mathbb{C}}$ by
\[
\mathsf{T}_{>c}=\{\mathbf{s}=(s_\alpha)_{\alpha\in\mathcal{A}}\mid\Re s_\alpha>\rho_\alpha-1+m_\alpha^{-1}+c\text{ for all }\alpha\in\mathcal{A}\}
\]
as in \cite[\S~7]{PSTVA21}.
\item We fix a positive constant $\delta$ such that $0<\delta<(3\max_{\alpha\in\mathcal{A}}m_\alpha)^{-1}$.
\end{enumerate}
\end{nota}

Next, following \cite[THEOREM~7.1]{PSTVA21}, we compute the local height integral $\widehat{\mathsf{H}}_{\mathbf{m},v}(\mathbf{0},\mathbf{s})$.
We recall that the tuple $\boldsymbol{\rho}=(\rho_\alpha)_{\alpha\in\mathcal{A}}$ is defined by
\[
-K_X=\sum_{\alpha\in\mathcal{A}}\rho_\alpha D_\alpha.
\]
We recall that $q_v=\#(\mathcal{O}_v/\mathfrak{m}_v)$.
Let $\tau_{X,v}$ be the Tamagawa measure on $X(F_v)$.
For each finite place $v\in\Omega_F^{<\infty}\setminus S$ satisfying Condition~\ref{con} and each $y\in\mathcal{X}(k_v)$, we use the notation $\mathrm{d}\mathbf{x}_v$ and $\mathrm{d}\mathbf{z}=\prod_{i=1}^n\mathrm{d}z_i$.
In particular, on $G(F_v)$, we have
\[
\mathrm{d}\mathbf{x}_v=\mathsf{H}_v(\mathbf{x}_v,\boldsymbol{\rho})\,\mathrm{d}\tau_{X,v}.
\]

\begin{prop}\label{Hvzero}
Let $\mathbf{s}\in\mathsf{T}_{>-\delta}$ and let $v\in\Omega_F^{<\infty}\setminus S$ be a finite place of $F$ satisfying conditions (1) and (3) of Lemma~\ref{N}.
Then we have
\[
\widehat{\mathsf{H}}_{\mathbf{m},v}(\mathbf{0},\mathbf{s})=\sum_{B\subseteq\mathcal{A}_v}\frac{\#\mathcal{D}_{v,B}^\circ(k_v)}{q_v^{n-\#B}}\prod_{\beta\in B}\left(1-\frac{1}{q_v}\right)\frac{q_v^{-m_{\alpha(\beta)}(s_{\alpha(\beta)}-\rho_{\alpha(\beta)}+1)}}{1-q_v^{-m_{\alpha(\beta)}(s_{\alpha(\beta)}-\rho_{\alpha(\beta)}+1)}}.
\]
\end{prop}

\begin{proof}
We have
\begin{align*}
\widehat{\mathsf{H}}_{\mathbf{m},v}(\mathbf{0},\mathbf{s})&=\int_{G(F_v)}\delta_{\mathbf{m},v}(\mathbf{x}_v)\mathsf{H}_v(\mathbf{x}_v,\mathbf{s})^{-1}\psi_{\mathbf{0},v}(\mathbf{x}_v)\,\mathrm{d}\mathbf{x}_v
\\
&=\int_{G(F_v)}\delta_{\mathbf{m},v}(\mathbf{x}_v)\mathsf{H}_v(\mathbf{x}_v,\mathbf{s})^{-1}\mathsf{H}_v(\mathbf{x}_v,\boldsymbol{\rho})\frac{\mathrm{d}\mathbf{x}_v}{\mathsf{H}_v(\mathbf{x}_v,\boldsymbol{\rho})}
\\
&=\int_{X(F_v)}\delta_{\mathbf{m},v}(\mathbf{x}_v)\mathsf{H}_v(\mathbf{x}_v,\mathbf{s}-\boldsymbol{\rho})^{-1}\,\mathrm{d}\tau_{X,v},
\end{align*}
where in the last integral the integrand is extended by zero to $X(F_v)\setminus G(F_v)$.
Using the disjoint union
\[
X(F_v)=\coprod_{B\subseteq\mathcal{A}_v}\coprod_{y\in\mathcal{D}_{v,B}^\circ(k_v)}\eta_v^{-1}(y),
\]
we decompose the above integral over $X(F_v)$ into integrals over the fibers $\eta_v^{-1}(y)$ as follows:
\begin{equation}\label{1}
\widehat{\mathsf{H}}_{\mathbf{m},v}(\mathbf{0},\mathbf{s})=\sum_{B\subseteq\mathcal{A}_v}\sum_{y\in\mathcal{D}_{v,B}^\circ(k_v)}\int_{\eta_v^{-1}(y)}\delta_{\mathbf{m},v}(\mathbf{x}_v)\mathsf{H}_v(\mathbf{x}_v,\mathbf{s}-\boldsymbol{\rho})^{-1}\,\mathrm{d}\tau_{X,v}.
\end{equation}

Let $B=\{\beta_1,\dots,\beta_r\}$, let $y\in\mathcal{D}_{v,B}^\circ(k_v)$, and let $(z_1,\dots,z_n)$ be the local analytic coordinates on $\eta_v^{-1}(y)$ fixed in Remark~\ref{Sal}.
Let $\mathbf{x}_v\in\eta_v^{-1}(y)$.
We write $z_i=z_i(\mathbf{x}_v)$ for each $1\leq i\leq r$, with the convention that $v(0)=\infty$.
By Remark~\ref{Hnp}~(1) and (2), we have
\begin{align*}
\delta_{\mathbf{m},v}(\mathbf{x}_v)&=\left\{
\begin{array}{ll}
1&\text{if }v(z_i)\in m_{\alpha(\beta_i)}\mathbb{Z}\text{ for all }i\in\{1,\dots,r\},\\
0&\text{otherwise},
\end{array}\right. \\
&=\prod_{i=1}^r\mathbf{1}_{\coprod_{k=1}^\infty\pi_v^{km_{\alpha(\beta_i)}}\mathcal{O}_v^{\times}}(z_i),
\end{align*}
and for $\mathbf{x}_v\in G(F_v)$, we have
\[
\mathsf{H}_v(\mathbf{x}_v,\mathbf{s}-\boldsymbol{\rho})^{-1}=\prod_{i=1}^r|z_i|_v^{s_{\alpha(\beta_i)}-\rho_{\alpha(\beta_i)}}.
\]
Hence, we have
\begin{align}
&\int_{\eta_v^{-1}(y)}\delta_{\mathbf{m},v}(\mathbf{x}_v)\mathsf{H}_v(\mathbf{x}_v,\mathbf{s}-\boldsymbol{\rho})^{-1}\,\mathrm{d}\tau_{X,v}\notag
\\
&=\int_{\mathfrak{m}_v^n}\prod_{i=1}^r\left(\mathbf{1}_{\coprod_{k=1}^\infty\pi_v^{km_{\alpha(\beta_i)}}\mathcal{O}_v^\times}(z_i)|z_i|_v^{s_{\alpha(\beta_i)}-\rho_{\alpha(\beta_i)}}\right)\,\mathrm{d}\mathbf{z}\notag
\\
&=\left(\prod_{i=1}^r\sum_{k=1}^\infty\int_{\pi_v^{km_{\alpha(\beta_i)}}\mathcal{O}_v^\times}|z_i|_v^{s_{\alpha(\beta_i)}-\rho_{\alpha(\beta_i)}}\,\mathrm{d}z_i\right)\left(\prod_{i=r+1}^n\int_{\mathfrak{m}_v}1\,\mathrm{d}z_i\right)\notag
\\
&=\frac{1}{q_v^{n-r}}\prod_{i=1}^r\sum_{k=1}^\infty\mu_v(\pi_v^{km_{\alpha(\beta_i)}}\mathcal{O}_v^\times)q_v^{-km_{\alpha(\beta_i)}(s_{\alpha(\beta_i)}-\rho_{\alpha(\beta_i)})}\notag
\\
&=\frac{1}{q_v^{n-r}}\prod_{i=1}^r\left(1-\frac{1}{q_v}\right)\sum_{k=1}^\infty q_v^{-km_{\alpha(\beta_i)}(s_{\alpha(\beta_i)}-\rho_{\alpha(\beta_i)}+1)}\notag
\\
&=\frac{1}{q_v^{n-r}}\prod_{i=1}^r\left(1-\frac{1}{q_v}\right)\frac{q_v^{-m_{\alpha(\beta_i)}(s_{\alpha(\beta_i)}-\rho_{\alpha(\beta_i)}+1)}}{1-q_v^{-m_{\alpha(\beta_i)}(s_{\alpha(\beta_i)}-\rho_{\alpha(\beta_i)}+1)}}\notag
\\
&=\frac{1}{q_v^{n-\#B}}\prod_{\beta\in B}\left(1-\frac{1}{q_v}\right)\frac{q_v^{-m_{\alpha(\beta)}(s_{\alpha(\beta)}-\rho_{\alpha(\beta)}+1)}}{1-q_v^{-m_{\alpha(\beta)}(s_{\alpha(\beta)}-\rho_{\alpha(\beta)}+1)}}\label{2}
\end{align}
for all $B\subseteq\mathcal{A}_v$ and all $y\in\mathcal{D}_{v,B}^\circ(k_v)$.
The fifth equality follows from the fact that $\mathbf{s}\in\mathsf{T}_{>-\delta}$ and
\[
m_{\alpha(\beta_i)}(\Re s_{\alpha(\beta_i)}-\rho_{\alpha(\beta_i)}+1)>0
\]
for all $1\leq i\leq r$.
Combining \eqref{1} and \eqref{2} completes the proof.
\end{proof}

To study the convergence of the height integral
\[
\widehat{\mathsf{H}}_{\mathbf{m}}(\mathbf{a},\mathbf{s})=\prod_{v\in\Omega_F}\widehat{\mathsf{H}}_{\mathbf{m},v}(\mathbf{a},\mathbf{s}),
\]
we recall the definitions of the relevant local and global zeta functions.

\begin{nota}\label{zeta}
\leavevmode
\begin{enumerate}
\item For each local field $K$, we define the local zeta function $\zeta_K$ of $K$ by
\[
\zeta_K(s)=\left\{
\begin{array}{ll}
s^{-1}&\text{if }K\text{ is Archimedean},\\
(1-q^{-s})^{-1}&\text{if }K\text{ is non-Archimedean},
\end{array}
\right.
\]
where $q$ is the cardinality of the residue field of $K$, as in \cite[\S~2.1]{PSTVA21}.
\item For each number field $K$, we define the Dedekind zeta function $\zeta_K$ of $K$ by
\[
\zeta_K(s)=\prod_{v\in\Omega_K^{<\infty}}\zeta_{K_v}(s)
\]
as in \cite[\S~2.1]{PSTVA21}.
\item For each $\alpha\in\mathcal{A}$, we define the function $\zeta_{F_\alpha,S^c}$ by
\[
\zeta_{F_\alpha,S^c}(s)=\prod_{v\in\Omega_F^{<\infty}\setminus S}\prod_{\beta\in\mathcal{A}_v(\alpha)}\zeta_{F_{v,\beta}}(s)
\]
as in \cite[\S~8.3]{PSTVA21}, where $F_{v,\beta}$ is the field of definition of one of the geometrically irreducible components of $D_{v,\beta}$ defined in Notation~\ref{D}~(3).
\end{enumerate}
\end{nota}

\begin{rem}\label{z}
Let $\alpha\in\mathcal{A}$.
By the proof of \cite[COROLLARY~7.5]{PSTVA21}, we have
\[
\prod_{w\mid v}(F_\alpha)_w\cong F_\alpha\otimes_FF_v=\prod_{\beta\in\mathcal{A}_v(\alpha)}F_{v,\beta}
\]
for all $v\in\Omega_F$, where $w\mid v$ means that $w$ runs over all places of $F_\alpha$ above $v$ and $(F_\alpha)_w$ is the $w$-adic completion of $F_\alpha$.
Hence, we have $\{(F_\alpha)_w\}_{w\mid v}=\{F_{v,\beta}\}_{\beta\in\mathcal{A}_v(\alpha)}$ and
\begin{align*}
\zeta_{F_\alpha}(s)&=\prod_{w\in\Omega_{F_\alpha}^{<\infty}}\zeta_{(F_\alpha)_w}(s)
\\
&=\prod_{v\in\Omega_F^{<\infty}}\prod_{w\mid v}\zeta_{(F_\alpha)_w}(s)
\\
&=\prod_{v\in\Omega_F^{<\infty}}\prod_{\beta\in\mathcal{A}_v(\alpha)}\zeta_{F_{v,\beta}}(s)
\\
&=\zeta_{F_\alpha,S^c}(s)\prod_{v\in S\cap\Omega_F^{<\infty}}\prod_{\beta\in\mathcal{A}_v(\alpha)}\zeta_{F_{v,\beta}}(s).
\end{align*}
\end{rem}

The following proposition is an analogue of \cite[PROPOSITION~7.4]{PSTVA21}.

\begin{prop}\label{Hvzero2}
Let $v$ be a finite place of $F$ such that $q_v>N$.
Then, uniformly for $\mathbf{s}\in\mathsf{T}_{>-\delta}$, we have
\[
\widehat{\mathsf{H}}_{\mathbf{m},v}(\mathbf{0},\mathbf{s})\prod_{\alpha\in\mathcal{A}}\prod_{\beta\in\mathcal{A}_v(\alpha)}\zeta_{F_{v,\beta}}(m_\alpha(s_\alpha-\rho_\alpha+1))^{-1}=1+O(q_v^{-1-\delta})\quad(q_v\to\infty).
\]
\end{prop}

\begin{proof}
Since $q_v>N$, the place $v$ satisfies the conditions of Lemma~\ref{N}.
For each subset $B=\{\beta_1,\dots,\beta_r\}$ of $\mathcal{A}_v$, we set
\begin{align*}
A_B(\mathbf{s})&:=\sum_{y\in\mathcal{D}_{v,B}^\circ(k_v)}\int_{\eta_v^{-1}(y)}\delta_{\mathbf{m},v}(\mathbf{x}_v)\mathsf{H}_v(\mathbf{x}_v,\mathbf{s}-\boldsymbol{\rho})^{-1}\,\mathrm{d}\tau_{X,v}
\\
&=\frac{\#\mathcal{D}_{v,B}^\circ(k_v)}{q_v^{n-\#B}}\prod_{\beta\in B}\left(1-\frac{1}{q_v}\right)\frac{q_v^{-m_{\alpha(\beta)}(s_{\alpha(\beta)}-\rho_{\alpha(\beta)}+1)}}{1-q_v^{-m_{\alpha(\beta)}(s_{\alpha(\beta)}-\rho_{\alpha(\beta)}+1)}}.
\end{align*}
By Proposition~\ref{Hvzero}, we have
\[
\widehat{\mathsf{H}}_{\mathbf{m},v}(\mathbf{0},\mathbf{s})=\sum_{B\subseteq\mathcal{A}_v}A_B(\mathbf{s}).
\]
We consider $A_B(\mathbf{s})$ in several cases.

\begin{enumerate}
\item Case $B=\varnothing$.

By Condition~\ref{con}~(3), the scheme $\mathcal{X}\otimes_{\mathcal{O}_{F,S}}\mathcal{O}_v$ contains $\mathbb{G}_{a,\mathcal{O}_v}^n$ as a dense open subscheme, and its complement is
\[
\bigcup_{\varnothing\neq B\subseteq\mathcal{A}_v}\mathcal{D}_{v,B}=\bigcup_{\varnothing\neq B\subseteq\mathcal{A}_v}\mathcal{D}_{v,B}^\circ.
\]
Hence, we have $\#\mathcal{D}_{v,\varnothing}^\circ(k_v)=\#\mathbb{G}_{a,\mathcal{O}_v}^n(k_v)=q_v^n$ and $A_\varnothing(\mathbf{s})=1$.
\item Case $r=1$ and $f_{v,\beta_1}\neq1$.

We recall that $f_{v,\beta_1}=[F_{v,\beta_1}:F_v]$.
By Condition~\ref{con}~(3), the scheme $D_{v,\beta_1}$ is smooth over $F_v$.
Since $f_{v,\beta_1}\neq1$, $D_{v,\beta_1}$ is not geometrically irreducible over $F_v$.
By Hensel's lemma, we have $\#\mathcal{D}_{v,B}(k_v)=0$ and $\#\mathcal{D}_{v,B}^\circ(k_v)=0$.
It follows that $A_B(\mathbf{s})=0$.
\item Case $r=1$ and $f_{v,\beta_1}=1$.

By Condition~\ref{con}~(3), $\mathcal{D}_{v,\beta_1}$ is smooth over $\mathcal{O}_v$.
The scheme $D_{v,\beta_1}$ is geometrically irreducible over $F_v$ since $f_{v,\beta_1}=1$.
Hence, $\mathcal{D}_{v,\beta_1}\otimes_{\mathcal{O}_v}k_v$ is also geometrically irreducible over $k_v$.
By the Lang--Weil estimates \cite{LW54}, uniformly for such $B$, we have
\[
\frac{\#\mathcal{D}_{v,B}^\circ(k_v)}{q_v^{n-\#B}}\left(1-\frac{1}{q_v}\right)=1+O(q_v^{-1/2})\quad(q_v\to\infty).
\]
Since
\[
-2m_{\alpha(\beta_1)}(\Re s_{\alpha(\beta_1)}-\rho_{\alpha(\beta_1)}+1)<2m_{\alpha(\beta_1)}\delta-2<-1-\delta
\]
for $\mathbf{s}\in\mathsf{T}_{>-\delta}$, we have
\[
\frac{q_v^{-2m_{\alpha(\beta_1)}(s_{\alpha(\beta_1)}-\rho_{\alpha(\beta_1)}+1)}}{1-q_v^{-m_{\alpha(\beta_1)}(s_{\alpha(\beta_1)}-\rho_{\alpha(\beta_1)}+1)}}=O(q_v^{-1-\delta})\quad(q_v\to\infty).
\]
It follows that
\begin{align*}
&A_B(\mathbf{s})
\\
&=\frac{\#\mathcal{D}_{v,B}^\circ(k_v)}{q_v^{n-\#B}}\left(1-\frac{1}{q_v}\right)\left(q_v^{-m_{\alpha(\beta_1)}(s_{\alpha(\beta_1)}-\rho_{\alpha(\beta_1)}+1)}+\frac{q_v^{-2m_{\alpha(\beta_1)}(s_{\alpha(\beta_1)}-\rho_{\alpha(\beta_1)}+1)}}{1-q_v^{-m_{\alpha(\beta_1)}(s_{\alpha(\beta_1)}-\rho_{\alpha(\beta_1)}+1)}}\right)
\\
&=q_v^{-m_{\alpha(\beta_1)}(s_{\alpha(\beta_1)}-\rho_{\alpha(\beta_1)}+1)}+O(q_v^{-1-\delta})\quad(q_v\to\infty).
\end{align*}
\item Case $r\geq2$.

By the Lang--Weil estimates, uniformly for such $B$, we have
\[
\frac{\#\mathcal{D}_{v,B}^\circ(k_v)}{q_v^{n-\#B}}\left(1-\frac{1}{q_v}\right)=O(1)\quad(q_v\to\infty).
\]
Since $r=\#B\geq2$, we have
\[
\prod_{\beta\in B}\left(1-\frac{1}{q_v}\right)\frac{q_v^{-m_{\alpha(\beta)}(s_{\alpha(\beta)}-\rho_{\alpha(\beta)}+1)}}{1-q_v^{-m_{\alpha(\beta)}(s_{\alpha(\beta)}-\rho_{\alpha(\beta)}+1)}}=O(q_v^{-1-\delta})\quad(q_v\to\infty)
\]
and
\[
-\sum_{\beta\in B}m_{\alpha(\beta)}(\Re s_{\alpha(\beta)}-\rho_{\alpha(\beta)}+1)<-1-\delta.
\]
It follows that $A_B(\mathbf{s})=O(q_v^{-1-\delta})$ as $q_v\to\infty$.
\end{enumerate}
Collecting these estimates, we obtain
\begin{equation}\label{3}
\widehat{\mathsf{H}}_{\mathbf{m},v}(\mathbf{0},\mathbf{s})=1+\sum_{\beta\in\mathcal{A}_v,f_{v,\beta}=1}q_v^{-m_{\alpha(\beta)}(s_{\alpha(\beta)}-\rho_{\alpha(\beta)}+1)}+O(q_v^{-1-\delta})\quad(q_v\to\infty).
\end{equation}
On the other hand, by expanding the local zeta factors, we have
\begin{equation}\label{4}
\prod_{\alpha\in\mathcal{A}}\prod_{\beta\in\mathcal{A}_v(\alpha)}\zeta_{F_{v,\beta}}(m_\alpha(s_\alpha-\rho_\alpha+1))^{-1}=1-\sum_{\beta\in\mathcal{A}_v,f_{v,\beta}=1}q_v^{-m_{\alpha(\beta)}(s_{\alpha(\beta)}-\rho_{\alpha(\beta)}+1)}+O(q_v^{-1-\delta}),
\end{equation}
as $q_v\to\infty$.
Combining the asymptotic formulae \eqref{3} and \eqref{4} completes the proof.
\end{proof}

Finally, we study the analytic properties of the height integral
\[
\widehat{\mathsf{H}}_{\mathbf{m}}(\mathbf{0},\mathbf{s})=\prod_{v\in\Omega_F}\widehat{\mathsf{H}}_{\mathbf{m},v}(\mathbf{0},\mathbf{s}).
\]

\begin{rem}\label{remzh}
By an argument similar to the proof of \cite[PROPOSITION~7.2]{PSTVA21} and \cite[PROPOSITION~7.3]{PSTVA21}, the function
\[
\mathbf{s}\mapsto\widehat{\mathsf{H}}_{\mathbf{m},v}(\mathbf{0},\mathbf{s})
\]
is holomorphic on $\mathsf{T}_{>-\delta}$ for all $v\in\Omega_F$.
\end{rem}

The following corollary is an analogue of \cite[COROLLARY~7.5]{PSTVA21}.

\begin{co}\label{fzero}
The infinite product
\begin{align*}
&f^{\mathbf{0}}(\mathbf{s})
\\
&:=\prod_{v\in\Omega_F^{<\infty}}\left(\widehat{\mathsf{H}}_{\mathbf{m},v}(\mathbf{0},\mathbf{s})\prod_{\alpha\in\mathcal{A}}\prod_{w\mid v}\zeta_{(F_\alpha)_w}(m_\alpha(s_\alpha-\rho_\alpha+1))^{-1}\right)\times\prod_{v\in\Omega_F^{\infty}}\widehat{\mathsf{H}}_{\mathbf{m},v}(\mathbf{0},\mathbf{s})
\\
&=\widehat{\mathsf{H}}_{\mathbf{m}}(\mathbf{0},\mathbf{s})\prod_{\alpha\in\mathcal{A}}\zeta_{F_\alpha}(m_\alpha(s_\alpha-\rho_\alpha+1))^{-1}
\end{align*}
converges absolutely and uniformly on compact subsets of $\mathsf{T}_{>-\delta}$ and defines a holomorphic function on $\mathsf{T}_{>-\delta}$.
In particular, the function
\[
\mathbf{s}\mapsto f^{\mathbf{0}}(\mathbf{s})\prod_{\alpha\in\mathcal{A}}\zeta_{F_\alpha}(m_\alpha(s_\alpha-\rho_\alpha+1))
\]
gives a meromorphic continuation of $\widehat{\mathsf{H}}_{\mathbf{m}}(\mathbf{0},\mathbf{s})$ to $\mathsf{T}_{>-\delta}$.
\end{co}

\begin{proof}
By Proposition~\ref{Hvzero2}, the above infinite product converges absolutely and uniformly on compact subsets of $\mathsf{T}_{>-\delta}$.
By Remark~\ref{remzh}, the function
\[
\mathbf{s}\mapsto\widehat{\mathsf{H}}_{\mathbf{m},v}(\mathbf{0},\mathbf{s})\prod_{\alpha\in\mathcal{A}}\prod_{w\mid v}\zeta_{(F_\alpha)_w}(m_\alpha(s_\alpha-\rho_\alpha+1))^{-1}
\]
is holomorphic on $\mathsf{T}_{>-\delta}$ for all $v\in\Omega_F^{<\infty}$.
Moreover, by Remark~\ref{remzh}, the function
\[
\mathbf{s}\mapsto
\widehat{\mathsf{H}}_{\mathbf{m},v}(\mathbf{0},\mathbf{s})
\]
is holomorphic on $\mathsf{T}_{>-\delta}$ for all $v\in\Omega_F^\infty$.
Hence, $f^{\mathbf{0}}(\mathbf{s})$ defines a holomorphic function on $\mathsf{T}_{>-\delta}$.
The final assertion follows from the identity
\[
\widehat{\mathsf{H}}_{\mathbf{m}}(\mathbf{0},\mathbf{s})=f^{\mathbf{0}}(\mathbf{s})\prod_{\alpha\in\mathcal{A}}\zeta_{F_\alpha}(m_\alpha(s_\alpha-\rho_\alpha+1)).\qedhere
\]
\end{proof}

For the rest of this subsection, we fix a big line bundle $L=\sum_{\alpha\in\mathcal{A}}\lambda_\alpha D_\alpha$ on $X$.
We note that $\lambda_\alpha>0$ for all $\alpha\in\mathcal{A}$ by Remark~\ref{Eff}~(1).
In the following, we study the rightmost pole of the function $s\mapsto\widehat{\mathsf{H}}_{\mathbf{m}}(\mathbf{0},sL)$.
The function $\zeta_{F_\alpha}(m_\alpha(s\lambda_\alpha-\rho_\alpha+1))$ has a unique pole when $m_\alpha(s\lambda_\alpha-\rho_\alpha+1)=1$, that is,
\[
s=\frac{\rho_\alpha-1+m_\alpha^{-1}}{\lambda_\alpha}.
\]

We recall the definitions of the $a$-invariant and the $b$-invariant of $(X,D_{\mathbf{m}})$ with respect to $L$ defined in \cite[\S~9]{PSTVA21}.

\begin{dfn}\label{ab}
\leavevmode
\begin{enumerate}
\item We define the \textbf{$a$-invariant} $a((X,D_{\mathbf{m}}),L)$ of $(X,D_{\mathbf{m}})$ with respect to $L$ by
\[
a=a((X,D_{\mathbf{m}}),L):=\max_{\alpha\in\mathcal{A}}\frac{\rho_\alpha-1+m_\alpha^{-1}}{\lambda_\alpha}.
\]
\item We define a subset $\mathcal{A}_{\mathbf{m}}(L)$ of $\mathcal{A}$ by
\[
\mathcal{A}_{\mathbf{m}}(L)=\left\{\alpha\in\mathcal{A}\,\middle|\,\frac{\rho_\alpha-1+m_\alpha^{-1}}{\lambda_\alpha}=a((X,D_{\mathbf{m}}),L)\right\}.
\]
\item We define the \textbf{$b$-invariant} $b(F,(X,D_{\mathbf{m}}),L)$ of $(X,D_{\mathbf{m}})$ with respect to $L$ by
\[
b=b(F,(X,D_{\mathbf{m}}),L)=\#\mathcal{A}_{\mathbf{m}}(L).
\]
\item We fix a positive constant $\delta'$ such that
\[
0<\delta'<\min\left\{\left\{a-\frac{\rho_\alpha-1+m_\alpha^{-1}}{\lambda_\alpha}\,\middle|\,\alpha\in\mathcal{A}\setminus\mathcal{A}_{\mathbf{m}}(L)\right\}\cup\left\{\frac{\delta}{\lambda_\alpha}\,\middle|\,\alpha\in\mathcal{A}\right\}\right\}.
\]
\end{enumerate}
\end{dfn}

\begin{rem}\label{remab}
We note that $sL\in\mathsf{T}_{>-\delta}$ if $\Re s>a-\delta'$ and that the function
\[
s\mapsto\prod_{\alpha\in\mathcal{A}}\zeta_{F_\alpha}(m_\alpha(s\lambda_\alpha-\rho_\alpha+1))
\]
has a rightmost pole at $s=a$, and its order is equal to $b$.
Hence, the function
\[
\widehat{\mathsf{H}}_{\mathbf{m}}(\mathbf{0},sL)=f^{\mathbf{0}}(sL)\prod_{\alpha\in\mathcal{A}}\zeta_{F_\alpha}(m_\alpha(s\lambda_\alpha-\rho_\alpha+1))
\]
is holomorphic on the half-plane $\Re s>a$, and if the function $\widehat{\mathsf{H}}_{\mathbf{m}}(\mathbf{0},sL)$ has a pole at $s=a$, then this pole is a rightmost pole and its order is at most $b$.
\end{rem}

The following proposition is an analogue of \cite[LEMMA~9.3]{PSTVA21}.

\begin{prop}\label{c}
The function $\widehat{\mathsf{H}}_{\mathbf{m}}(\mathbf{0},sL)$ has a rightmost pole at $s=a$ of order $b$.
Moreover, the constant
\[
c:=\lim_{s\to a}(s-a)^b\widehat{\mathsf{H}}_{\mathbf{m}}(\mathbf{0},sL)
\]
is positive.
\end{prop}

\begin{proof}
By Corollary~\ref{fzero}, we have
\[
\widehat{\mathsf{H}}_{\mathbf{m}}(\mathbf{0},sL)=f^{\mathbf{0}}(sL)\prod_{\alpha\in\mathcal{A}}\zeta_{F_\alpha}(m_\alpha(s\lambda_\alpha-\rho_\alpha+1))
\]
as meromorphic functions on the half-plane $\Re s>a-\delta'$.
By the definitions of $a$ and $b$, we have
\begin{align*}
&\lim_{s\to a}(s-a)^b\prod_{\alpha\in\mathcal{A}}\zeta_{F_\alpha}(m_\alpha(s\lambda_\alpha-\rho_\alpha+1))
\\
&=\left(\prod_{\alpha\in\mathcal{A}_{\mathbf{m}}(L)}\lim_{s\to a}(s-a)\zeta_{F_\alpha}(m_\alpha(s\lambda_\alpha-\rho_\alpha+1))\right)\left(\prod_{\alpha\in\mathcal{A}\setminus\mathcal{A}_{\mathbf{m}}(L)}\zeta_{F_\alpha}(m_\alpha(a\lambda_\alpha-\rho_\alpha+1))\right).
\end{align*}
Since the residue of the Dedekind zeta function at 1 is positive, for each $\alpha\in\mathcal{A}_{\mathbf{m}}(L)$, we have
\[
\lim_{s\to a}(s-a)\zeta_{F_\alpha}(m_\alpha(s\lambda_\alpha-\rho_\alpha+1))=\frac{1}{m_\alpha\lambda_\alpha}\Res_{z=1}\zeta_{F_\alpha}(z)>0.
\]
Moreover, we have
\[
m_\alpha(a\lambda_\alpha-\rho_\alpha+1)>1,\quad\zeta_{F_\alpha}(m_\alpha(a\lambda_\alpha-\rho_\alpha+1))>0
\]
for each $\alpha\in\mathcal{A}\setminus\mathcal{A}_{\mathbf{m}}(L)$.
Hence, we have
\[
\lim_{s\to a}(s-a)^b\prod_{\alpha\in\mathcal{A}}\zeta_{F_\alpha}(m_\alpha(s\lambda_\alpha-\rho_\alpha+1))>0.
\]
Thus, it suffices to show that the following limit is positive:
\begin{align*}
\lim_{t\to a}f^{\mathbf{0}}(tL)&=\left(\prod_{v\in\Omega_F^\infty}\widehat{\mathsf{H}}_{\mathbf{m},v}(\mathbf{0},aL)\right)
\\
&\quad\times\left(\prod_{\substack{v\in\Omega_F^{<\infty}\\ q_v\leq N}}\widehat{\mathsf{H}}_{\mathbf{m},v}(\mathbf{0},aL)\prod_{\alpha\in\mathcal{A}}\prod_{w\mid v}\zeta_{(F_\alpha)_w}(m_\alpha(a\lambda_\alpha-\rho_\alpha+1))^{-1}\right)
\\
&\quad\times\left(\lim_{t\to a}\prod_{\substack{v\in\Omega_F^{<\infty}\\ q_v>N}}\widehat{\mathsf{H}}_{\mathbf{m},v}(\mathbf{0},tL)\prod_{\alpha\in\mathcal{A}}\prod_{w\mid v}\zeta_{(F_\alpha)_w}(m_\alpha(t\lambda_\alpha-\rho_\alpha+1))^{-1}\right).
\end{align*}
Here, the limit is taken as $t$ tends to $a$ through real values.
For each $\alpha\in\mathcal{A}$, $\mathsf{f}_\alpha$ is a section of $\mathcal{O}(D_\alpha)$.
Then we have $\|\mathsf{f}_\alpha(\mathbf{x}_v)\|_v>0$ for all $\alpha\in\mathcal{A}$, $v\in\Omega_F$, and $\mathbf{x}_v\in G(F_v)$.
Moreover, $\delta_{\mathbf{m},v}$ is equal to $1$ on a non-empty open subset of $G(F_v)$.
It follows that
\[
\widehat{\mathsf{H}}_{\mathbf{m},v}(\mathbf{0},tL)=\int_{G(F_v)}\delta_{\mathbf{m},v}(\mathbf{x}_v)\prod_{\alpha\in\mathcal{A}}\|\mathsf{f}_\alpha(\mathbf{x}_v)\|_v^{t\lambda_\alpha}\,\mathrm{d}\mathbf{x}_v>0
\]
for each real number $t>a-\delta'$ and $v\in\Omega_F$.
Hence, we have
\begin{align*}
0&<\prod_{v\in\Omega_F^\infty}\widehat{\mathsf{H}}_{\mathbf{m},v}(\mathbf{0},aL),
\\
0&<\prod_{\substack{v\in\Omega_F^{<\infty}\\ q_v\leq N}}\widehat{\mathsf{H}}_{\mathbf{m},v}(\mathbf{0},aL)\prod_{\alpha\in\mathcal{A}}\prod_{w\mid v}\zeta_{(F_\alpha)_w}(m_\alpha(a\lambda_\alpha-\rho_\alpha+1))^{-1}.
\end{align*}
Moreover, by Proposition~\ref{Hvzero2}, the limit
\[
\lim_{t\to a}\prod_{\substack{v\in\Omega_F^{<\infty}\\ q_v>N}}\widehat{\mathsf{H}}_{\mathbf{m},v}(\mathbf{0},tL)\prod_{\alpha\in\mathcal{A}}\prod_{w\mid v}\zeta_{(F_\alpha)_w}(m_\alpha(t\lambda_\alpha-\rho_\alpha+1))^{-1}
\]
exists and is non-zero.
Since
\[
\widehat{\mathsf{H}}_{\mathbf{m},v}(\mathbf{0},tL)\prod_{\alpha\in\mathcal{A}}\prod_{w\mid v}\zeta_{(F_\alpha)_w}(m_\alpha(t\lambda_\alpha-\rho_\alpha+1))^{-1}>0
\]
for each $v\in\Omega_F^{<\infty}$, we have
\[
\lim_{t\to a}\prod_{\substack{v\in\Omega_F^{<\infty}\\ q_v>N}}\widehat{\mathsf{H}}_{\mathbf{m},v}(\mathbf{0},tL)\prod_{\alpha\in\mathcal{A}}\prod_{w\mid v}\zeta_{(F_\alpha)_w}(m_\alpha(t\lambda_\alpha-\rho_\alpha+1))^{-1}>0.
\]
This completes the proof.
\end{proof}

\subsection{The height integrals when $\mathbf{a}\neq\mathbf{0}$}
We continue to use the notation fixed in the previous subsection.
In this subsection, for $\mathbf{a}\in G(F)\setminus\{\mathbf{0}\}$, we study the local height integral
\[
\widehat{\mathsf{H}}_{\mathbf{m},v}(\mathbf{a},\mathbf{s})=\int_{G(F_v)}\delta_{\mathbf{m},v}(\mathbf{x}_v)\mathsf{H}_v(\mathbf{x}_v,-\mathbf{s})\psi_{\mathbf{a},v}(\mathbf{x}_v)\,\mathrm{d}\mathbf{x}_v
\]
at $v\in\Omega_F$ and its Euler product
\[
\widehat{\mathsf{H}}_{\mathbf{m}}(\mathbf{a},\mathbf{s})=\prod_{v\in\Omega_F}\widehat{\mathsf{H}}_{\mathbf{m},v}(\mathbf{a}_v,\mathbf{s})
\]
following \cite[\S~8]{PSTVA21}.
We note that $\widehat{\mathsf{H}}_{\mathbf{m}}(\mathbf{a},\mathbf{s})=0$ when $\mathbf{a}\notin\Lambda_X$ by Proposition~\ref{Hzero}.
Thus, in what follows, it suffices to consider $\mathbf{a}\in\Lambda_X\setminus\{\mathbf{0}\}$.

We recall some results on the integrals of characters over $\mathcal{O}_v^\times$, and introduce some notation.
The following statement is proved in \cite[Lemma~10.3]{CLT02}.

\begin{rem}\label{psi}
Let $v$ be a finite place of $F$ such that $q_v>N$ and let $u\in\mathcal{O}_v^\times$.
\begin{enumerate}
\item We have
\[
\int_{\mathcal{O}_v^\times}\psi_v\left(\frac{ux_v}{\pi_v}\right)\,\mathrm{d}x_v=-\frac{1}{q_v}.
\]
\item Let $i$ and $d$ be non-negative integers such that $d\neq0$ and $(i,d)\neq(1,1)$.
Then we have
\[
\int_{\mathcal{O}_v^\times}\psi_v\left(\frac{ux_v^d}{\pi_v^{id}}\right)\,\mathrm{d}x_v=0.
\]
\end{enumerate}
\end{rem}

The following notation is based on \cite[before Proposition~5.5]{CLT02}.
\begin{nota}
For each $\mathbf{a}\in\Lambda_X\setminus\{\mathbf{0}\}$, we define a subset $\mathcal{A}^0(\mathbf{a})$ of $\mathcal{A}$ by
\[
\mathcal{A}^0(\mathbf{a})=\{\alpha\in\mathcal{A}\mid d_\alpha(f_{\mathbf{a}})=0\}.
\]
\end{nota}

First, we consider the case $j_v(\mathbf{a})=0$.
The following proposition is an analogue of \cite[PROPOSITION~8.1]{PSTVA21}.

\begin{prop}\label{Hva}
We have
\[
\widehat{\mathsf{H}}_{\mathbf{m},v}(\mathbf{a},\mathbf{s})\prod_{\alpha\in\mathcal{A}^0(\mathbf{a})}\prod_{\beta\in\mathcal{A}_v(\alpha)}\zeta_{F_{v,\beta}}(m_\alpha(s_\alpha-\rho_\alpha+1))^{-1}=1+O(q_v^{-1-\delta})
\]
uniformly for $\mathbf{a}\in\Lambda_X\setminus\{\mathbf{0}\}$, $\mathbf{s}\in\mathsf{T}_{>-\delta}$, and $v\in\Omega_F^{<\infty}$ such that $q_v>N$ and $j_v(\mathbf{a})=0$.
\end{prop}

\begin{proof}
Let $\mathbf{a}=(a_1,\dots,a_n)\in\Lambda_X\setminus\{\mathbf{0}\}$, let $\mathbf{s}\in\mathsf{T}_{>-\delta}$, and let $v\in\Omega_F^{<\infty}$ such that $q_v>N$ and $j_v(\mathbf{a})=0$.
We recall that $q_v=\#(\mathcal{O}_v/\mathfrak{m}_v)$ and
\[
j_v(\mathbf{a})=\min\{v(a_1),\dots,v(a_n)\},
\]
and we use the measures $\mathrm{d}\mathbf{x}_v$ and $\mathrm{d}\tau_{X,v}$ fixed in the previous subsection.
In particular, on $G(F_v)$, we have
\[
\mathrm{d}\mathbf{x}_v=\mathsf{H}_v(\mathbf{x}_v,\boldsymbol{\rho})\,\mathrm{d}\tau_{X,v}.
\]
By an argument similar to that in the proof of Proposition~\ref{Hvzero}, we have
\[
\widehat{\mathsf{H}}_{\mathbf{m},v}(\mathbf{a},\mathbf{s})=\sum_{B\subseteq\mathcal{A}_v}\sum_{y\in\mathcal{D}_{v,B}^\circ(k_v)}\int_{\eta_v^{-1}(y)}\delta_{\mathbf{m},v}(\mathbf{x}_v)\mathsf{H}_v(\mathbf{x}_v,\mathbf{s}-\boldsymbol{\rho})^{-1}\psi_{\mathbf{a},v}(\mathbf{x}_v)\,\mathrm{d}\tau_{X,v}.
\]
For each subset $B=\{\beta_1,\dots,\beta_r\}$ of $\mathcal{A}_v$, we set
\[
A_B^{\mathbf{a}}(\mathbf{s}):=\sum_{y\in\mathcal{D}_{v,B}^\circ(k_v)}\int_{\eta_v^{-1}(y)}\delta_{\mathbf{m},v}(\mathbf{x}_v)\mathsf{H}_v(\mathbf{x}_v,\mathbf{s}-\boldsymbol{\rho})^{-1}\psi_{\mathbf{a},v}(\mathbf{x}_v)\,\mathrm{d}\tau_{X,v}.
\]
We consider $A_B^{\mathbf{a}}(\mathbf{s})$ in several cases.
\begin{enumerate}
\item Case $B=\varnothing$.

We note that $\mathbf{a}\in G(\mathcal{O}_v)$ since $j_v(\mathbf{a})=0$.
If $\mathbf{x}_v\in\eta_v^{-1}(\mathcal{D}_{v,\varnothing}^\circ(k_v))=G(\mathcal{O}_v)$, then we have $f_{\mathbf{a}}(\mathbf{x}_v)\in\mathcal{O}_v$ and $\psi_{\mathbf{a},v}(\mathbf{x}_v)=\psi_v(f_{\mathbf{a}}(\mathbf{x}_v))=1$.
Therefore, we have
\[
A_{\varnothing}^{\mathbf{a}}(\mathbf{s})=\sum_{y\in\mathcal{D}_{v,\varnothing}^\circ(k_v)}\int_{\eta_v^{-1}(y)}\delta_{\mathbf{m},v}(\mathbf{x}_v)\mathsf{H}_v(\mathbf{x}_v,\mathbf{s}-\boldsymbol{\rho})^{-1}\,\mathrm{d}\tau_{X,v}=A_{\varnothing}(\mathbf{s})=1.
\]
The last equality follows from case (1) of the proof of Proposition~\ref{Hvzero2}.
\item Case $r=1$ and $f_{v,\beta_1}\neq1$.

As in case (2) of the proof of Proposition~\ref{Hvzero2}, we have $\#\mathcal{D}_{v,B}^\circ(k_v)=0$.
It follows that $A_B^{\mathbf{a}}(\mathbf{s})=0$.

\item Case $r=1$, $f_{v,\beta_1}=1$, and $d_{\alpha(\beta_1)}(f_{\mathbf{a}})=0$.

As in \cite[\S~8.1.1]{PSTVA21}, we have $\psi_{\mathbf{a},v}(\mathbf{x}_v)=1$.
Hence, we obtain
\begin{align*}
A_B^{\mathbf{a}}(\mathbf{s})&=\sum_{y\in\mathcal{D}_{v,B}^\circ(k_v)}\int_{\eta_v^{-1}(y)}\delta_{\mathbf{m},v}(\mathbf{x}_v)\mathsf{H}_v(\mathbf{x}_v,\mathbf{s}-\boldsymbol{\rho})^{-1}\,\mathrm{d}\tau_{X,v}
\\
&=A_B(\mathbf{s})
\\
&=q_v^{-m_{\alpha(\beta_1)}(s_{\alpha(\beta_1)}-\rho_{\alpha(\beta_1)}+1)}+O(q_v^{-1-\delta}).
\end{align*}
The last estimate follows from case (3) of the proof of Proposition~\ref{Hvzero2}.
\item Case $r=1$, $f_{v,\beta_1}=1$, and $d_{\alpha(\beta_1)}(f_{\mathbf{a}})\neq0$.

First, we calculate the integral
\[
\int_{\eta_v^{-1}(y)}\delta_{\mathbf{m},v}(\mathbf{x}_v)\mathsf{H}_v(\mathbf{x}_v,\mathbf{s}-\boldsymbol{\rho})^{-1}\psi_{\mathbf{a},v}(\mathbf{x}_v)\,\mathrm{d}\tau_{X,v}
\]
for each $y\in\mathcal{D}_{v,B}^\circ(k_v)$.
We recall that $(z_1,\dots,z_n)\colon\eta_v^{-1}(y)\to F_v^n$ are the local analytic coordinates as in Remark~\ref{Sal} and $\mathrm{d}\mathbf{z}=\prod_{i=1}^n\mathrm{d}z_i$ is the measure on $\mathfrak{m}_v^n$.
By Remark~\ref{Hnp}, we have
\begin{align*}
\delta_{\mathbf{m},v}(\mathbf{x}_v)&=\left\{
\begin{array}{ll}
1&\text{if }v(z_1)\in m_{\alpha(\beta_1)}\mathbb{Z},\\
0&\text{if }v(z_1)\notin m_{\alpha(\beta_1)}\mathbb{Z},
\end{array}\right. \\
\mathsf{H}_v(\mathbf{x}_v,\mathbf{s}-\boldsymbol{\rho})^{-1}&=|z_1(\mathbf{x}_v)|_v^{s_{\alpha(\beta_1)}-\rho_{\alpha(\beta_1)}},
\end{align*}
for all $\mathbf{x}_v\in\eta_v^{-1}(y)\cap G(F_v)$.
We note that $\delta_{\mathbf{m},v}(\mathbf{x}_v)=0$ for $\mathbf{x}_v\in X(F_v)\setminus G(F_v)$.
We divide the argument into cases in order to evaluate the character factor.
\begin{enumerate}
\item Case $y\in(\mathcal{D}_{v,B}^\circ\setminus\mathcal{E}(f_{\mathbf{a}}))(k_v)$.

By Remark~\ref{Hnp}~(3), there exists $u\in\mathcal{O}_v^\times$ such that, for all $\mathbf{x}_v\in\eta_v^{-1}(y)\cap G(F_v)$,
\[
\psi_{\mathbf{a},v}(\mathbf{x}_v)=\psi_v\left(\frac{u\pi_v^{j_v(\mathbf{a})}}{z_1^{d_{\alpha(\beta_1)}(f_{\mathbf{a}})}}\right)=\psi_v\left(\frac{u}{z_1^{d_{\alpha(\beta_1)}(f_{\mathbf{a}})}}\right).
\]
Hence, we have
\begin{align*}
&\int_{\eta_v^{-1}(y)}\delta_{\mathbf{m},v}(\mathbf{x}_v)\mathsf{H}_v(\mathbf{x}_v,\mathbf{s}-\boldsymbol{\rho})^{-1}\psi_{\mathbf{a},v}(\mathbf{x}_v)\,\mathrm{d}\tau_{X,v}
\\
&=\int_{\mathfrak{m}_v^n}\mathbf{1}_{\coprod_{k=1}^\infty\pi_v^{km_{\alpha(\beta_1)}}\mathcal{O}_v^\times}(z_1)|z_1|_v^{s_{\alpha(\beta_1)}-\rho_{\alpha(\beta_1)}}\psi_v\left(\frac{u}{z_1^{d_{\alpha(\beta_1)}(f_{\mathbf{a}})}}\right)\,\mathrm{d}\mathbf{z}
\\
&=\left(\sum_{k=1}^\infty\int_{\pi_v^{km_{\alpha(\beta_1)}}\mathcal{O}_v^\times}|\pi_v^{km_{\alpha(\beta_1)}}|_v^{s_{\alpha(\beta_1)}-\rho_{\alpha(\beta_1)}}\psi_v\left(\frac{u}{z_1^{d_{\alpha(\beta_1)}(f_{\mathbf{a}})}}\right)\,\mathrm{d}z_1\right)\left(\prod_{i=2}^n\int_{\mathfrak{m}_v}1\,\mathrm{d}z_i\right)
\\
&=\frac{1}{q_v^{n-1}}\sum_{k=1}^\infty q_v^{-km_{\alpha(\beta_1)}(s_{\alpha(\beta_1)}-\rho_{\alpha(\beta_1)})}\int_{\pi_v^{km_{\alpha(\beta_1)}}\mathcal{O}_v^\times}\psi_v\left(\frac{u}{z_1^{d_{\alpha(\beta_1)}(f_{\mathbf{a}})}}\right)\,\mathrm{d}z_1.
\end{align*}
For each $k$, we have
\begin{align*}
&\int_{\pi_v^{km_{\alpha(\beta_1)}}\mathcal{O}_v^\times}\psi_v\left(\frac{u}{z_1^{d_{\alpha(\beta_1)}(f_{\mathbf{a}})}}\right)\,\mathrm{d}z_1
\\
&=q_v^{-km_{\alpha(\beta_1)}}\int_{\mathcal{O}_v^\times}\psi_v\left(\frac{u}{\pi_v^{km_{\alpha(\beta_1)}d_{\alpha(\beta_1)}(f_{\mathbf{a}})}z_1^{d_{\alpha(\beta_1)}(f_{\mathbf{a}})}}\right)\,\mathrm{d}z_1
\\
&=q_v^{-km_{\alpha(\beta_1)}}\int_{\mathcal{O}_v^\times}\psi_v\left(\frac{uz_1^{d_{\alpha(\beta_1)}(f_{\mathbf{a}})}}{\pi_v^{km_{\alpha(\beta_1)}d_{\alpha(\beta_1)}(f_{\mathbf{a}})}}\right)\,\mathrm{d}z_1
\\
&=\left\{\begin{array}{ll}
-q_v^{-km_{\alpha(\beta_1)}-1}=-q_v^{-2}&\text{if }k=m_{\alpha(\beta_1)}=d_{\alpha(\beta_1)}(f_{\mathbf{a}})=1,
\\
0&\text{otherwise},
\end{array}\right.,
\end{align*}
where the first and second equalities follow from the changes of variables $z_1\mapsto\pi_v^{km_{\alpha(\beta_1)}}z_1$ and $z_1\mapsto z_1^{-1}$, respectively, and the third equality is obtained by applying Remark~\ref{psi} to $(km_{\alpha(\beta_1)},d_{\alpha(\beta_1)}(f_{\mathbf{a}}))$.
Hence, we have
\[
\int_{\eta_v^{-1}(y)}\delta_{\mathbf{m},v}(\mathbf{x}_v)\mathsf{H}_v(\mathbf{x}_v,\mathbf{s}-\boldsymbol{\rho})^{-1}\psi_{\mathbf{a},v}(\mathbf{x}_v)\,\mathrm{d}\tau_{X,v}=0
\]
if $(m_{\alpha(\beta_1)},d_{\alpha(\beta_1)}(f_{\mathbf{a}}))\neq(1,1)$.
If $(m_{\alpha(\beta_1)},d_{\alpha(\beta_1)}(f_{\mathbf{a}}))=(1,1)$, then we have
\begin{align*}
&\left|\int_{\eta_v^{-1}(y)}\delta_{\mathbf{m},v}(\mathbf{x}_v)\mathsf{H}_v(\mathbf{x}_v,\mathbf{s}-\boldsymbol{\rho})^{-1}\psi_{\mathbf{a},v}(\mathbf{x}_v)\,\mathrm{d}\tau_{X,v}\right|
\\
&=\left|\frac{1}{q_v^{n-1}}q_v^{-s_{\alpha(\beta_1)}+\rho_{\alpha(\beta_1)}}q_v^{-2}\right|
\\
&=q_v^{-n-1-\Re s_{\alpha(\beta_1)}+\rho_{\alpha(\beta_1)}}
\\
&\leq q_v^{-n-1+\delta}\leq q_v^{-n-2/3},
\end{align*}
where the first inequality follows from the facts that $m_{\alpha(\beta_1)}=1$ and $\mathbf{s}\in\mathsf{T}_{>-\delta}$, and the second inequality follows from the fact that $\delta<1/3$.
\item Case $y\in(\mathcal{D}_{v,B}^\circ\cap\mathcal{E}(f_{\mathbf{a}}))(k_v)$.

We have
\begin{align*}
&\left|\int_{\eta_v^{-1}(y)}\delta_{\mathbf{m},v}(\mathbf{x}_v)\mathsf{H}_v(\mathbf{x}_v,\mathbf{s}-\boldsymbol{\rho})^{-1}\psi_{\mathbf{a},v}(\mathbf{x}_v)\,\mathrm{d}\tau_{X,v}\right|
\\
&\leq\int_{\eta_v^{-1}(y)}|\delta_{\mathbf{m},v}(\mathbf{x}_v)\mathsf{H}_v(\mathbf{x}_v,\mathbf{s}-\boldsymbol{\rho})^{-1}\psi_{\mathbf{a},v}(\mathbf{x}_v)|\,\mathrm{d}\tau_{X,v}
\\
&=\int_{\eta_v^{-1}(y)}\delta_{\mathbf{m},v}(\mathbf{x}_v)\mathsf{H}_v(\mathbf{x}_v,(\Re s_\alpha)_{\alpha\in\mathcal{A}}-\boldsymbol{\rho})^{-1}\,\mathrm{d}\tau_{X,v}
\\
&=\frac{1}{q_v^{n-1}}\left(1-\frac{1}{q_v}\right)\frac{q_v^{-m_{\alpha(\beta_1)}(\Re s_{\alpha(\beta_1)}-\rho_{\alpha(\beta_1)}+1)}}{1-q_v^{-m_{\alpha(\beta_1)}(\Re s_{\alpha(\beta_1)}-\rho_{\alpha(\beta_1)}+1)}}
\\
&\leq\frac{2q_v^{m_{\alpha(\beta_1)}\delta-1}}{q_v^{n-1}}<2q_v^{-n+1/3},
\end{align*}
where the second equality follows from \eqref{2} in the proof of Proposition~\ref{Hvzero}, and the third inequality follows from the fact that $\mathbf{s}\in\mathsf{T}_{>-\delta}$ and $\delta<(3m_{\alpha(\beta_1)})^{-1}$.
\end{enumerate}
Therefore, we have
\begin{align*}
|A_B^{\mathbf{a}}(\mathbf{s})|&\leq\sum_{y\in(\mathcal{D}_{v,B}^\circ\setminus\mathcal{E}(f_{\mathbf{a}}))(k_v)}\left|\int_{\eta_v^{-1}(y)}\delta_{\mathbf{m},v}(\mathbf{x}_v)\mathsf{H}_v(\mathbf{x}_v,\mathbf{s}-\boldsymbol{\rho})^{-1}\psi_{\mathbf{a},v}(\mathbf{x}_v)\,\mathrm{d}\tau_{X,v}\right|
\\
&+\sum_{y\in(\mathcal{D}_{v,B}^\circ\cap\mathcal{E}(f_{\mathbf{a}}))(k_v)}\left|\int_{\eta_v^{-1}(y)}\delta_{\mathbf{m},v}(\mathbf{x}_v)\mathsf{H}_v(\mathbf{x}_v,\mathbf{s}-\boldsymbol{\rho})^{-1}\psi_{\mathbf{a},v}(\mathbf{x}_v)\,\mathrm{d}\tau_{X,v}\right|
\\
&\leq\#\mathcal{D}_{v,B}(k_v)q_v^{-n-2/3}+\#(\mathcal{D}_{v,B}\cap\mathcal{E}(f_{\mathbf{a}}))(k_v)q_v^{-n+1/3}.
\end{align*}
By the Lang--Weil estimates, uniformly for $\mathbf{a}$, we have
\[
\#\mathcal{D}_{v,B}(k_v)=O(q_v^{n-1}),\quad\#(\mathcal{D}_{v,B}\cap\mathcal{E}(f_{\mathbf{a}}))(k_v)=O(q_v^{n-2}).
\]
Hence, we have $|A_B^{\mathbf{a}}(\mathbf{s})|=O(q_v^{-5/3})$.
\item Case $r\geq2$.

By case (4) of the proof of Proposition~\ref{Hvzero2}, we have
\begin{align*}
|A_B^{\mathbf{a}}(\mathbf{s})|&\leq\sum_{y\in\mathcal{D}_{v,B}^\circ(k_v)}\left|\int_{\eta_v^{-1}(y)}\delta_{\mathbf{m},v}(\mathbf{x}_v)\mathsf{H}_v(\mathbf{x}_v,\mathbf{s}-\boldsymbol{\rho})^{-1}\psi_{\mathbf{a},v}(\mathbf{x}_v)\,\mathrm{d}\tau_{X,v}\right|
\\
&\leq\sum_{y\in\mathcal{D}_{v,B}^\circ(k_v)}\int_{\eta_v^{-1}(y)}\delta_{\mathbf{m},v}(\mathbf{x}_v)\mathsf{H}_v(\mathbf{x}_v,(\Re s_\alpha)_{\alpha\in\mathcal{A}}-\boldsymbol{\rho})^{-1}\,\mathrm{d}\tau_{X,v}
\\
&=A_B((\Re s_\alpha)_{\alpha\in\mathcal{A}})=O(q_v^{-1-\delta}).
\end{align*}
\end{enumerate}

Collecting these estimates, we obtain
\begin{equation}\label{5}
\widehat{\mathsf{H}}_{\mathbf{m},v}(\mathbf{a},\mathbf{s})=1+\sum_{\substack{\beta\in\mathcal{A}_v,f_{v,\beta}=1,\\ \alpha(\beta)\in\mathcal{A}^0(\mathbf{a})}}q_v^{-m_{\alpha(\beta)}(s_{\alpha(\beta)}-\rho_{\alpha(\beta)}+1)}+O(q_v^{-1-\delta}).
\end{equation}
On the other hand, by expanding the local zeta factors, we have
\begin{equation}\label{6}
\prod_{\alpha\in\mathcal{A}^0(\mathbf{a})}\prod_{\beta\in\mathcal{A}_v(\alpha)}\zeta_{F_{v,\beta}}(m_\alpha(s_\alpha-\rho_\alpha+1))^{-1}=1-\sum_{\substack{\beta\in\mathcal{A}_v,f_{v,\beta}=1,\\ \alpha(\beta)\in\mathcal{A}^0(\mathbf{a})}}q_v^{-m_{\alpha(\beta)}(s_{\alpha(\beta)}-\rho_{\alpha(\beta)}+1)}+O(q_v^{-1-\delta}),
\end{equation}
as $q_v\to\infty$.
Combining the asymptotic formulae \eqref{5} and \eqref{6} completes the proof.
\end{proof}

Next, we bound the local height integral $\widehat{\mathsf{H}}_{\mathbf{m},v}(\mathbf{a},\mathbf{s})$ when $j_v(\mathbf{a})\neq0$.

\begin{prop}\label{Hvajneq0}
We have
\[
\widehat{\mathsf{H}}_{\mathbf{m},v}(\mathbf{a},\mathbf{s})\prod_{\alpha\in\mathcal{A}^0(\mathbf{a})}\prod_{\beta\in\mathcal{A}_v(\alpha)}\zeta_{F_{v,\beta}}(m_\alpha(s_\alpha-\rho_\alpha+1))^{-1}=O(1)
\]
uniformly for $\mathbf{a}\in\Lambda_X\setminus\{\mathbf{0}\}$, $\mathbf{s}\in\mathsf{T}_{>-\delta}$, and $v\in\Omega_F^{<\infty}$ such that $q_v>N$ and $j_v(\mathbf{a})\neq0$.
\end{prop}

\begin{proof}
Since $\mathbf{s}\in\mathsf{T}_{>-\delta}$ and $\delta<(3\max_{\alpha\in\mathcal{A}}m_\alpha)^{-1}$, we have
\begin{align}
&\left|\prod_{\alpha\in\mathcal{A}^0(\mathbf{a})}\prod_{\beta\in\mathcal{A}_v(\alpha)}\zeta_{F_{v,\beta}}(m_\alpha(s_\alpha-\rho_\alpha+1))^{-1}\right|\notag
\\
&\leq\prod_{\alpha\in\mathcal{A}^0(\mathbf{a})}\prod_{\beta\in\mathcal{A}_v(\alpha)}|1-q_v^{-f_{v,\beta}m_\alpha(s_\alpha-\rho_\alpha+1)}|\notag
\\
&\leq\prod_{\alpha\in\mathcal{A}^0(\mathbf{a})}\prod_{\beta\in\mathcal{A}_v(\alpha)}(1+q_v^{f_{v,\beta}(m_{\alpha}\delta-1)})\leq 2^{\#\overline{\mathcal{A}}},\label{7}
\end{align}
where $\overline{\mathcal{A}}$ is the set of geometrically irreducible components of $\bigcup_{\alpha\in\mathcal{A}}D_\alpha$.
By Proposition~\ref{Hvzero}, we have
\begin{align*}
|\widehat{\mathsf{H}}_{\mathbf{m},v}(\mathbf{a},\mathbf{s})|&=\left|\int_{G(F_v)}\delta_{\mathbf{m},v}(\mathbf{x}_v)\mathsf{H}_v(\mathbf{x}_v,-\mathbf{s})\psi_{\mathbf{a},v}(\mathbf{x}_v)\,\mathrm{d}\mathbf{x}_v\right|
\\
&\leq\int_{G(F_v)}\delta_{\mathbf{m},v}(\mathbf{x}_v)\mathsf{H}_v(\mathbf{x}_v,(-\Re s_{\alpha})_{\alpha\in\mathcal{A}})\,\mathrm{d}\mathbf{x}_v
\\
&=\widehat{\mathsf{H}}_{\mathbf{m},v}(\mathbf{0},(\Re s_\alpha)_{\alpha\in\mathcal{A}})
\\
&=\sum_{B\subseteq\mathcal{A}_v}\frac{\#\mathcal{D}_{v,B}^\circ(k_v)}{q_v^{n-\#B}}\prod_{\beta\in B}\left(1-\frac{1}{q_v}\right)\frac{1}{q_v^{m_{\alpha(\beta)}(\Re s_{\alpha(\beta)}-\rho_{\alpha(\beta)}+1)}-1}
\\
&<\sum_{B\subseteq\mathcal{A}_v}\frac{\#\mathcal{D}_{v,B}^\circ(k_v)}{q_v^{n-\#B}}\prod_{\beta\in B}\frac{1}{q_v^{2/3}-1}
\\
&\leq\sum_{B\subseteq\mathcal{A}_v}\frac{\#\mathcal{D}_{v,B}^\circ(k_v)}{q_v^{n-\#B}}\prod_{\beta\in B}\frac{1}{2^{2/3}-1}<4^{\#\overline{\mathcal{A}}}\sum_{B\subseteq\mathcal{A}_v}\frac{\#\mathcal{D}_{v,B}^\circ(k_v)}{q_v^{n-\#B}}.
\end{align*}
By the Lang--Weil estimates, we have
\[
\frac{\#\mathcal{D}_{v,B}^\circ(k_v)}{q_v^{n-\#B}}=O(1).
\]
Hence, we obtain
\begin{equation}\label{8}
\widehat{\mathsf{H}}_{\mathbf{m},v}(\mathbf{a},\mathbf{s})=O(1)
\end{equation}
uniformly for $\mathbf{a}\in\Lambda_X\setminus\{\mathbf{0}\}$, $\mathbf{s}\in\mathsf{T}_{>-\delta}$, and $v\in\Omega_F^{<\infty}$ such that $q_v>N$ and $j_v(\mathbf{a})\neq0$.
Combining \eqref{7} and \eqref{8} completes the proof.
\end{proof}

We evaluate the local height integrals $\widehat{\mathsf{H}}_{\mathbf{m},v}(\mathbf{a},\mathbf{s})$ at finite places $v$ with $q_v\leq N$ and at infinite places $v$.
We first treat the places $v\notin S$.

\begin{nota}
For each $\mathbf{a}=(a_1,\dots,a_n)\in G(F)\setminus\{\mathbf{0}\}$ and $v\in\Omega_F$, we set
\[
H_v(\mathbf{a})=\max\{|a_1|_v,\dots,|a_n|_v\}.
\]
We fix a norm $\|\cdot\|$ on the real vector space $G(F\otimes_{\mathbb{Q}}\mathbb{R})$.
We also denote the composition of $\|\cdot\|\colon G(F\otimes_{\mathbb{Q}}\mathbb{R})\to\mathbb{R}$ and the diagonal embedding $G(F)\to G(F\otimes_{\mathbb{Q}}\mathbb{R})$ by $\|\cdot\|$.
For each $v\in\Omega_F^\infty$, we fix a positive constant $c_v$ such that $H_v(\mathbf{a})\leq c_v(1+\|\mathbf{a}\|)$ for all $\mathbf{a}\in G(F)$, and we set $c_\infty=\prod_{v\in\Omega_F^{\infty}}c_v$.
\end{nota}

\begin{rem}\label{Hvbad}
By an argument similar to that in \cite[COROLLARY~3.3.7]{CLT12} and \cite[PROPOSITION~8.3]{PSTVA21}, we have the following statement.
\begin{enumerate}
\item Let $v\in\Omega_F^{<\infty}\setminus S$ and let $\mathbf{a}\in\Lambda_X\setminus\{\mathbf{0}\}$.
Then the function
\[
\mathbf{s}\mapsto\widehat{\mathsf{H}}_{\mathbf{m},v}(\mathbf{a},\mathbf{s})
\]
is holomorphic on the domain
\[
U_{\mathbf{a}}:=\{\mathbf{s}=(s_\alpha)_{\alpha\in\mathcal{A}}\in\Pic(X)_{\mathbb{C}}\mid\Re s_\alpha>\rho_\alpha-1\text{ for all }\alpha\in\mathcal{A}^0(\mathbf{a})\}.
\]
\item There exist positive numbers $\varepsilon$ and $C_1$ such that
\[
|\widehat{\mathsf{H}}_{\mathbf{m},v}(\mathbf{a},\mathbf{s})|\leq C_1(1+\|\mathbf{a}\|)^\varepsilon
\]
holds for all $\mathbf{a}\in\Lambda_X\setminus\{\mathbf{0}\}$, all $v\in\Omega_F^{<\infty}\setminus S$, and all $\mathbf{s}\in U_{\mathbf{a}}$.
\end{enumerate}
\end{rem}

We note that the domain $U_{\mathbf{a}}$ contains $\mathsf{T}_{>-\delta}$ for all $\mathbf{a}\in\Lambda_X\setminus\{\mathbf{0}\}$.
Hence, for $v\in\Omega_F^{<\infty}\setminus S$, the function $\mathbf{s}\mapsto\widehat{\mathsf{H}}_{\mathbf{m},v}(\mathbf{a},\mathbf{s})$ is holomorphic on $\mathsf{T}_{>-\delta}$.

If $v\in S$, then the local height integrals at $v$ coincide with the corresponding local height integrals in \cite{CLT02} and \cite{PSTVA21}.
We recall the following statement.

\begin{rem}\label{HvS}
\leavevmode
\begin{enumerate}
\item Let $v\in\Omega_F$ and let $\mathbf{a}\in\Lambda_X\setminus\{\mathbf{0}\}$.
By \cite[PROPOSITION~8.4]{PSTVA21}, the function
\[
\mathbf{s}\mapsto\widehat{\mathsf{H}}_{\mathbf{m},v}(\mathbf{a},\mathbf{s})
\]
is holomorphic on $\mathsf{T}_{>-\delta}$.
\item Let $l$ be a positive integer.
By \cite[Proposition~10.1]{CLT02}, there exist positive constants $C_l'$ and $M_l'$ such that
\[
\left|\prod_{v\in S}\widehat{\mathsf{H}}_{\mathbf{m},v}(\mathbf{a},\mathbf{s})\right|\leq\frac{C_l'(1+|\mathbf{s}|)^{M_l'}}{(1+\|\mathbf{a}\|)^l}
\]
holds for all $\mathbf{a}\in\Lambda_X\setminus\{\mathbf{0}\}$ and all $\mathbf{s}\in\mathsf{T}_{>-\delta}$.
\end{enumerate}
\end{rem}

\begin{rem}
The statement of \cite[PROPOSITION~8.4]{PSTVA21} is true, but it is not sufficient to prove \cite[PROPOSITION~9.1]{PSTVA21} because, in general, the statement that
\[
\sum_{\mathbf{a}\in\Lambda_X}\frac{1}{(1+H_\infty(\mathbf{a}))^N}
\]
converges for sufficiently large $N$ is false, where
\[
H_\infty(\mathbf{a})=\prod_{v\in\Omega_F^\infty}\max\{|a_1|_v,\dots,|a_n|_v\}.
\]
Indeed, if $n=1$, by the product formula, we have
\[
H_\infty(a)=\prod_{v\in\Omega_F^{<\infty}}|a|_v^{-1}=H_\infty(ua)
\]
for all $a\in\Lambda_X\setminus\{\mathbf{0}\}$ and all $u\in\mathcal{O}_F^\times$.
Hence, the series
\[
\sum_{a\in\Lambda_X}\frac{1}{(1+H_\infty(a))^N}
\]
diverges if $n=1$ and $\#\mathcal{O}_F^\times=\infty$.
This occurs, for example, when $F=\mathbb{Q}(\sqrt{5})$.
Therefore, we have to consider $\|\mathbf{a}\|$ instead of $H_\infty(\mathbf{a})$.
\end{rem}

Finally, we study the height integral
\[
\widehat{\mathsf{H}}_{\mathbf{m}}(\mathbf{a},\mathbf{s})=\prod_{v\in\Omega_F}\widehat{\mathsf{H}}_{\mathbf{m},v}(\mathbf{a},\mathbf{s})
\]
using ideas similar to those in \cite[\S~8.1]{PSTVA21}.
We introduce the following notation.

\begin{nota}
For each $\mathbf{a}\in\Lambda_X\setminus\{\mathbf{0}\}$ and each $i=1,2,3$, we define a subset $S_i^{\mathbf{a}}$ of $\Omega_F$ as follows:
\begin{enumerate}
\item $S_1^{\mathbf{a}}:=\{v\in\Omega_F^{<\infty}\setminus S\mid q_v\leq N\}$,
\item $S_2^{\mathbf{a}}:=\{v\in\Omega_F^{<\infty}\setminus S\mid q_v>N,j_v(\mathbf{a})\neq0\}$,
\item $S_3^{\mathbf{a}}:=\{v\in\Omega_F^{<\infty}\setminus S\mid q_v>N,j_v(\mathbf{a})=0\}$.
\end{enumerate}
Since $S_1^{\mathbf{a}}$ does not depend on $\mathbf{a}$, we sometimes write $S_1$ for $S_1^{\mathbf{a}}$.
\end{nota}

We note that $\Omega_F=S\sqcup S_1^{\mathbf{a}}\sqcup S_2^{\mathbf{a}}\sqcup S_3^{\mathbf{a}}$ for all $\mathbf{a}\in\Lambda_X\setminus\{\mathbf{0}\}$.
By Lemma~\ref{N}~(4), we have $j_v(\mathbf{a})>0$ if $v\in S_2^{\mathbf{a}}$.
The following proposition is an analogue of \cite[PROPOSITION~8.5]{PSTVA21}.

\begin{prop}\label{Hhat}
\leavevmode
\begin{enumerate}
\item Let $\mathbf{a}\in\Lambda_X\setminus\{\mathbf{0}\}$.
Then the function
\[
\mathbf{s}\mapsto\widehat{\mathsf{H}}_{\mathbf{m}}(\mathbf{a},\mathbf{s})\prod_{\alpha\in\mathcal{A}^0(\mathbf{a})}\zeta_{F_\alpha,S^c}(m_\alpha(s_\alpha-\rho_\alpha+1))^{-1}
\]
is holomorphic on $\mathsf{T}_{>-\delta}$, where the function $\zeta_{F_\alpha,S^c}$ is as in Notation~\ref{zeta}~(3).
\item Let $r$ be a positive integer.
Then there exist positive numbers $c_r$ and $M_r$ such that
\[
\left|\widehat{\mathsf{H}}_{\mathbf{m}}(\mathbf{a},\mathbf{s})\prod_{\alpha\in\mathcal{A}^0(\mathbf{a})}\zeta_{F_\alpha,S^c}(m_\alpha(s_\alpha-\rho_\alpha+1))^{-1}\right|\leq\frac{c_r(1+|\mathbf{s}|)^{M_r}}{(1+\|\mathbf{a}\|)^{r}}
\]
for all $\mathbf{a}\in\Lambda_X\setminus\{\mathbf{0}\}$ and all $\mathbf{s}\in\mathsf{T}_{>-\delta}$.
\end{enumerate}
\end{prop}

\begin{proof}
\begin{enumerate}
\item For each $i=1,2,3$, we define a function $f_i^{\mathbf{a}}$ by
\[
f_i^{\mathbf{a}}(\mathbf{s})=\prod_{v\in S_i^{\mathbf{a}}}\left(\widehat{\mathsf{H}}_{\mathbf{m},v}(\mathbf{a},\mathbf{s})\prod_{\alpha\in\mathcal{A}^0(\mathbf{a})}\prod_{\beta\in\mathcal{A}_v(\alpha)}\zeta_{F_{v,\beta}}(m_\alpha(s_\alpha-\rho_\alpha+1))^{-1}\right).
\]
Then we have
\begin{align*}
&\widehat{\mathsf{H}}_{\mathbf{m}}(\mathbf{a},\mathbf{s})\prod_{\alpha\in\mathcal{A}^0(\mathbf{a})}\zeta_{F_\alpha,S^c}(m_\alpha(s_\alpha-\rho_\alpha+1))^{-1}
\\
&=\left(\prod_{\alpha\in\mathcal{A}^0(\mathbf{a})}\prod_{v\in\Omega_F^{<\infty}\setminus S}\prod_{\beta\in\mathcal{A}_v(\alpha)}\zeta_{F_{v,\beta}}(m_\alpha(s_\alpha-\rho_\alpha+1))^{-1}\right)\left(\prod_{v\in\Omega_F}\widehat{\mathsf{H}}_{\mathbf{m},v}(\mathbf{a},\mathbf{s})\right)
\\
&=\left(\prod_{i=1}^3\prod_{v\in S_i^{\mathbf{a}}}\widehat{\mathsf{H}}_{\mathbf{m},v}(\mathbf{a},\mathbf{s})\prod_{\alpha\in\mathcal{A}^0(\mathbf{a})}\prod_{\beta\in\mathcal{A}_v(\alpha)}\zeta_{F_{v,\beta}}(m_\alpha(s_\alpha-\rho_\alpha+1))^{-1}\right)\left(\prod_{v\in S}\widehat{\mathsf{H}}_{\mathbf{m},v}(\mathbf{a},\mathbf{s})\right)
\\
&=\left(\prod_{i=1}^3f_i^{\mathbf{a}}(\mathbf{s})\right)\left(\prod_{v\in S}\widehat{\mathsf{H}}_{\mathbf{m},v}(\mathbf{a},\mathbf{s})\right).
\end{align*}
Since $\#S_1^{\mathbf{a}}<\infty$ and $\#S_2^{\mathbf{a}}<\infty$, the functions $f_1^{\mathbf{a}}$ and $f_2^{\mathbf{a}}$ are holomorphic on $\mathsf{T}_{>-\delta}$.
By Proposition~\ref{Hva}, the infinite product
\[
f_3^{\mathbf{a}}(\mathbf{s})=\prod_{v\in S_3^{\mathbf{a}}}\left(\widehat{\mathsf{H}}_{\mathbf{m},v}(\mathbf{a},\mathbf{s})\prod_{\alpha\in\mathcal{A}^0(\mathbf{a})}\prod_{\beta\in\mathcal{A}_v(\alpha)}\zeta_{F_{v,\beta}}(m_\alpha(s_\alpha-\rho_\alpha+1))^{-1}\right)
\]
converges absolutely and uniformly on compact subsets of $\mathsf{T}_{>-\delta}$.
Moreover, by Remark~\ref{Hvbad}~(1), each factor of the above infinite product is holomorphic on $\mathsf{T}_{>-\delta}$.
Hence, the above infinite product defines a holomorphic function on $\mathsf{T}_{>-\delta}$.
By Remark~\ref{HvS}~(1), the function $\prod_{v\in S}\widehat{\mathsf{H}}_{\mathbf{m},v}(\mathbf{a},\mathbf{s})$ is holomorphic on $\mathsf{T}_{>-\delta}$.
Therefore, the function
\[
\mathbf{s}\mapsto\widehat{\mathsf{H}}_{\mathbf{m}}(\mathbf{a},\mathbf{s})\prod_{\alpha\in\mathcal{A}^0(\mathbf{a})}\zeta_{F_\alpha,S^c}(m_\alpha(s_\alpha-\rho_\alpha+1))^{-1}
\]
is holomorphic on $\mathsf{T}_{>-\delta}$.
\item We evaluate functions $f_1^{\mathbf{a}}(\mathbf{s})$, $f_2^{\mathbf{a}}(\mathbf{s})$, $f_3^{\mathbf{a}}(\mathbf{s})$, and $\prod_{v\in S}\widehat{\mathsf{H}}_{\mathbf{m},v}(\mathbf{a},\mathbf{s})$.
We fix positive constants $\varepsilon$ and $C_1$ as in Remark~\ref{Hvbad}.
By Remark~\ref{Hvbad}~(2), for all $\mathbf{a}\in\Lambda_X\setminus\{\mathbf{0}\}$ and all $\mathbf{s}\in\mathsf{T}_{>-\delta}$, we have
\[
\left|\prod_{v\in S_1^{\mathbf{a}}}\widehat{\mathsf{H}}_{\mathbf{m},v}(\mathbf{a},\mathbf{s})\right|\leq C_1^{\#S_1}(1+\|\mathbf{a}\|)^{\varepsilon\#S_1}.
\]
By a calculation similar to that leading to \eqref{7} in the proof of Proposition~\ref{Hvajneq0}, we have
\[
\left|\prod_{v\in S_1^{\mathbf{a}}}\prod_{\alpha\in\mathcal{A}^0(\mathbf{a})}\prod_{\beta\in\mathcal{A}_v(\alpha)}\zeta_{F_{v,\beta}}(m_\alpha(s_\alpha-\rho_\alpha+1))^{-1}\right|\leq 2^{(\#S_1)(\#\overline{\mathcal{A}})}.
\]
It follows that
\begin{equation}\label{9}
|f_1^{\mathbf{a}}(\mathbf{s})|\leq(C_1 2^{\#\overline{\mathcal{A}}})^{\#S_1}(1+\|\mathbf{a}\|)^{\varepsilon\#S_1}.
\end{equation}
By Proposition~\ref{Hvajneq0}, we have a constant $C_2>1$ such that
\[
\left|\widehat{\mathsf{H}}_{\mathbf{m},v}(\mathbf{a},\mathbf{s})\prod_{\alpha\in\mathcal{A}^0(\mathbf{a})}\prod_{\beta\in\mathcal{A}_v(\alpha)}\zeta_{F_{v,\beta}}(m_\alpha(s_\alpha-\rho_\alpha+1))^{-1}\right|<C_2
\]
for all $\mathbf{s}\in\mathsf{T}_{>-\delta}$, $\mathbf{a}\in\Lambda_X\setminus\{\mathbf{0}\}$, and $v\in S_2^{\mathbf{a}}$.
It follows that
\[
|f_2^{\mathbf{a}}(\mathbf{s})|\leq C_2^{\#S_2^{\mathbf{a}}}.
\]
Since $\Lambda_X$ is a finitely generated $\mathcal{O}_F$-module, we have
\[
R=\min\left\{\inf_{\mathbf{a}\in\Lambda_X\setminus\{\mathbf{0}\}}\inf_{v\in\Omega_F^{<\infty}}q_v^{j_v(\mathbf{a})},1\right\}>0.
\]
We have $q_v^{j_v(\mathbf{a})}\geq R$ if $v\in S_1\cup(S\cap\Omega_F^{<\infty})$, $q_v^{j_v(\mathbf{a})}\geq2$ if $v\in S_2^{\mathbf{a}}$, and $q_v^{j_v(\mathbf{a})}=1$ if $v\in S_3^{\mathbf{a}}$.
Hence, we obtain
\begin{align*}
2^{\#S_2^{\mathbf{a}}}R^{\#(S\cup S_1)}&\leq\prod_{v\in\Omega_F^{<\infty}}q_v^{j_v(\mathbf{a})}
\\
&=\prod_{v\in\Omega_F^{<\infty}}q_v^{\min\{v(a_1),\dots,v(a_n)\}}
\\
&=\prod_{v\in\Omega_F^{<\infty}}\max\{|a_1|_v,\dots,|a_n|_v\}^{-1}
\\
&\leq\prod_{v\in\Omega_F^\infty}\max\{|a_1|_v,\dots,|a_n|_v\}
\\
&=\prod_{v\in\Omega_F^\infty}H_v(\mathbf{a})
\\
&\leq c_\infty(1+\|\mathbf{a}\|)^{[F:\mathbb{Q}]},
\end{align*}
where the second inequality follows from the product formula.
Therefore, we have
\[
\#S_2^{\mathbf{a}}\leq[F:\mathbb{Q}]\log_2(1+\|\mathbf{a}\|)+\log_2c_\infty-\#(S\cup S_1)\log_2R
\]
and
\begin{align}
|f_2^{\mathbf{a}}(\mathbf{s})|&\leq C_2^{[F:\mathbb{Q}]\log_2(1+\|\mathbf{a}\|)+\log_2c_\infty-\#(S\cup S_1)\log_2R}\notag
\\
&=(1+\|\mathbf{a}\|)^{[F:\mathbb{Q}]\log_2C_2}C_2^{\log_2c_\infty-\#(S\cup S_1)\log_2R}.\label{10}
\end{align}
By Proposition~\ref{Hva}, the infinite product defining $f_3^{\mathbf{a}}(\mathbf{s})$ converges uniformly for $\mathbf{s}\in\mathsf{T}_{>-\delta}$ and $\mathbf{a}\in\Lambda_X\setminus\{\mathbf{0}\}$.
Hence, we have
\begin{equation}\label{11}
C_3:=\sup_{\mathbf{a}\in\Lambda_X\setminus\{\mathbf{0}\}}\sup_{\mathbf{s}\in\mathsf{T}_{>-\delta}}|f_3^{\mathbf{a}}(\mathbf{s})|<\infty.
\end{equation}
Let $l_r$ be a positive integer such that
\[
l_r>r+[F:\mathbb{Q}]\log_2C_2+\varepsilon\#S_1.
\]
By Remark~\ref{HvS}~(2), we have
\begin{equation}\label{12}
\left|\prod_{v\in S}\widehat{\mathsf{H}}_{\mathbf{m},v}(\mathbf{a},\mathbf{s})\right|\leq\frac{C_{l_r}'(1+|\mathbf{s}|)^{M_{l_r}'}}{(1+\|\mathbf{a}\|)^{l_r}}.
\end{equation}
Combining \eqref{9}, \eqref{10}, \eqref{11}, and \eqref{12}, we have 
\begin{align*}
&\left|\widehat{\mathsf{H}}_{\mathbf{m}}(\mathbf{a},\mathbf{s})\prod_{\alpha\in\mathcal{A}^0(\mathbf{a})}\zeta_{F_\alpha,S^c}(m_\alpha(s_\alpha-\rho_\alpha+1))^{-1}\right|
\\
&=\left(\prod_{i=1}^3|f_i^{\mathbf{a}}(\mathbf{s})|\right)\left|\prod_{v\in S}\widehat{\mathsf{H}}_{\mathbf{m},v}(\mathbf{a},\mathbf{s})\right|
\\
&\leq(C_1 2^{\#\overline{\mathcal{A}}})^{\# S_1}C_2^{\log_2c_\infty-\#(S\cup S_1)\log_2R}C_3C_{l_r}'\frac{(1+|\mathbf{s}|)^{M_{l_r}'}}{(1+\|\mathbf{a}\|)^{l_r-[F:\mathbb{Q}]\log_2C_2-\varepsilon\#S_1}}
\\
&<(C_1 2^{\#\overline{\mathcal{A}}})^{\# S_1}C_2^{\log_2c_\infty-\#(S\cup S_1)\log_2R}C_3C_{l_r}'\frac{(1+|\mathbf{s}|)^{M_{l_r}'}}{(1+\|\mathbf{a}\|)^r}
\end{align*}
for all $\mathbf{a}\in\Lambda_X\setminus\{\mathbf{0}\}$ and $\mathbf{s}\in\mathsf{T}_{>-\delta}$.
Therefore, the proposition holds for
\begin{align*}
c_r&=(C_1 2^{\#\overline{\mathcal{A}}})^{\# S_1}C_2^{\log_2c_\infty-\#(S\cup S_1)\log_2R}C_3C_{l_r}',
\\
M_r&=M_{l_r}'.\qedhere
\end{align*}
\end{enumerate}
\end{proof}

We show the third condition of the Poisson summation formula.

\begin{prop}\label{conv}
For all $\mathbf{s}\in\mathsf{T}_{>0}$, the series
\[
\sum_{\mathbf{a}\in\Lambda_X}\widehat{\mathsf{H}}_{\mathbf{m}}(\mathbf{a},\mathbf{s})
\]
converges absolutely.
Moreover, if $\Re\mathbf{s}$ is sufficiently large, then the series
\[
\sum_{\mathbf{a}\in G(F)}\widehat{\mathsf{H}}_{\mathbf{m}}(\mathbf{a},\mathbf{s})
\]
converges absolutely.
\end{prop}

\begin{proof}
By \cite[\S~3.5.2]{CLT12}, there exists a positive integer $r$ such that the series
\[
\sum_{\mathbf{a}\in\Lambda_X}\frac{1}{(1+\|\mathbf{a}\|)^r}
\]
converges.
We fix such a positive integer $r$.
Let $\mathbf{s}=(s_\alpha)_{\alpha\in\mathcal{A}}\in\mathsf{T}_{>0}$, and let $c_r$ and $M_r$ be positive constants as in Proposition~\ref{Hhat}.
We define a positive constant $c_{\mathbf{s}}$ by
\[
c_{\mathbf{s}}:=\max_{B\subseteq\mathcal{A}}\prod_{\alpha\in B}|\zeta_{F_\alpha,S^c}(m_\alpha(s_\alpha-\rho_\alpha+1))|<\infty.
\]
Then we have
\[
|\widehat{\mathsf{H}}_{\mathbf{m}}(\mathbf{a},\mathbf{s})|\leq\frac{c_rc_{\mathbf{s}}(1+\|\mathbf{s}\|)^{M_r}}{(1+\|\mathbf{a}\|)^{r}}
\]
for all $\mathbf{a}\in\Lambda_X\setminus\{\mathbf{0}\}$.
Therefore, we obtain
\[
\sum_{\mathbf{a}\in\Lambda_X\setminus\{\mathbf{0}\}}|\widehat{\mathsf{H}}_{\mathbf{m}}(\mathbf{a},\mathbf{s})|\leq c_rc_{\mathbf{s}}(1+|\mathbf{s}|)^{M_r}\sum_{\mathbf{a}\in\Lambda_X}\frac{1}{(1+\|\mathbf{a}\|)^r}<\infty.
\]
The term corresponding to $\mathbf{a}=\mathbf{0}$ is finite for every
$\mathbf{s}\in\mathsf{T}_{>0}$.
Hence the first assertion follows.
The last assertion follows from Proposition~\ref{Hzero}.
\end{proof}

In the following, we analyze the poles of the height zeta function.
We recall that
\[
\mathcal{A}^0(\mathbf{a})=\{\alpha\in\mathcal{A}\mid d_\alpha(f_{\mathbf{a}})=0\}
\]
for $\mathbf{a}\in G(F)\setminus\{\mathbf{0}\}$.
We introduce the following notation.

\begin{nota}
\leavevmode
\begin{enumerate}
\item We define a set $\mathscr{A}$ by
\[
\mathscr{A}:=\{\mathcal{A}^0(\mathbf{a})\mid\mathbf{a}\in\Lambda_X\setminus\{\mathbf{0}\}\}.
\]
\item For each $\mathcal{A}'\in\mathscr{A}$, we define a subset $\Lambda_{\mathcal{A}'}$ of $\Lambda_X$ by
\[
\Lambda_{\mathcal{A}'}:=\{\mathbf{a}\in\Lambda_X\setminus\{\mathbf{0}\}\mid\mathcal{A}'=\mathcal{A}^0(\mathbf{a})\}.
\]
\end{enumerate}
\end{nota}
We note that $\Lambda_X\setminus\{\mathbf{0}\}=\coprod_{\mathcal{A}'\in\mathscr{A}}\Lambda_{\mathcal{A}'}$.
By the Poisson summation formula, we have
\[
\mathsf{Z}_{\mathbf{m}}(\mathbf{s})=\sum_{\mathbf{a}\in\Lambda_X}\widehat{\mathsf{H}}_{\mathbf{m}}(\mathbf{a},\mathbf{s})=\widehat{\mathsf{H}}_{\mathbf{m}}(\mathbf{0},\mathbf{s})+\sum_{\mathcal{A}'\in\mathscr{A}}\sum_{\mathbf{a}\in\Lambda_{\mathcal{A}'}}\widehat{\mathsf{H}}_{\mathbf{m}}(\mathbf{a},\mathbf{s}).
\]

\begin{prop}\label{phihol}
Let $\mathcal{A}'\in\mathscr{A}$.
Then the function
\[
\phi_{\mathcal{A}'}(\mathbf{s}):=\sum_{\mathbf{a}\in\Lambda_{\mathcal{A}'}}\widehat{\mathsf{H}}_{\mathbf{m}}(\mathbf{a},\mathbf{s})\prod_{\alpha\in\mathcal{A}'}\zeta_{F_\alpha}(m_\alpha(s_\alpha-\rho_\alpha+1))^{-1}
\]
is holomorphic on $\mathsf{T}_{>-\delta}$.
\end{prop}

\begin{proof}
By Proposition~\ref{Hhat}~(1), the function
\[
\phi_{\mathbf{a}}(\mathbf{s}):=\widehat{\mathsf{H}}_{\mathbf{m}}(\mathbf{a},\mathbf{s})\prod_{\alpha\in\mathcal{A}'}\zeta_{F_\alpha,S^c}(m_\alpha(s_\alpha-\rho_\alpha+1))^{-1}
\]
is holomorphic on $\mathsf{T}_{>-\delta}$ for all $\mathbf{a}\in\Lambda_{\mathcal{A}'}$.
By Proposition~\ref{Hhat}~(2), the series
\[
\sum_{\mathbf{a}\in\Lambda_{\mathcal{A}'}}\phi_{\mathbf{a}}(\mathbf{s})
\]
converges absolutely and uniformly on compact subsets of $\mathsf{T}_{>-\delta}$, and defines a holomorphic function on $\mathsf{T}_{>-\delta}$.
We recall that
\[
\zeta_{F_\alpha}(s)=\zeta_{F_\alpha,S^c}(s)\prod_{v\in S\setminus\Omega_F^{\infty}}\prod_{\beta\in\mathcal{A}_v(\alpha)}\zeta_{F_{v,\beta}}(s)
\]
for each $\alpha\in\mathcal{A}$ as in Remark~\ref{z}.
Hence, we obtain
\[
\phi_{\mathcal{A}'}(\mathbf{s})=\sum_{\mathbf{a}\in\Lambda_{\mathcal{A}'}}\phi_{\mathbf{a}}(\mathbf{s})\prod_{\alpha\in\mathcal{A}'}\prod_{v\in S\setminus\Omega_F^\infty}\prod_{\beta\in\mathcal{A}_v(\alpha)}\zeta_{F_{v,\beta}}(m_\alpha(s_\alpha-\rho_\alpha+1))^{-1},
\]
and it follows that $\phi_{\mathcal{A}'}$ is holomorphic on $\mathsf{T}_{>-\delta}$.
\end{proof}

We set
\begin{align*}
f(\mathbf{s})&:=\sum_{\mathbf{a}\in\Lambda_X}\widehat{\mathsf{H}}_{\mathbf{m}}(\mathbf{a},\mathbf{s})\prod_{\alpha\in\mathcal{A}}\zeta_{F_\alpha}(m_\alpha(s_\alpha-\rho_\alpha+1))^{-1}
\\
&=\widehat{\mathsf{H}}_{\mathbf{m}}(\mathbf{0},\mathbf{s})\prod_{\alpha\in\mathcal{A}}\zeta_{F_\alpha}(m_\alpha(s_\alpha-\rho_\alpha+1))^{-1}+\sum_{\mathcal{A}'\in\mathscr{A}}\phi_{\mathcal{A}'}(\mathbf{s})\prod_{\alpha\in\mathcal{A}\setminus\mathcal{A}'}\zeta_{F_\alpha}(m_\alpha(s_\alpha-\rho_\alpha+1))^{-1}.
\end{align*}
For $\mathbf{s}\in\Pic(X)_{\mathbb{C}}$ with $\Re\mathbf{s}$ sufficiently large, we have
\begin{equation}\label{13}
\mathsf{Z}_{\mathbf{m}}(\mathbf{s})=f(\mathbf{s})\prod_{\alpha\in\mathcal{A}}\zeta_{F_\alpha}(m_\alpha(s_\alpha-\rho_\alpha+1)).
\end{equation}
By Corollary~\ref{fzero} and Proposition~\ref{phihol}, the function $f$ is meromorphic on $\mathsf{T}_{>-\delta}$.
In particular, \eqref{13} gives a meromorphic continuation of $\mathsf{Z}_{\mathbf{m}}$ to $\mathsf{T}_{>-\delta}$.

Let $L=\sum_{\alpha\in\mathcal{A}}\lambda_\alpha D_\alpha$ be a big line bundle on $X$.
We study the rightmost pole of $\mathsf{Z}_{\mathbf{m}}(sL)$ to apply the Tauberian theorem.
We recall that
\begin{align*}
a&=\max_{\alpha\in\mathcal{A}}\frac{\rho_\alpha-1+m_\alpha^{-1}}{\lambda_\alpha},
\\
\mathcal{A}_{\mathbf{m}}(L)&=\left\{\alpha\in\mathcal{A}\middle|\frac{\rho_\alpha-1+m_\alpha^{-1}}{\lambda_\alpha}=a\right\},
\\
b&=\#\mathcal{A}_{\mathbf{m}}(L),
\end{align*}
and fix a positive constant $\delta'$ such that
\[
0<\delta'<\min\left\{\left\{a-\frac{\rho_\alpha-1+m_\alpha^{-1}}{\lambda_\alpha}\,\middle|\,\alpha\in\mathcal{A}\setminus\mathcal{A}_{\mathbf{m}}(L)\right\}\cup\left\{\frac{\delta}{\lambda_\alpha}\,\middle|\,\alpha\in\mathcal{A}\right\}\right\}.
\]
We proved that $\widehat{\mathsf{H}}_{\mathbf{m}}(\mathbf{0},sL)$ has a rightmost pole at $s=a$ of order $b$, and the constant $c=\lim_{s\to a}(s-a)^b\widehat{\mathsf{H}}_{\mathbf{m}}(\mathbf{0},sL)$ is positive in Proposition~\ref{c}.
We show that the height integral $\widehat{\mathsf{H}}_{\mathbf{m}}(\mathbf{a},\mathbf{s})$ does not contribute to the rightmost pole of $\mathsf{Z}_{\mathbf{m}}$ for all $\mathbf{a}\in\Lambda_X\setminus\{\mathbf{0}\}$.
We recall that a $\mathbb{Q}$-Cartier divisor on $X$ is rigid if its Iitaka dimension (see \cite[Definition~2.1.3]{Laz04} for the definition) is zero.
The following statement is analogous to the corresponding statement in \cite[\S~9.1]{PSTVA21}.

\begin{prop}\label{Ha}
We assume that $aL+K_X+D_{\mathbf{m}}$ is rigid.
Then, for all $\mathbf{a}\in\Lambda_X\setminus\{\mathbf{0}\}$, the function
\[
s\mapsto(s-a)^{b-1}\widehat{\mathsf{H}}_{\mathbf{m}}(\mathbf{a},sL)
\]
is holomorphic on the half-plane $\Re s>a-\delta'$.
\end{prop}

\begin{proof}
If the function $s\mapsto\widehat{\mathsf{H}}_{\mathbf{m}}(\mathbf{a},sL)$ is holomorphic on the half-plane $\Re s>a-\delta'$, then the assertion is immediate.
We may assume that the above function is not holomorphic on the half-plane $\Re s>a-\delta'$.
Since $\delta'<\delta/\lambda_\alpha$ for all $\alpha\in\mathcal{A}$, we have $sL\in\mathsf{T}_{>-\delta}$ for $\Re s>a-\delta'$.
We define a holomorphic function $\phi_{\mathbf{a}}$ on $\mathsf{T}_{>-\delta}$ as in the proof of Proposition~\ref{phihol}.
Then we have
\[
\widehat{\mathsf{H}}_{\mathbf{m}}(\mathbf{a},sL)=\phi_{\mathbf{a}}(sL)\prod_{\alpha\in\mathcal{A}^0(\mathbf{a})}\zeta_{F_\alpha,S^c}(m_\alpha(s\lambda_\alpha-\rho_\alpha+1)).
\]
We note that the function $\phi_{\mathbf{a}}(sL)$ is holomorphic on the half-plane $\Re s>a-\delta'$.
If $\mathcal{A}^0(\mathbf{a})\cap\mathcal{A}_{\mathbf{m}}(L)=\varnothing$, then, for all $\alpha\in\mathcal{A}^0(\mathbf{a})$, we have
\[
\Re(m_\alpha(s\lambda_\alpha-\rho_\alpha+1))>1
\]
and the function $\widehat{\mathsf{H}}_{\mathbf{m}}(\mathbf{a},\mathbf{s})$ is holomorphic on $\mathsf{T}_{>-\delta'}$.
This is a contradiction.
Hence, the function
\[
s\mapsto\prod_{\alpha\in\mathcal{A}^0(\mathbf{a})}\zeta_{F_\alpha,S^c}(m_\alpha(s\lambda_\alpha-\rho_\alpha+1))
\]
has a rightmost pole at $s=a$ of order $b':=\#(\mathcal{A}^0(\mathbf{a})\cap\mathcal{A}_{\mathbf{m}}(L))$, and is holomorphic at all other points in the half-plane $\Re s>a-\delta'$.
Hence, the function $s\mapsto\widehat{\mathsf{H}}_{\mathbf{m}}(\mathbf{a},sL)$ has a rightmost pole at $s=a$ and is holomorphic at all other points in the half-plane $\Re s>a-\delta'$.
We denote by $b''$ the order of the pole of $\widehat{\mathsf{H}}_{\mathbf{m}}(\mathbf{a},sL)$ at $s=a$.
Since $aL+K_X+D_{\mathbf{m}}$ is rigid, we have $\mathcal{A}_{\mathbf{m}}(L)\nsubseteq\mathcal{A}^0(\mathbf{a})$ as in \cite[\S~9.1]{PSTVA21}.
It follows that 
\[
b''\leq b'=\#(\mathcal{A}^0(\mathbf{a})\cap\mathcal{A}_{\mathbf{m}}(L))<\#\mathcal{A}_{\mathbf{m}}(L)=b.
\]
Therefore, the function
\[
s\mapsto(s-a)^{b-1}\widehat{\mathsf{H}}_{\mathbf{m}}(\mathbf{a},sL)
\]
is holomorphic on the half-plane $\Re s>a-\delta'$.
\end{proof}

Finally, we determine the rightmost pole of $\mathsf{Z}_{\mathbf{m}}(sL)$.

\begin{prop}\label{Haneq0}
Suppose that $aL+K_X+D_{\mathbf{m}}$ is rigid.
Then the function
\[
s\mapsto(s-a)^{b-1}\sum_{\mathbf{a}\in\Lambda_X\setminus\{\mathbf{0}\}}\widehat{\mathsf{H}}_{\mathbf{m}}(\mathbf{a},sL)
\]
is holomorphic on the half-plane $\Re s>a-\delta'$.
Moreover, the function $\mathsf{Z}_{\mathbf{m}}(sL)$ has a rightmost pole at $s=a$ of order $b$, and we have
\[
\lim_{s\to a}(s-a)^b\mathsf{Z}_{\mathbf{m}}(sL)=\lim_{s\to a}(s-a)^b\widehat{\mathsf{H}}_{\mathbf{m}}(\mathbf{0},sL)=c>0.
\]
\end{prop}

\begin{proof}
Since $aL+K_X+D_{\mathbf{m}}$ is rigid, for all $\mathcal{A}'\in\mathscr{A}$, we have $\#(\mathcal{A}'\cap\mathcal{A}_{\mathbf{m}}(L))\leq b-1$ as in \cite[\S~9.1]{PSTVA21}.
Hence, we have
\begin{align*}
&(s-a)^{b-1}\sum_{\mathbf{a}\in\Lambda_X\setminus\{\mathbf{0}\}}\widehat{\mathsf{H}}_{\mathbf{m}}(\mathbf{a},sL)
\\
&=(s-a)^{b-1}\sum_{\mathcal{A}'\in\mathscr{A}}\sum_{\mathbf{a}\in\Lambda_{\mathcal{A}'}}\widehat{\mathsf{H}}_{\mathbf{m}}(\mathbf{a},sL)
\\
&=(s-a)^{b-1}\sum_{\mathcal{A}'\in\mathscr{A}}\sum_{\mathbf{a}\in\Lambda_{\mathcal{A}'}}\phi_{\mathbf{a}}(sL)\prod_{\alpha\in\mathcal{A}'}\zeta_{F_\alpha,S^c}(m_\alpha(s\lambda_\alpha-\rho_\alpha+1))
\\
&=\sum_{\mathcal{A}'\in\mathscr{A}}(s-a)^{b-1}\left(\sum_{\mathbf{a}\in\Lambda_{\mathcal{A}'}}\phi_{\mathbf{a}}(sL)\right)\prod_{\alpha\in\mathcal{A}'}\zeta_{F_\alpha,S^c}(m_\alpha(s\lambda_\alpha-\rho_\alpha+1))
\\
&=\sum_{\mathcal{A}'\in\mathscr{A}}(s-a)^{b-1-\#(\mathcal{A}'\cap\mathcal{A}_{\mathbf{m}}(L))}\phi_{\mathcal{A}'}(sL)
\\
&\qquad\times\left(\prod_{\alpha\in\mathcal{A}'\setminus\mathcal{A}_{\mathbf{m}}(L)}\zeta_{F_\alpha}(m_\alpha(s\lambda_\alpha-\rho_\alpha+1))\right)\left(\prod_{\alpha\in\mathcal{A}'\cap\mathcal{A}_{\mathbf{m}}(L)}(s-a)\zeta_{F_\alpha}(m_\alpha(s\lambda_\alpha-\rho_\alpha+1))\right).
\end{align*}
By Proposition~\ref{phihol}, the definition of $\mathcal{A}_{\mathbf{m}}(L)$, and the fact that the Dedekind zeta function $\zeta_{F_\alpha}(s)$ has a single pole at $s=1$, the function
\[
s\mapsto(s-a)^{b-1}\sum_{\mathbf{a}\in\Lambda_X\setminus\{\mathbf{0}\}}\widehat{\mathsf{H}}_{\mathbf{m}}(\mathbf{a},sL)
\]
is holomorphic on the half-plane $\Re s>a-\delta'$.

By Proposition~\ref{c}, the function
\[
(s-a)^b\mathsf{Z}_{\mathbf{m}}(sL)=(s-a)^b\widehat{\mathsf{H}}_{\mathbf{m}}(\mathbf{0},sL)+(s-a)^b\sum_{\mathbf{a}\in\Lambda_X\setminus\{\mathbf{0}\}}\widehat{\mathsf{H}}_{\mathbf{m}}(\mathbf{a},sL)
\]
is holomorphic on the half-plane $\Re s>a-\delta'$, and we have
\[
\lim_{s\to a}(s-a)^b\mathsf{Z}_{\mathbf{m}}(sL)=\lim_{s\to a}(s-a)^b\left(\sum_{\mathbf{a}\in\Lambda_X}\widehat{\mathsf{H}}_{\mathbf{m}}(\mathbf{a},sL)\right)=\lim_{s\to a}(s-a)^b\widehat{\mathsf{H}}_{\mathbf{m}}(\mathbf{0},sL)=c>0.
\]
Thus $s\mapsto\mathsf{Z}_{\mathbf{m}}(sL)$ has a rightmost pole at $s=a$
of order $b$.
\end{proof}

By Landau's theorem, the Dirichlet series
\[
\mathsf{Z}_{\mathbf{m}}(sL)=\sum_{\mathbf{x}\in(\mathcal{X},\mathcal{D}_{\mathbf{m}})^{\mathrm{D}}(\mathcal{O}_{F,S})}\frac{1}{H_L(\mathbf{x})^s}
\]
converges for $\Re s>a$.
Applying the Tauberian theorem to the function $s\mapsto\mathsf{Z}_{\mathbf{m}}(sL)$, we obtain the main theorem.

\begin{theorem}
Let $F$ be a number field, let $G=\mathbb{G}_{a,F}^n$ be the $n$-dimensional vector group over $F$, let $X$ be a smooth projective equivariant compactification of $G$, and let $D=X\setminus G=\bigcup_{\alpha\in\mathcal{A}}D_\alpha$ be the decomposition into irreducible components.
We assume that $\sum_{\alpha\in\mathcal{A}}D_\alpha$ is a strict normal crossings divisor.
Let $(X,D_{\mathbf{m}})$ be a Campana orbifold over $F$ associated with $\mathbf{m}=(m_\alpha)_{\alpha\in\mathcal{A}}\in\mathbb{Z}_{>0}^{\mathcal{A}}$, let $S$ be a finite set of places of $F$ containing all infinite places, and let $(\mathcal{X},\mathcal{D}_{\mathbf{m}})$ be a good integral model of $(X,D_{\mathbf{m}})$ away from $S$.
For each $\alpha\in\mathcal{A}$, we fix an adelic metric on $\mathcal{O}(D_\alpha)$.
Let $L=\sum_{\alpha\in\mathcal{A}}\lambda_\alpha D_\alpha$ be a big line bundle on $X$ equipped with an adelic metric as in \cite[\S~1.3]{Pey95}, and let $H_L\colon G(F)\to\mathbb{R}_{>0}$ be the height function associated with $L$.
If $aL+K_X+D_{\mathbf{m}}$ is rigid, then we have the asymptotic formula
\[
\#\{P\in(\mathcal{X},\mathcal{D}_{\mathbf{m}})^{\mathrm{D}}(\mathcal{O}_{F,S})\mid H_L(P)\leq B\}\sim\frac{c}{a(b-1)!}B^a(\log B)^{b-1}\quad(B\to\infty).
\]
Here $(\mathcal{X},\mathcal{D}_{\mathbf{m}})^{\mathrm{D}}(\mathcal{O}_{F,S})$ is the set of Darmon $\mathcal{O}_{F,S}$-points on $(\mathcal{X},\mathcal{D}_{\mathbf{m}})$, and the constants $a$, $b$, and $c$ are defined by
\begin{align*}
a&=\max_{\alpha\in\mathcal{A}}\frac{\rho_\alpha-1+m_\alpha^{-1}}{\lambda_\alpha},
\\
b&=\#\left\{\alpha\in\mathcal{A}\,\middle|\,\frac{\rho_\alpha-1+m_\alpha^{-1}}{\lambda_\alpha}=a\right\},
\\
c&=\lim_{s\to a}(s-a)^b\widehat{\mathsf{H}}_{\mathbf{m}}(\mathbf{0},sL),
\end{align*}
where $\boldsymbol{\rho}=(\rho_\alpha)_{\alpha\in\mathcal{A}}$ corresponds to the anticanonical bundle on $X$, and the function $\widehat{\mathsf{H}}_{\mathbf{m}}(\mathbf{0},sL)$ is defined in Definition~\ref{HI}~(2).
\end{theorem}

We note that the constants $a$ and $b$ coincide with those in \cite{PSTVA21} corresponding to Campana $\mathcal{O}_{F,S}$-points on $(\mathcal{X},\mathcal{D}_{\mathbf{m}})$, but the constant $c$ differs in general.

\section{Compatibility with Manin's conjecture for $\mathcal{M}$-points}

In this section, we show that Theorem~\ref{main} is compatible with Manin's conjecture for $\mathcal{M}$-points formulated in \cite{Moe26a}.
The notion of $\mathcal{M}$-points was introduced in \cite[Definition~3.12]{Moe24} and \cite[Definition~3.12]{Moe26a}; it generalizes semi-integral points, including Darmon points.

We continue to use the notation fixed in the previous section.
In particular, $F$ is a number field and $X$ is a smooth projective equivariant compactification of the $n$-dimensional vector group $G=\mathbb{G}_{a,F}^n$ over $F$.
Moreover,
\[
X\setminus G=\bigcup_{\alpha\in\mathcal{A}}D_\alpha
\]
is the decomposition into irreducible components, and we assume that the boundary divisor
\[
D=\sum_{\alpha\in\mathcal{A}}D_\alpha
\]
is a strict normal crossings divisor.
We fix a tuple $\mathbf{m}=(m_\alpha)_{\alpha\in\mathcal{A}}\in\mathbb{Z}_{>0}^{\mathcal{A}}$ of positive integers and define a $\mathbb{Q}$-divisor $D_{\mathbf{m}}$ on $X$ by
\[
D_{\mathbf{m}}=\sum_{\alpha\in\mathcal{A}}\left(1-\frac{1}{m_\alpha}\right)D_\alpha.
\]
We fix a finite set $S$ of places of $F$ containing all infinite places and a good integral model $(\mathcal{X},\mathcal{D}_{\mathbf{m}})$ of $(X,D_{\mathbf{m}})$ away from $S$.
For each $\alpha\in\mathcal{A}$, $\mathcal{O}(D_\alpha)$ denotes the line bundle associated with $D_\alpha$, and we fix a section $\mathsf{f}_\alpha$ of $\mathcal{O}(D_\alpha)$ and a smooth adelic metric on $\mathcal{O}(D_\alpha)$.
Recall that $\overline{\mathbb{N}}=\mathbb{Z}_{\geq0}\cup\{\infty\}$.
We set
\[
\mathfrak{M}=\prod_{\alpha\in\mathcal{A}}m_\alpha\mathbb{Z}_{\geq0}\subseteq\overline{\mathbb{N}}^{\mathcal{A}},\quad M=((D_\alpha)_{\alpha\in\mathcal{A}},\mathfrak{M}),\quad\mathcal{M}=((\mathcal{D}_\alpha)_{\alpha\in\mathcal{A}},\mathfrak{M}).
\]
By \cite[Definition~3.24]{Moe26a}, $\mathcal{M}$-points on $(\mathcal{X},\mathcal{M})$ over $\mathcal{O}_{F,S}$ coincide with Darmon points on $(\mathcal{X},\mathcal{D}_{\mathbf{m}})$.

\subsection{The $a$-invariant and the $b$-invariant}

In this subsection, we show that Theorem~\ref{main} is compatible with Manin's conjecture for $\mathcal{M}$-points formulated in \cite[Conjecture~1.4]{Moe26a} with respect to the $a$-invariant and the $b$-invariant.
We recall the Picard group and the effective cone of a pair to state the $a$-invariant and the $b$-invariant defined in \cite[Definition~4.37]{Moe26a}.

\begin{dfn}
\leavevmode
\begin{enumerate}
\item Following \cite[Definition~4.10]{Moe26a}, we define the \textbf{group of divisors} on $(X,M)$ by
\[
\Div(X,M):=\Div(G)\oplus\bigoplus_{\alpha\in\mathcal{A}}\mathbb{Z}\tilde{D}_\alpha,
\]
where $\tilde{D}_\alpha$ is a formal symbol.
Let $D'\in\Div(X,M)$. We say that $D'$ is an \textbf{effective divisor} on $(X,M)$ if there exist an effective divisor $D_G'$ on $G$ and a tuple of non-negative integers $(n_\alpha)_{\alpha\in\mathcal{A}}$ such that
\[
D'=D_G'+\sum_{\alpha\in\mathcal{A}}n_\alpha\tilde{D}_\alpha.
\]
\item We use the natural identification
\[
\Div(X)=\Div(G)\oplus\bigoplus_{\alpha\in\mathcal A}\mathbb ZD_\alpha.
\]
Following \cite[Definition~4.15]{Moe26a}, we define the \textbf{pullback}
\[
\pr_M^*\colon \Div(X)\to\Div(X,M)
\]
by
\[
\pr_M^*\left(D_G+\sum_{\alpha\in\mathcal A}n_\alpha D_\alpha\right)=D_G+\sum_{\alpha\in\mathcal A}n_\alpha m_\alpha\tilde D_\alpha,
\]
where $D_G\in\Div(G)$ and $n_\alpha\in\mathbb{Z}$.
\end{enumerate}
\end{dfn}

\begin{rem}
By \cite[Example~4.28]{Moe26a}, $\Div(X,M)$ is computed as above in our setting.
For every $D_G\in\Div(G)$, no $D_\alpha$ is contained in the support of the Zariski closure of $D_G$ in $X$.
By \cite[Remark~4.16]{Moe26a}, the homomorphism $\pr_M^{*}$ defined above coincides with the pullback defined in \cite[Definition~4.15]{Moe26a}.
\end{rem}

Following \cite[Definition~4.23]{Moe26a}, we define the Picard group of $(X,M)$.

\begin{dfn}
Let $\PDiv(X)$ be the group of principal divisors on $X$.
\begin{enumerate}
\item We define the \textbf{Picard group} $\Pic(X,M)$ of $(X,M)$ by
\[
\Pic(X,M)=\Div(X,M)/\pr_M^{*}(\PDiv(X)),
\]
where $\pr_M^{*}(\PDiv(X))$ is the image of $\PDiv(X)$ under the composite map
\[
\pr_M^{*}\colon\PDiv(X)\to\Div(X)\to\Div(X,M).
\]
\item Following \cite[Definition~4.35]{Moe26a}, we define the \textbf{effective cone} $\Eff^1(X,M)$ of $(X,M)$ by the cone in $\Pic(X,M)_{\mathbb{R}}$ generated by the classes of effective divisors on $(X,M)$, and we define the \textbf{pseudo-effective cone} $\overline{\Eff}^1(X,M)$ of $(X,M)$ by the topological closure of $\Eff^1(X,M)$ in $\Pic(X,M)_\mathbb{R}$.
\end{enumerate}
\end{dfn}

In our setting, the Picard group and the pseudo-effective cone admit explicit descriptions.

\begin{prop}\label{PicXM}
\leavevmode
\begin{enumerate}
\item We have the canonical isomorphism
\[
\Pic(X,M)\cong\bigoplus_{\alpha\in\mathcal{A}}\mathbb{Z}\tilde{D}_\alpha.
\]
\item We have
\[
\Eff^1(X,M)=\overline{\Eff}^1(X,M)=\bigoplus_{\alpha\in\mathcal{A}}\mathbb{R}_{\geq0}\tilde{D}_\alpha.
\]
\end{enumerate}
\end{prop}

\begin{proof}
\begin{enumerate}
\item We consider the composite map
\[
\PDiv(X)\hookrightarrow\Div(X)\to\Div(X,M)=\Div(G)\oplus\bigoplus_{\alpha\in\mathcal{A}}\mathbb{Z}\tilde{D}_\alpha\to\Div(G),
\]
where the second map is $\pr_M^{*}$ and the last map is the projection.
This composite map is the restriction map from principal divisors on $X$ to divisors on $G$.
Since $\Div(G)=\PDiv(G)$ and $F(X)=F(G)$, this map is surjective.
The injectivity follows from the fact that $\mathcal{O}(G)^\times=F^\times$.
Hence, we have
\[
\PDiv(X)\cong\PDiv(G)=\Div(G).
\]
Moreover, we have the exact sequence
\[
0\to\PDiv(X)\to\Div(G)\oplus\bigoplus_{\alpha\in\mathcal{A}}\mathbb{Z}\tilde{D}_\alpha\to\Pic(X,M)\to0.
\]
It follows that
\[
\Pic(X,M)\cong\bigoplus_{\alpha\in\mathcal{A}}\mathbb{Z}\tilde{D}_\alpha.
\]
\item This follows from the facts that
\[
\overline{\Eff}^1(X)=\bigoplus_{\alpha\in\mathcal{A}}\mathbb{R}_{\geq0}D_\alpha
\]
and that the pullback $\pr_M^{*}\colon\Pic(X)\to\Pic(X,M)$ is given by $D_\alpha\mapsto m_\alpha\tilde{D}_\alpha$.\qedhere
\end{enumerate}
\end{proof}

As in the usual formulation of Manin's conjecture, the formulation of Manin's conjecture for $\mathcal{M}$-points uses the canonical divisor class of a pair to define the $a$-invariant and the $b$-invariant.

\begin{dfn}
Following \cite[Definition~4.33]{Moe26a}, we define the canonical divisor class of $(X,M)$ by
\[
K_{(X,M)}=\pr_M^{*}K_X+\sum_{\alpha\in\mathcal{A}}(m_\alpha-1)\tilde{D}_\alpha.
\]
\end{dfn}

Let
\[
K_X=-\sum_{\alpha\in\mathcal{A}}\rho_\alpha D_\alpha
\]
as in Remark~\ref{Eff}~(2).
Then we have
\begin{align*}
K_{(X,M)}&=\pr_M^{*}K_X+\sum_{\alpha\in\mathcal{A}}(m_\alpha-1)\tilde{D}_\alpha
\\
&=\sum_{\alpha\in\mathcal{A}}(-\rho_\alpha\pr_M^{*}(D_\alpha)+(m_\alpha-1)\tilde{D}_\alpha)
\\
&=\sum_{\alpha\in\mathcal{A}}(-\rho_\alpha m_\alpha+m_\alpha-1)\tilde{D}_\alpha
\end{align*}
in $\Pic(X,M)$.
The $a$-invariant and the $b$-invariant of pairs are defined in \cite[Definition~4.37]{Moe26a}.
In our setting, they are given as follows.

\begin{dfn}
Let $L=\sum_{\alpha\in\mathcal{A}}\lambda_\alpha D_\alpha$ be a big $\mathbb{Q}$-divisor class on $X$.
We define the \textbf{$a$-invariant} of $(X,M)$ with respect to $L$ by
\[
a((X,M),L):=\inf\{t\in\mathbb{R}\mid t\pr_M^{*}L+K_{(X,M)}\in\overline{\Eff}^1(X,M)\}.
\]
We define the \textbf{$b$-invariant} $b(F,(X,M),L)$ of $(X,M)$ with respect to $L$ by the codimension of the minimal supported face of $\overline{\Eff}^1(X,M)$ containing $a((X,M),L)\pr_M^{*}L+K_{(X,M)}$.
\end{dfn}

Recall that in Definition~\ref{ab}, we defined 
\begin{align*}
a((X,D_{\mathbf{m}}),L)&=\max_{\alpha\in\mathcal{A}}\frac{\rho_\alpha-1+m_\alpha^{-1}}{\lambda_\alpha},
\\
\mathcal{A}_{\mathbf{m}}(L)&=\left\{\alpha\in\mathcal{A}\middle|a((X,D_{\mathbf{m}}),L)=\frac{\rho_\alpha-1+m_\alpha^{-1}}{\lambda_\alpha}\right\},
\\
b(F,(X,D_{\mathbf{m}}),L)&=\#\mathcal{A}_{\mathbf{m}}(L).
\end{align*}
The following proposition asserts that the $a$-invariant and the $b$-invariant in Theorem~\ref{main} coincide with the corresponding invariants defined in \cite{Moe26a}.

\begin{prop}\label{abM}
We have
\[
a((X,D_{\mathbf{m}}),L)=a((X,M),L),\quad b(F,(X,D_{\mathbf{m}}),L)=b(F,(X,M),L).
\]
In particular, Theorem~\ref{main} is compatible with Manin's conjecture for $\mathcal{M}$-points formulated in \cite[Conjecture~1.4]{Moe26a} with respect to the $a$-invariant and the $b$-invariant.
\end{prop}

\begin{proof}
We have
\begin{align*}
t\pr_M^{*}L+K_{(X,M)}&=\sum_{\alpha\in\mathcal{A}}(tm_\alpha\lambda_\alpha-\rho_\alpha m_\alpha+m_\alpha-1)\tilde{D}_\alpha
\\
&=\sum_{\alpha\in\mathcal{A}}\lambda_\alpha m_\alpha\left(t-\frac{\rho_\alpha-1+m_\alpha^{-1}}{\lambda_\alpha}\right)\tilde{D}_\alpha.
\end{align*}
Hence, for a real number $t$, the $\mathbb{R}$-divisor class $t\pr_M^{*}L+K_{(X,M)}$ is contained in $\overline{\Eff}^1(X,M)$ if and only if
\[
t\geq\frac{\rho_\alpha-1+m_\alpha^{-1}}{\lambda_\alpha}
\]
for all $\alpha\in\mathcal{A}$.
Hence, we have
\[
a((X,M),L)=\max_{\alpha\in\mathcal{A}}\frac{\rho_\alpha-1+m_\alpha^{-1}}{\lambda_\alpha}=a((X,D_{\mathbf{m}}),L).
\]
The assertion for the $b$-invariant follows from the fact that the minimal supported face of $\overline{\Eff}^1(X,M)$ containing $a\pr_M^{*}L+K_{(X,M)}$ is generated by $\mathbb{R}_{\geq0}\tilde{D}_\alpha$ for $\alpha\in\mathcal{A}\setminus\mathcal{A}_{\mathbf{m}}(L)$.
Hence its codimension is $\#\mathcal{A}_{\mathbf{m}}(L)$.
\end{proof}

\subsection{The leading constant}

In this subsection, we compare the leading constant
\[
c=\lim_{s\to a}(s-a)^b\widehat{\mathsf{H}}_{\mathbf{m}}(\mathbf{0},sL)
\]
in Theorem~\ref{main} with the leading constant predicted in \cite[Conjecture~1.4]{Moe26a}.
We first compute the $\alpha$-constant.
Let $L=\sum_{\alpha\in\mathcal{A}}\lambda_\alpha D_\alpha$ be a big line bundle on $X$, and let
\[
a=a((X,M),L)=a((X,D_{\mathbf{m}}),L),\quad b=b(F,(X,M),L)=b(F,(X,D_{\mathbf{m}}),L)
\]
be the $a$-invariant and the $b$-invariant of $(X,M)$ with respect to $L$, respectively.
As in Theorem~\ref{main}, we assume that $aL+K_X+D_{\mathbf{m}}$ is a rigid divisor.
First, following \cite[Definition~8.1]{Moe26a}, we introduce $M^\circ$ in order to compute the $\alpha$-constant in our setting.

\begin{dfn}
We set
\[
\mathfrak{M}^\circ=\prod_{\alpha\in\mathcal{A}_{\mathbf{m}}(L)}m_\alpha\mathbb{Z}_{\geq0}\times\prod_{\alpha\in\mathcal{A}\setminus\mathcal{A}_{\mathbf{m}}(L)}\{0\}\subseteq\overline{\mathbb{N}}^{\mathcal{A}}
\]
and $M^\circ=((D_\alpha)_{\alpha\in\mathcal{A}},\mathfrak{M}^\circ)$.
\end{dfn}

\begin{rem}
In our setting, the open subvariety $X^\circ$ of $X$ appearing in \cite[Definition~8.1]{Moe26a} coincides with $X$.
Indeed, by \cite[Theorem~1.7]{Moe26a}, the $\mathbb{R}$-divisor class $a\pr_M^{*}L+K_{(X,M)}$ is rigid in the sense of \cite[Definition~4.44]{Moe26a}.
Moreover,
\[
\sum_{\alpha\in\mathcal{A}\setminus\mathcal{A}_{\mathbf{m}}(L)}\lambda_\alpha m_\alpha\left(a-\frac{\rho_\alpha-1+m_\alpha^{-1}}{\lambda_\alpha}\right)\tilde{D}_\alpha
\]
is the unique effective divisor on $(X,M)$ representing the class $a\pr_M^{*}L+K_{(X,M)}$.
The restriction of this effective divisor to $U=G=X\setminus\bigcup_{\alpha\in\mathcal{A}}D_\alpha$ is zero.
It follows that $X=X^\circ$.
\end{rem}

\begin{rem}\label{PicXMcirc}
We define $\Div(X,M^\circ)$, $\pr_{M^\circ}^{*}\colon\Div(X)\to\Div(X,M^\circ)$, $\Pic(X,M^\circ)$, and $\Eff^1(X,M^\circ)$ as in \cite[Definition~4.10, Definition~4.15, Definition~4.23, and Definition~4.35]{Moe26a}, respectively.
For each $\alpha\in\mathcal{A}$, we set
\[
m_\alpha^\circ=\left\{
\begin{array}{ll}
m_\alpha&\text{if }\alpha\in\mathcal{A}_{\mathbf{m}}(L),
\\
\infty&\text{if }\alpha\notin\mathcal{A}_{\mathbf{m}}(L)
\end{array}
\right.
\] 
and we set $\mathbf{m}^\circ=(m_\alpha^{\circ})_{\alpha\in\mathcal{A}}$.
Then the pair $(X,M^\circ)$ corresponds to the Darmon points on $(\mathcal{X},\mathcal{D}_{\mathbf{m}^\circ})$ in the sense of \cite[Definition~3.24]{Moe26a}.
Hence, we have
\[
\Div(X,M^\circ)\cong\Div(G)\oplus\bigoplus_{\alpha\in\mathcal{A}_{\mathbf{m}}(L)}\mathbb{Z}\tilde{D}_\alpha
\]
and the pullback $\pr_{M^\circ}^{*}\colon\Div(X)\to\Div(X,M^\circ)$ is given by
\[
D_G+\sum_{\alpha\in\mathcal{A}}n_\alpha D_\alpha\mapsto D_G+\sum_{\alpha\in\mathcal{A}_{\mathbf{m}}(L)}n_\alpha m_\alpha\tilde{D}_\alpha,
\]
where $D_G\in\Div(G)$ and $(n_\alpha)_{\alpha\in\mathcal{A}}$ is a tuple of integers.
By an argument similar to that in the proof of Proposition~\ref{PicXM}, we obtain
\begin{align*}
\Pic(X,M^\circ)&\cong\bigoplus_{\alpha\in\mathcal{A}_{\mathbf{m}}(L)}\mathbb{Z}\tilde{D}_\alpha,
\\
\Eff^1(X,M^\circ)&=\bigoplus_{\alpha\in\mathcal{A}_{\mathbf{m}}(L)}\mathbb{R}_{\geq0}\tilde{D}_\alpha.
\end{align*}
\end{rem}

We now compute the $\alpha$-constant defined in \cite[Definition~8.2]{Moe26a} in our setting.
By Remark~\ref{PicXMcirc}, we have $\#\Pic(X,M^\circ)_{\mathrm{torsion}}=1$.

\begin{dfn}
We define the \textbf{$\alpha$-constant} of $(X,M)$ with respect to $L$ by
\[
\alpha((X,M),L)=\frac{1}{\#\Pic(X,M^\circ)_{\mathrm{torsion}}}\int_{\Lambda^{\vee}}e^{-\langle\pr_{M^\circ}^{*}(L),\mathbf{x}\rangle}\,\mathrm{d}\mathbf{x}=\int_{\Lambda^{\vee}}e^{-\langle\pr_{M^\circ}^{*}(L),\mathbf{x}\rangle}\,\mathrm{d}\mathbf{x}.
\]
\end{dfn}

In our setting, we obtain
\begin{align*}
\alpha((X,M),L)&=\int_{\mathbb{R}_{\geq0}^{\mathcal{A}_{\mathbf{m}}(L)}}\exp\left(\left\langle-\pr_{M^\circ}^{*}\left(\sum_{\alpha\in\mathcal{A}_{\mathbf{m}}(L)}\lambda_\alpha D_\alpha\right),(x_\alpha)_{\alpha\in\mathcal{A}_{\mathbf{m}}(L)}\right\rangle\right)\,\mathrm{d}\mathbf{x}
\\
&=\int_{\mathbb{R}_{\geq0}^{\mathcal{A}_{\mathbf{m}}(L)}}\exp\left(-\sum_{\alpha\in\mathcal{A}_{\mathbf{m}}(L)}m_\alpha\lambda_\alpha x_\alpha\right)\,\mathrm{d}\mathbf{x}
\\
&=\prod_{\alpha\in\mathcal{A}_{\mathbf{m}}(L)}\int_0^\infty e^{-m_\alpha\lambda_\alpha x_\alpha}\,\mathrm{d}x_\alpha
\\
&=\prod_{\alpha\in\mathcal{A}_{\mathbf{m}}(L)}\frac{1}{m_\alpha\lambda_\alpha}.
\end{align*}
Next, we describe the Tamagawa constant defined in \cite[Definition~8.6]{Moe26a} in our setting.
We recall the Artin $L$-function and the algebraic Brauer group in our setting.
In Notation~\ref{D}~(2), for each $\alpha\in\mathcal{A}$, we have a decomposition into geometrically irreducible components
\[
D_\alpha\otimes_F\overline{F}=\bigcup_{\overline{\alpha}\in\overline{\mathcal{A}}(\alpha)}D_{\overline{\alpha}}.
\]
We set
\begin{align*}
\mathfrak{M}_{\overline{F}}^\circ&=\prod_{\alpha\in\mathcal{A}_{\mathbf{m}}(L)}\prod_{\overline{\alpha}\in\overline{\mathcal{A}}(\alpha)}m_\alpha\mathbb{Z}_{\geq0},
\\
M_{\overline{F}}^\circ&=((D_{\overline{\alpha}})_{\overline{\alpha}\in\overline{\mathcal{A}}(\alpha),\alpha\in\mathcal{A}_{\mathbf{m}}(L)},\mathfrak{M}_{\overline{F}}^\circ).
\end{align*}
By \cite[Remark~4.13]{Moe26a}, we obtain the $\Gal(\overline{F}/F)$-module $\Pic(X_{\overline{F}},M_{\overline{F}}^\circ)$.
Let
\[
L(\Pic(X_{\overline{F}},M_{\overline{F}}^\circ),s)=\prod_{v\in\Omega_F}L_v(\Pic(X_{\overline{F}},M_{\overline{F}}^\circ),s)
\]
be the Artin $L$-function with respect to $\Pic(X_{\overline{F}},M_{\overline{F}}^\circ)$, where $L_v(\Pic(X_{\overline{F}},M_{\overline{F}}^\circ),s)=1$ for $v\in\Omega_F^\infty$.
As in \cite[\S~8.1]{Moe26a}, we set
\[
L^{*}(\Pic(X_{\overline{F}},M_{\overline{F}}^\circ),1)=\lim_{s\to1}(s-1)^bL(\Pic(X_{\overline{F}},M_{\overline{F}}^\circ),s)
\]
and set
\[
\lambda_v=L_v(\Pic(X_{\overline{F}},M_{\overline{F}}^\circ),1)
\]
for each $v\in\Omega_F$.

In our setting, the algebraic Brauer group $\Br_1(X,M^\circ)$ satisfies $\Br_1(X,M^\circ)/\Br F=0$.
Indeed, by definition, we have
\[
\Br_1(X,M^\circ)\subseteq\Br(X,M^\circ)\subseteq\Br(G)=\Br(F).
\]
Hence, the Brauer sum appearing in \cite[Definition~8.3]{Moe26a} reduces to the contribution of the trivial class.

Let $\tau_{X,v}$ denote the Tamagawa measure on $X(F_v)$ defined in \cite[\S~2.1.8]{CLT10} and let $\mathrm{d}\mathbf{x}_v$ be the Haar measure on $G(F_v)$ fixed in \S~2.
As stated in Remark~\ref{Sal}~(1), on $G(F_v)$, we have
\[
\mathrm{d}\mathbf{x}_v=\mathsf{H}_v(\mathbf{x}_v,-K_X)\mathrm{d}\tau_{X,v}.
\]
Following \cite[Definition~8.6]{Moe26a}, we define the Tamagawa constant.

\begin{dfn}
We define the \textbf{Tamagawa constant} by
\[
\tau(F,S,(\mathcal{X},\mathcal{M}),L)=L^{*}(\Pic(X_{\overline{F}},M_{\overline{F}}^\circ),1)\prod_{v\in\Omega_F}\lambda_v^{-1}\int_{G(F_v)}\mathsf{H}_v(\mathbf{x}_v,aL+K_X)^{-1}\delta_{\mathbf{m},v}(\mathbf{x}_v)\,\mathrm{d}\tau_{X,v}.
\]
\end{dfn}

To prove that Theorem~\ref{main} is compatible with the leading constant predicted in \cite[Conjecture~1.4]{Moe26a}, it suffices to show that
\[
\frac{c}{a(b-1)!}=\frac{\alpha((X,M),L)\tau(F,S,(\mathcal{X},\mathcal{M}),L)}{a((X,M),L)(b(F,(X,M),L)-1)!}.
\]
By Proposition~\ref{abM}, this is equivalent to showing that
\[
c=\alpha((X,M),L)\tau(F,S,(\mathcal{X},\mathcal{M}),L),
\]
where
\[
c=\lim_{s\to a}(s-a)^b\widehat{\mathsf{H}}_{\mathbf{m}}(\mathbf{0},sL)=\lim_{s\to a}(s-a)^b\int_{G(\mathbb{A}_F)}\delta_{\mathbf{m}}(\mathbf{x})\mathsf{H}(\mathbf{x},-sL)\,\mathrm{d}\mathbf{x}
\]
is the positive constant defined in Proposition~\ref{c}.
The proof of the following proposition is based on \cite[LEMMA~9.3]{PSTVA21}.

\begin{prop}
With the above notation, we have
\[
c=\alpha((X,M),L)\tau(F,S,(\mathcal{X},\mathcal{M}),L).
\]
In particular, Theorem~\ref{main} is compatible with Manin's conjecture for $\mathcal{M}$-points formulated in \cite[Conjecture~1.4]{Moe26a}.
\end{prop}

\begin{proof}
We recall that the function $\delta_{\mathbf{m}}=\prod_{v\in\Omega_F}\delta_{\mathbf{m},v}$ is defined in Definition~\ref{deltam}.
We denote by $\zeta_{F_\alpha}^*(1)$ the residue of the zeta function $\zeta_{F_\alpha}(s)$ at $s=1$.
Then we have
\begin{align*}
c&=\lim_{s\to a}(s-a)^b\int_{G(\mathbb{A}_F)}\delta_{\mathbf{m}}(\mathbf{x})\mathsf{H}(\mathbf{x},-sL)\,\mathrm{d}\mathbf{x}
\\
&=\lim_{s\to a}(s-a)^b\prod_{v\in\Omega_F}\int_{G(F_v)}\delta_{\mathbf{m},v}(\mathbf{x}_v)\mathsf{H}_v(\mathbf{x}_v,-sL-K_X)\,\mathrm{d}\tau_{X,v}
\\
&=\lim_{s\to a}(s-a)^b\prod_{\alpha\in\mathcal{A}_{\mathbf{m}}(L)}\zeta_{F_\alpha}(m_\alpha(s\lambda_\alpha-\rho_\alpha+1))
\\
&\,\times\lim_{s\to a}\prod_{v\in\Omega_F}\left(\prod_{\alpha\in\mathcal{A}_{\mathbf{m}}(L)}\prod_{w\mid v}\zeta_{(F_\alpha)_w}(m_\alpha(s\lambda_\alpha-\rho_\alpha+1))\right)^{-1}\int_{G(F_v)}\delta_{\mathbf{m},v}(\mathbf{x}_v)\mathsf{H}_v(\mathbf{x}_v,-sL-K_X)\,\mathrm{d}\tau_{X,v}
\\
&=\left(\prod_{\alpha\in\mathcal{A}_{\mathbf{m}}(L)}\frac{1}{m_\alpha\lambda_\alpha}\right)\left(\prod_{\alpha\in\mathcal{A}_{\mathbf{m}}(L)}\zeta_{F_\alpha}^{*}(1)\right)
\\
&\,\times\left(\prod_{v\in\Omega_F}\left(\left(\prod_{\alpha\in\mathcal{A}_{\mathbf{m}}(L)}\prod_{w\mid v}\zeta_{(F_\alpha)_w}(1)^{-1}\right)\int_{G(F_v)}\delta_{\mathbf{m},v}(\mathbf{x}_v)\mathsf{H}_v(\mathbf{x}_v,-aL-K_X)\,\mathrm{d}\tau_{X,v}\right)\right)
\\
&=\alpha((X,M),L)\left(\prod_{\alpha\in\mathcal{A}_{\mathbf{m}}(L)}\zeta_{F_\alpha}^{*}(1)\right)
\\
&\,\times\left(\prod_{v\in\Omega_F}\left(\prod_{\alpha\in\mathcal{A}_{\mathbf{m}}(L)}\prod_{w\mid v}\zeta_{(F_\alpha)_w}(1)^{-1}\right)\int_{G(F_v)}\delta_{\mathbf{m},v}(\mathbf{x}_v)\mathsf{H}_v(\mathbf{x}_v,-aL-K_X)\,\mathrm{d}\tau_{X,v}\right),
\end{align*}
where for each $\alpha\in\mathcal{A}$, $w\mid v$ means that $w$ runs over all places of $F_\alpha$ lying above $v$, and $(F_\alpha)_w$ is the completion of $F_\alpha$ with respect to $w$.
On the other hand, we have
\begin{align*}
&L^{*}(\Pic(X_{\overline{F}},M_{\overline{F}}^\circ),1)^{-1}\prod_{\alpha\in\mathcal{A}_{\mathbf{m}}(L)}\zeta_{F_\alpha}^{*}(1)
\\
&=\left(\lim_{s\to1}(s-1)^bL(\Pic(X_{\overline{F}},M_{\overline{F}}^\circ),s)\right)^{-1}\left(\lim_{s\to1}(s-1)^b\prod_{\alpha\in\mathcal{A}_{\mathbf{m}}(L)}\zeta_{F_\alpha}(s)\right)
\\
&=\lim_{s\to1}L(\Pic(X_{\overline{F}},M_{\overline{F}}^\circ),s)^{-1}\prod_{\alpha\in\mathcal{A}_{\mathbf{m}}(L)}\zeta_{F_\alpha}(s)
\\
&=\lim_{s\to1}\prod_{v\in\Omega_F}\left(L_v(\Pic(X_{\overline{F}},M_{\overline{F}}^\circ),s)^{-1}\prod_{\alpha\in\mathcal{A}_{\mathbf{m}}(L)}\prod_{w\mid v}\zeta_{(F_\alpha)_w}(s)\right)
\\
&=\prod_{v\in\Omega_F}\lambda_v^{-1}\prod_{\alpha\in\mathcal{A}_{\mathbf{m}}(L)}\prod_{w\mid v}\zeta_{(F_\alpha)_w}(1).
\end{align*}
Therefore, we obtain
\[
c=\alpha((X,M),L)L^{*}(\Pic(X_{\overline{F}},M_{\overline{F}}^\circ),1)\prod_{v\in\Omega_F}\lambda_v^{-1}\int_{G(F_v)}\mathsf{H}_v(\mathbf{x}_v,aL+K_X)^{-1}\delta_{\mathbf{m},v}(\mathbf{x}_v)\,\mathrm{d}\tau_{X,v}.\qedhere
\]
\end{proof}

\end{document}